# A COMPUTATIONAL ALGORTHM FOR THE HARDY FUNCTION $Z(t)$, UTLISING SUB-SEQUENCES OF GENERALISED CUBIC GAUSS SUMS, WITH AN OVERALL OPERATIONAL COMPLEXITY $O\left((t/\varepsilon_t)^{[0.25,\,0.3]}\{log(t)\}^{2+o(1)}\right)$ FOR $t \in [10^{23-35}]$

D. M. Lewis†

*Department of Mathematics, University of Liverpool,*

*M&O Building, Peach St, Liverpool L69 7ZL.*

(Tel: *+44 151 794 4014.* Email: *D.M.Lewis@liverpool.ac.uk*)

A. R. Brereton

*University of California, Berkeley, USA*

(Email: *ashbre@ibl.gov*)

**Abstract**

In 2011 G. A. Hiary devised a computational algorithm for the Hardy function $Z(t)$, requiring just $O\left(t^{1/3}\{log(t)\}^{\kappa}\right)$ operations. This compares to $O\left(\sqrt{t}\right)$ operations necessary for computing $Z(t)$ using the classical Riemann-Siegel formula. The methodology involved the sub-division of the Riemann-Siegel formula into sequences of quadratic Gauss/exponential sums of various lengths $N$. Such sums can be computed rapidly, in $O(log(N))$ operations, using standard recursive schemes. More recently, the principal author developed a similar algorithm with an $O\left((t/\varepsilon_t)^{1/3}\{log(t)\}^2\right)$ operational count, accurate to $\varepsilon_t$ in the relative error. Although constructively analogous, the sub-division into quadratic sums was applied to a different asymptotic formula for $Z(t)$, giving the new algorithm an original formulation. This paper presents a significant extension of these ideas. The main theoretical result is an asymptotic expression for $Z(t)$ in terms of sub-sequences of generalised, $m^{th}$-order, Gauss sums of progressively increasing length. Computationally, the main focus falls upon the cubic Gauss sum formulation. The particular parameterisation of these cubic sums makes them amenable to rapid computation, utilising a recursive scheme similar to those implemented for quadratic sums. The net result is a computational algorithm for $Z(t)$ with a reduced $O\left((t/\varepsilon_t)^{[0.25,\,0.3]}\{log(t)\}^{2+o(1)}\right)$ operational count for $t \in [10^{23-35}]$, to high accuracy. Sample computations lend practical support to these findings.

† Corresponding author

## 1. Introduction

Interest in novel methodologies for the computation of the Riemann zeta function $\zeta(s)$ along the critical line $s = 1/2 + \mathrm{i}t$ remains strong, given its fundamental importance to the distribution of the primes. Practically, computations along the critical line come down to the evaluation of Hardy's $Z$-function [12], denoted by $Z(t) \in \mathbb{R}$. This is defined by

$$Z(t) = \zeta(1/2 + \mathrm{i}t)\left[\frac{\Gamma(1/4-\mathrm{i}t/2)}{\Gamma(1/4+\mathrm{i}t/2)}\pi^{\mathrm{i}t}\right]^{-1/2}, \tag{1}$$

or sometimes by the equivalent form $Z(t) = \zeta(1/2 + \mathrm{i}t)e^{\mathrm{i}\theta(t)}$, where

$$\theta(t) = Im\{Log\Gamma(1/4 + \mathrm{i}\,t/2)\} - \frac{t}{2}log(\pi) \approx \frac{t}{2}\left\{log\left(\frac{t}{2\pi}\right) - 1\right\} - \frac{\pi}{8} + \frac{1}{48t} + \cdots \tag{2}$$

is the Riemann-Siegel theta function. Historically such computations for $Z(t)$ were initiated by Riemann himself with his development of the Riemann-Siegel formula [4, 5], a method which remained pretty much state of the art right into the twenty-first century [3] and still retains many practical advantages today. This specifies

$$Z(t) = 2\sum_{N=1}^{N_t}\frac{cos(\theta(t) - tlog(N))}{\sqrt{N}} - (-1)^{N_t}\left(\frac{2\pi}{t}\right)^{1/4}\sum_{r=0}^{M}(-1)^r\left(\frac{2\pi}{t}\right)^{r/2}\Psi_{\mathrm{r}}(p) + R_M(t). \tag{3}$$

Here $N_t = \left\lfloor\sqrt{t/2\pi}\right\rfloor, p = frac\left(\sqrt{t/2\pi}\right)$, $\Psi_0(p) = cos\{2\pi(p^2 - p - 1/16)\}/cos\{2\pi p\}$ and $\Psi_{\mathrm{r}}(p)$ are combinations of the derivatives of the function $\Psi_0(p)$. (Bounds on the cut-off remainder term $R_M(t)$ are discussed in [5], but typically $|R_2(t)| < 0.011t^{-7/4}\ \forall\ t > 200$.) Algorithmically the Riemann-Siegel formula is an $O(\sqrt{t})ops_{s;c}$ operational method (specifically $N_t$ *trigonometric* operations, where the notation $ops_{s;c}$ is used to denote the number of basic $\pm,\times$ operations needed for a single sine/cosine evaluation). It is only in the last decade that the formulation of genuinely faster algorithms in terms of reductions in operational count (as opposed to faster function evaluations brought about by increases in processor speed) have been achieved. Specifically, the work of [2, 7, 8] describes the development of two related algorithms of operational complexity of $O\left(t^{1/3+o(1)}\right)$ in time and $O\left(t^{4/13+o(1)}\right)$ in time and space. These operational counts pertain for a single evaluation of $Z(t)$. Utilising methodology first proposed by [18], the papers [2, 7] also discuss how the $O\left(t^{1/3+o(1)}\right)$ algorithm can be augmented to produce simultaneous multi-evaluations of $Z(t + j\delta)$ at little extra operational expense (here $j\delta \ll t$, typically $j = 1..1000$ and $\delta = 0.04$). More recently, [14] presented the derivation and application of an similar $O\left(t^{1/3+o(1)}\right)$ algorithm for the computation of $Z(t)$. This was developed from an alternative representation of the Hardy function as a sum of Kummer's functions [19, Chap 13], given by

$$Z(t) = \frac{2e^{\mathrm{i}\theta(t)-3\pi t/4-\mathrm{i}\pi/8}\pi^{5/4+\mathrm{i}t/2}}{\Gamma(1/4+\mathrm{i}t/2)(1+e^{-2\pi t})} \sum_{\alpha\in 2\mathbb{N}+1} \alpha\Phi\left(\frac{3}{4}-\frac{\mathrm{i}t}{2},\frac{3}{2},\frac{\pi\alpha^2\mathrm{i}}{4}\right). \tag{4}$$

Here $\Phi(a,b,z) \equiv {}_1F_1(a,b,z)$ denotes Kummer's (or confluent hypergeometric) function. Full details of the derivation of (4), which arises from the definition of $\zeta(s)$ based upon the Riemann-Siegel Integral formula [4, section 7.9; 12, eq. 4.11], are presented in [13].

The algorithms developed by [2, 7] rely primarily on the observation (made by D. R. Heath-Brown [see 7] on analysis presented in [24, section 5.2]) that subsections, or 'blocks', of terms that constitute the main part of the Riemann-Siegel formula could be re-written in the form of a series of quadratic Gaussian sums, albeit with very specific parameters. The algorithm of [14] rests on an analogous observation but applied to asymptotic approximations for the terms of (4). A general quadratic Gaussian sum of length $N$ is defined by

$$S_N(\Phi_1,\Phi_2) = \sum_{k=0}^{N} exp[2\pi\mathrm{i}(\Phi_1 k + \Phi_2 k^2)], \tag{5}$$

where $|\Phi_1| \leq 1/2$ and $|\Phi_2| \leq 1/4$ are fixed coefficients. (Strictly speaking the algorithms of [2, 7, 8] utilise expressions involving quadratic *exponential* sums, which are effectively Gauss sums with amplitudes $(k/N)^j$, but the distinction is academic in what follows.) Typically the Gaussian sums themselves can be evaluated very rapidly, by means of algorithms that utilise the Poisson summation process [11, 12] to recursively reduce the initial sum of length $N = O(t^{1/6})$ down to a short 'kernel' sum of length $\leq K \ll N$, in $n_K = O(log(N/K))$ iterations. The kernel sum is then evaluated exactly and the iterations reversed, to give a final estimate for the original Gauss sum in $O(K + log(N/K))ops_{s;c}$, a vast saving when compared to computing (5) directly. Further details of the implementation of such procedures are presented in the papers of [7, 14] and a generalised version of the methodology will be presented in Section 3 of this work. Since the Hardy function basically requires the computation of $O(t^{1/3})$ such sums, one arrives at the $O(t^{1/3+o(1)})$ operational algorithms of [7, 14].

The obvious question that arises from this body of work is how much further it can be pushed towards the formulation of even faster operational algorithms? Beyond the practical advantages, such advances would have theoretical ramifications. Formulations leading to significant reductions in operational complexity for the Hardy function have the potential to reduce the exponent bounding its growth $Z(t) = O(t^c)$, from its current best estimate of $c = 32/205 + \varepsilon$ [10, 12], found using the exponent pair method [11]. The difficulties with going beyond the computational methods discussed above lie in the fact that the recursive scheme to evaluate the Gauss sums touched on above, relies upon appropriate (approximate) quadratic reciprocity formulae [1, 21, 22]. Essentially such formulae express any, *suitably long*, quadratic Gauss sum in terms of another, *shorter, quadratic* Gauss sum, plus some small terms (typically sums of error function type integrals, which in practice die off very rapidly).

In order to improve upon the $O(t^{1/3+o(1)})$ methods, ideally one would like to substitute cubic, or indeed higher order Gauss sums, in place of the quadratic sums. Theoretically, cubic or higher order sums could be extended to much greater lengths, meaning that the actual numbers appearing in their corresponding formulation of $Z(t)$ would be significantly reduced, resulting in a faster algorithm. The trouble with this idea is that no analogous reciprocity formulae exists for higher order Gauss sums, which renders their rapid computation by a recursive scheme of the type utilised for quadratic sums problematic. However, this initial pessimistic appraisal is somewhat premature. Crucially (see Section 2 below), the specific higher order Gauss sums that appear in the formulation for $Z(t)$ possess the structural characteristic that their higher order coefficients are small, and conform to a regular asymptotic pattern, specifically $\Phi_l = O(t^{(2-l)/2})$ for $l \geq 3$. Together with corresponding constraints on their respective lengths, this feature makes them much more amenable to rapid computation than general sums of this class. Indeed [7] was able to utilise this idea to formulate routines for the rapid computation of certain cubic Gauss sums with small cubic coefficients, $|\Phi_3| < 0.01/N^2$, routines which form the basis for the theoretically faster $(t^{4/13+o(1)})$ algorithm. (Complications arising from its pre-processing requirements [2, 7], means that, in reality, the $O(t^{4/13+O(1)})$ algorithm offers little practical improvement, at least for computationally feasible values of $t$.) However, a feature of these routines is the retention of the reciprocity characteristic inherent in the corresponding quadratic algorithms, i.e. the original cubic Gauss sum is estimated by means of a shorter *but still cubic* Gauss sum. This imposes limitations on the computational savings that can be achieved and more importantly is unnecessarily restrictive. If in the process of its iterative reduction, the original cubic sum (or indeed initial sums of higher orders) can be generalized down to kernel sums of higher orders, then it can be computed much more rapidly. It is the *raison d'être* of this work that the structural characteristics, mentioned above, inherent in the higher order coefficients $\Phi_{3..m}$ of the generalised Gauss sums making up $Z(t)$, are sufficient to allow for their computation in around the same order of operations (specifically see section 3) as those required for the computation of a typical quadratic sum. This in turn would allow for the development of much faster operational algorithms for $Z(t)$ beyond the $O(t^{1/3+o(1)})$ methods already available. The first step in this process is to devise a specific formulation that approximates the Hardy function in terms of generalized Gauss sums (of *any* prescribed order $m \geq 3$), which in turn are suitable for rapid computation. Even without the motivation of potential computational savings, such a formulation is likely to be of some theoretical interest and is the subject of the next section.

## 2. Formulating Hardy's $Z$-Function into Generalized $m^{th}$-order Gauss Sums

This section provides a natural extension to the work of [14, section 2] which demonstrated a structural connection between the Hardy function and a whole series of sub-sequences based upon pairs of generalized *quadratic* Gauss sums $S_{M_t^-}(\Phi_1, \Phi_1)$. Each successive quadratic Gauss sum pair must conform to strict limits on their respective *collection sizes* $M_t \in 2\mathbb{N} + 1$, which in turn fix their respective sum lengths $M_t^- = (M_t - 1)/2$. The aim of this section is to establish an analogous structural connection between the Hardy function and series of generalized $m^{th}$-order Gauss sums $S_{M_t^-}(\Phi_{1..m})$

$$S_{M_t^-}(\Phi_{1..m}) = \sum_{k=0}^{M_t^-} exp[2\pi \mathrm{i}(\Phi_1 k + \Phi_2 k^2 + \Phi_3 k^3 + \cdots + \Phi_m k^m)]. \qquad (6)$$

Increases in the order $m$ permit corresponding increases in the collection size, which means that the number of $S_{M_t^-}(\Phi_{1..m})$ making up $Z(t)$ is very much smaller than in the quadratic case.

The starting point of this work is the exact formulation of $Z(t)$ in terms of a sum of Kummer's functions given by (4) above. When the Kummer's functions are expressed in terms of their respective Euler's integral formulation [16, eq. 4.2(1)] one obtains the following representation [14, eq. 19]

$$Z(t)=\mathcal{H}(t)e^{-\mathrm{i}\pi/8}\left(\frac{\pi}{2^5 t}\right)^{\frac{1}{4}} \sum_{\alpha\in 2\mathbb{N}+1} \alpha \int_0^1 e^{\mathrm{i}\pi\alpha^2 x/4} \frac{\exp[\mathrm{i}(t/2)log\{(1-x)/x\}]}{[x(1-x)]^{1/4}} dx, \qquad (7)$$

where [14, eqs. 10 & 11a]

$$\mathcal{H}(t) = \frac{e^{-[1/32t^2+O(t^{-4})]}e^{2/(12t^2+3)+O(t^{-4})}}{(1+e^{-2\pi t})(1+1/4t^2)^{1/4}e^{[1/24t^2+O(t^{-4})]}} \approx 1 + \frac{1}{32t^2}. \qquad (8)$$

Lengthy analysis of the integrals making up (7) [13, appendix A] in the limit as $t \to \infty$ demonstrates that the infinite series does not formally converge as it stands. However, further analysis shows that the number of summands can be terminated at any odd integer $N_\alpha \geq INT_O(t/\pi) + 2$ and beyond this termination point the series can be continued analytically and estimated using the Euler-Maclaurin summation technique. (The notation $N/INT_{O,\,E}$ will be used to denote the round and floor functions, to the nearest *O*dd and *E*ven integers respectively.) Application of standard contour integration methods to the integrals appearing in (7) [13, appendix A; 14, section 2.1] also demonstrates that all the terms with $\alpha < a = \sqrt{8t/\pi}$ are of $O\left(e^{-O(t)}\right)$. This means that their combined total gives only an exponentially small contribution to $Z(t)$ and in effect acts as an exponential cut-off to the lower end of the series for $\alpha < a$. Hence the really crucial terms which actually determine the value of $Z(t)$ lie within the range $\alpha \in [a, N_\alpha]$, which means that the series (7) can be truncated to a finite sum of the form

$$Z(t)=\mathcal{H}(t)e^{-\mathrm{i}\pi/8}\left(\frac{\pi}{2^5 t}\right)^{\frac{1}{4}} \sum_{\substack{\alpha>a \\ odd}}^{N_\alpha} \alpha \int_0^1 e^{\mathrm{i}\pi\alpha^2 x/4} \frac{\exp[\mathrm{i}(t/2)log\{(1-x)/x\}]}{[x(1-x)]^{1/4}} dx + R_{EM} + T_{a+\varepsilon}. \quad (9)$$

Here $R_{EM}$ represents the Euler-Maclaurin remainder terms [13, eqs. A73-74] and $T_{a+\varepsilon}$ is a transition term (specifically eq. 129 below), which is applicable on those occasions when the parameter $a$ happens to lie within $\pm\varepsilon = O\left(t^{-1/6}\right)$ of an odd integer. In such cases the commencement of the exponential cut-off is anticipated by the appearance of the transition term $T_{a+\varepsilon}$, derived [13, eq. A65] from the corresponding integral in (9). The main sum itself then begins at either $\alpha = INT_O(a) + 4$ if $\varepsilon \geq 0$, or $\alpha = INT_O(a) + 2$ if $\varepsilon < 0$. But the exact details of $R_{EM}$ and $T_{a+\varepsilon}$ have no bearing on the following analysis.

2.1 *Combining together collections of integrals that form the main sum (9) for* $Z(t)$.

The first step in the methodology is divide the terms making up the main sum of (9) into a series of *collections*, each characterised by a *pivot* integer $\alpha_E \in 2\mathbb{N}$ (which lies at the midpoint of the collection) and a collection size $M_t \in 2\mathbb{N}+1$, with $O(M_t) < O(\alpha_E)$. Specifically a single collection consists of those terms with summands $\alpha = \{\alpha_E \pm 1, \alpha_E \pm 3, \ldots, \alpha_E \pm M_t\}$, combined together in pairs, spaced equidistant about the pivot. It is then easy to show that the sum of a collection satisfies

$$\sum_{\substack{j=1\\ odd}}^{M_t} (\alpha_E \pm j) \int_0^1 e^{\mathrm{i}\pi(\alpha_E \pm j)^2 x/4} \frac{\exp[\mathrm{i}(t/2)log\{(1-x)/x\}]}{[x(1-x)]^{1/4}} dx =$$

$$2\sum_{\substack{j=1\\ odd}}^{M_t} \int_0^1 e^{\mathrm{i}\pi\alpha_E^2 x/4} \frac{\exp[\mathrm{i}(t/2)log\{(1-x)/x\}]}{[x(1-x)]^{1/4}} \left\{\alpha_E \cos\left(\frac{\pi j \alpha_E x}{2}\right) + \mathrm{i} j \sin\left(\frac{\pi j \alpha_E x}{2}\right)\right\} e^{\mathrm{i}\pi j^2 x/4} dx. \quad (10)$$

The fact that the main sum commences at $\alpha > a$ rather than $\alpha = 1$ is critical. Since every possible pivot $\alpha_E \geq O(t^{1/2})$, this means that $O(\alpha_E^2) > O(\alpha_E j) > O(j^2)\ \forall\, j \in [1, M_t]$, which in turn implies that (asymptotically) the exponential phase of the integrand in (10) is dominated by the $e^{\mathrm{i}\pi\alpha_E^2 x/4}$ term. This observation provides both a means for the estimation of these integrals and imposes specific limitations on the largest permissible collection size $M_t(\alpha_E)$ for a given pivot integer. Concentrating on the integral containing the $\cos(\pi j \alpha_E x/2)$ term (the analysis is analogous for the $\sin(\pi j \alpha_E x/2)$ term, but since its amplitude $j \ll \alpha_E$ this will not prove significant) one can split the integration range into two equal sub-sections for $x \in [0, 1/2] \cup [1/2, 1]$ and then reflect the second section onto the first using the substitution $y = 1 - x$, to give

$$\int_0^1 e^{\mathrm{i}\pi(\alpha_E^2 + j^2)x/4} \frac{\exp[\mathrm{i}(t/2)log\{(1-x)/x\}]}{[x(1-x)]^{1/4}} \cos\left(\frac{\pi j \alpha_E x}{2}\right) dx = B_C(\alpha_E, t, j) + e^{\mathrm{i}\pi(\alpha_E + j)^2/4}\overline{B_C(\alpha_E, t, j)}. \quad (11)$$

The integral $B_C(\alpha_E, t, j)$, defined by

$$B_C(\alpha_E, t, j) = \int_1^\infty \frac{\exp[\mathrm{i}\pi\alpha_E^2/4(w+1) + \mathrm{i}(t/2)log(w)]}{w^{1/4}(w+1)^{3/2}} \cos\left(\frac{\pi\alpha_E j}{2(w+1)}\right) e^{\mathrm{i}\pi j^2/4(w+1)} dw, \quad (12)$$

is obtained following the substitution $x = (w+1)^{-1}$ for $x \in [0, 1/2]$ [13, appendix A; 14].

An approximation for integral $B_C(\alpha, t, j)$ was derived using standard contour integration methods in [13, 14]. The same basic procedure will be followed here, but extended to a much greater degree of precision, which is necessary in order to formulate $Z(t)$ in terms of generalized Gauss sums. The asymptotically dominant part of the exponential phase possesses a saddle point situated at $w = pc(\alpha_E)$, where $pc$ is the "portcullis" variable defined (for any real $\alpha \geq a$) by

$$pc(\alpha) = 2(\alpha/a)^2\left[1 + (1 - (a/\alpha)^2)^{1/2}\right] - 1 \in [1, \infty),$$

$$\Rightarrow \pi\alpha^2/2 = t(pc + 1)^2/pc \Rightarrow \alpha/a = (pc + 1)/\left(2\sqrt{pc}\right). \tag{13}$$

This "portcullis" variable acts as a mathematical gateway, linking the summands over $\alpha$ in formula (9) to the summands over $N$ in (3), as established in [13] (and see section 4.1 below). It also, as will be established, specifies the particular coefficients $\Phi_{1..m}$ of the Gaussian sums (6) which constitute the Hardy function. The value of $B_C(\alpha, t, j)$ itself is completely dominated by the behaviour of the integrand along the contour defined by $w = pc(\alpha_E) + ue^{i\pi/4}$, with $u \in \left[-(pc+1)/\sqrt{2}, \infty\right]$, passing through the saddle. Along this path

$$B_C(\alpha_E, t, j) = e^{i\pi/4} \int_{\frac{-(pc+1)}{\sqrt{2}}}^{\infty} e^{d_2 u^2} \left[f(u) e^{i\pi j^2/4(pc+1+ue^{i\pi/4})}\right] \cos\left(\frac{\pi\alpha_E j}{2(pc + 1 + ue^{i\pi/4})}\right) du + O\left(e^{-O(t)}\right). \tag{14}$$

Here the (real) coefficient $d_2 = -t(pc - 1)/4(pc)^2(pc + 1) < 0$, and

$$f(u) = \frac{exp\left[\sum_{i=3}^{\infty} d_i(pc)u^i + i\left\{g(u) - \frac{arctan\left(u/\left(pc\sqrt{2} + u\right)\right)}{4} - \frac{3arctan\left(u/\left(pc\sqrt{2} + \sqrt{2} + u\right)\right)}{2}\right\}\right]}{\left[(pc)^2 + u^2 + u\sqrt{2}pc\right]^{1/8}\left[(pc + 1)^2 + u^2 + u\sqrt{2}(pc + 1)\right]^{3/4}},$$

where

$$g(u) = \frac{(\pi\alpha_E^2/4)\left(pc + 1 + u/\sqrt{2}\right)}{\left[(pc + 1)^2 + u^2 + u\sqrt{2}(pc + 1)\right]} + \frac{tlog\left((pc)^2 + u^2 + u\sqrt{2}pc\right)}{4}. \tag{15a, b}$$

Up to this point the methodology is identical to that leading to the connection between the Hardy function and sub-sequences of generalized *quadratic* Gauss sums [14, eq. 54]. This was achieved by first expanding the exponential phases of the integrand in (14) in powers of $u$, truncating the expansion beyond quadratic powers $u^2$ and then formulating an estimate in terms of an asymptotic expansion of a standard integral. In order to establish an analogous connection to generalized $m^{th}$-order Gauss sums, this procedure must be extended to include all terms up to $O(u^m)$, an operation which presents extra technical difficulties.

2.2 *Extended analysis of integral* $B_C(\alpha_E, t, j)$

The first additional consideration concerns the coefficients $d_{i\geq 3}(pc)$ which form the sum $\sum_{i=3}^{\infty} d_i(pc)u^i$ in (15a). In the quadratic case this sum could be ignored, since the bounds on collection size were such as to render its contribution to $B_C(\alpha_E, t, j)$ asymptotically small in comparison with the main estimate of [14, eq. 47]. However, to extend the methodology beyond quadratic case its contribution must be included, up to and including terms of $O(u^m)$, and the form of these coefficients established. The $d_{i\geq 3}(pc)$ represent the (real part) Taylor coefficients of the asymptotically dominant portion of the exponential phase of the integrand in (12), namely

$$\frac{(\pi\alpha_E^2/4)(u/\sqrt{2})}{\left[(pc+1)^2+u^2+u\sqrt{2}(pc+1)\right]}-\frac{t}{2}arctan\left(\frac{u/\sqrt{2}}{pc+u/\sqrt{2}}\right)+\mathrm{i}g(u), \tag{16}$$

expanded around $u = 0$, the saddle point (where $\pi\alpha_E^2/2 = t(pc+1)^2/pc$). Let $h_i(pc) = d_i(pc) + \mathrm{i}g_i(pc)$ represent the combined real and imaginary parts of these coefficients. Then it is relatively easy to show that

$$h_l(pc) = \frac{\mathrm{i}\left[e^{-3\pi\mathrm{i}/4}\right]^l t}{(2l)pc^l(pc+1)^{l-1}}\left[lpc^{l-1}-(pc+1)^{l-1}\right], \qquad l = 1,2,\dots \tag{17}$$

with $h_1 = 0$ (by definition of the saddle) and $h_2 = d_2$ above.

One can rewrite (14) in the form $(I^+ + I^-)/2$, where

$$I^{\pm}(A,B,C) = e^{\mathrm{i}\pi/4}\int_{\frac{-C}{\sqrt{2}}}^{\infty} e^{-Au^2}\left[f(u)e^{\mathrm{i}\pi j^2/4(C+ue^{\mathrm{i}\pi/4})}\right]e^{\pm\mathrm{i}Bj/(C+ue^{\mathrm{i}\pi/4})}du,$$

$$\text{with}\quad A \equiv |d_2| = \frac{t(pc-1)}{4(pc)^2(pc+1)}, \qquad B \equiv \frac{\pi\alpha_E}{2} = \sqrt{\frac{\pi t}{2pc}}(pc+1), \qquad C \equiv pc+1. \qquad (18a,b)$$

In order to estimate the integrals $I^{\pm}(A,B,C)$, it is necessary to expand the integrand in (18a) as a power series in $u$. But clearly such an expansion cannot be valid over the entirety of the semi-infinite integration range. The crucial next step is to make the supposition that the bulk of the integrand (18a) is actually concentrated very close to the saddle point over a much reduced range $\in [-X, X]$, within which such an expansion would be valid. Here $X \in \left[0, C/\sqrt{2}\right)$ is some, as yet unspecified, parameter which must be identified. Furthermore, one needs to demonstrate that beyond this range, when $u \geq |X|$, the remainder of the integral would be negligibly small because the integrand itself would be dying off exponentially.

Suppose, for the moment, that such a suitable parameter $X$ does indeed exist. Since $X < C/\sqrt{2}$, binomial expansions of the phase factors $\left(1 + ue^{\mathrm{i}\pi/4}/C\right)^{-1}$ are valid and the integrand in (18a) can itself be expanded about $u = 0$ to give

$$f(0)e^{\pm \mathrm{i}B^{\pm}j/C}\left[e^{\mp \mathrm{i}B^{\pm}jue^{\mathrm{i}\pi/4}/C^2}\right]\left[e^{-(A\pm B^{\pm}j/C^3)u^2}\right]exp\left[\sum_{l=3}^{\infty}\left(h_l \pm \frac{\mathrm{i}B^{\pm}j\left(-e^{\mathrm{i}\pi/4}\right)^l}{C^{l+1}}\right)u^l\right] + O(u) \tag{19}$$

where $j = 1,3, \ldots, M_t$ and

$$B^{\pm} = B\left[1 \pm \frac{\pi j}{4B}\right] = \frac{\pi\alpha_E}{2}\left[1 \pm \frac{j}{2\alpha_E}\right] = \sqrt{\frac{\pi t}{2pc}}(pc+1) \pm \frac{\pi j}{4}. \tag{20}$$

In equation (19) the $O(u)$ factor represents those terms arising from the expansion of the denominator of the function $f(u)$ given by (15a). Their details are not particularly important since their contribution to integral (18a) is asymptotically smaller [see 14] than that made by the main term of (19). At the endpoints $u = \pm X$, the supposition demands that the real part of the phase should be less than zero. To establish this, first consider the case when $(pc - 1) \geq O(1)$, by far and away the most common scenario, and concentrate, initially, on just the terms of $O(u)$ and $O(u^2)$ contributing to the phase of (19). (Minor modifications to deal with the case when $pc$ approaches unity are discussed towards the end of section 2.4.) Since $O(\alpha_E^2) > O(\alpha_E j)\ \forall\ j \in [1, M_t]$, the quadratic coefficient $A = O(t/pc^2) \gg B^{\pm}j/C^3 = O\left(j\sqrt{t}/pc^{5/2}\right)$. This suggests that in order for the real part of the phase to be negative at $u = \pm X$ one would require quadratic term to exceed the linear term, i.e.

$$AX^2 \geq max\left|\mathrm{Re}\left\{\frac{B^{\pm}jXe^{\mathrm{i}\pi/4}}{C^2}\right\}\right| \approx \left|\frac{BM_tX}{\sqrt{2}C^2}\right| \Rightarrow X \geq \frac{BM_t}{\sqrt{2}AC^2} = 2\sqrt{\frac{\pi}{t}}\frac{pc^{3/2}M_t}{(pc-1)}. \tag{21}$$

Substituting $=$ for $\geq$ in (21) defines $X$ in terms of the (as yet undetermined) collection size $M_t$. Of course the real part of the phase (19) also contains contributions from the higher $O\left(u^{l\geq 3}\right)$ terms, which in theory have the potential to materially alter this estimate. However, in practice that will not be the case. Comparison of the ratio of higher order terms $h_l u^l$ (since $h_l \gg B^{\pm}j/C^{l+1}$) to $Au^2$ at $u = \pm X$ and utilising expressions (17) established for the $h_l$, shows that

$$\left|\frac{h_l u^l}{A u^2}\right|_{u=X} \sim \begin{cases} 2X^{l-2}(pc+1)/[l(pc-1)pc^{l-2}] \text{ if } & pc \lesssim \left(e^{log(l)/(l-1)}-1\right)^{-1} \\ 2X^{l-2}pc/[(pc-1)(pc+1)^{l-2}] \text{ if } & pc \gtrsim \left(e^{log(l)/(l-1)}-1\right)^{-1} \end{cases} \quad l \geq 3 \qquad (22)$$

Examination of the terms on the right-hand side of (22) imposes the condition that $X \leq C/2$ to ensure that they are less than unity (a slightly stricter criterion than originally proposed). So for example, substituting $X = C/2$ into (22) gives

$$\left|\frac{h_l u^l}{A u^2}\right|_{u=X} < \begin{cases} 2[(pc+1)/2pc]^{l-2}(pc+1)/l(pc-1) \text{ if } & pc \lesssim \left(e^{log(l)/(l-1)}-1\right)^{-1} \\ 2[1/2]^{l-2}\, pc/(pc-1) \qquad\qquad \text{ if } & pc \gtrsim \left(e^{log(l)/(l-1)}-1\right)^{-1} \end{cases} \qquad (23)$$

All the terms in (23) are less than unity for all $l \geq 4$, if $(pc-1) \geq O(1)$. In the case $l = 3$, substitution of the exact expression for $h_3$ yields a factor $(2pc^2 - 2pc - 1)/3pc(pc-1) < 1$ too. This analysis demonstrates that within the range $\in [-X,\ X]$ the higher order terms $h_l u^l$ are both decreasing rapidly and are dominated by the quadratic term. Which means that their collective contribution is simply not large enough to materially alter the conclusion that, provided $X$ satisfies constraint specified by (21), the real part of the phase of integrand (18a) will be less than zero when $u \geq |X|$. Hence the corresponding integral is dominated by that part of the integrand concentrated within the range $[-X,\ X]$. Numerical calculations of the approximation to the real phase of (18a) based on just the linear and quadratic terms, compared to the exact value computed $u = \pm X$ (using the collection size established below), confirm these conclusions. The quadratic approximation is almost exact at endpoint $u = -X$ for all $t$, $pc$ and $m$ values, whilst at $u = X$ it contributes about 60% of the total. Crucially the real part of the phase is *always negative* at and beyond these endpoints, as predicted.

Now in order to finally establish $X$, it remains to specify a value for the collection size $M_t(pc, m)$ in (21). This can be done by continuing the expansion in (19) up to powers $u^m$ and then cutting the series off beyond $l = m$, by ensuring the $(m+1)^{th}$ term is negligible at $u = \pm X$. Specifically, one requires

$$|h_{m+1}X^{m+1}| \leq \varepsilon_t, \qquad (24)$$

where $0 < \varepsilon_t \ll 1$ is a small error tolerance parameter, which for the moment will simply assumed to satisfy $\varepsilon_t \underset{t\to\infty}{\longrightarrow} 0$. (In practice, $\varepsilon_t$ fixes the upper bound on the *relative error* that accrues when utilising the recursive algorithm *MGS* devised in section 3 to rapidly compute the associated Gaussian sum. Typically, one would set $\varepsilon_t \leq 0.01$ for a large value of $t$.) Substituting condition (21) into (24) gives

$$\left| h_{m+1} \left\{ \frac{BM_t}{\sqrt{2}AC^2} \right\}^{m+1} \right| \leq \varepsilon_t. \tag{25}$$

The complicated nature of the $h_{m+1}$ coefficients makes this constraint somewhat awkward to formulate and it is easier to approximate them using $A$ and $C$. From (22) one can see that for very large $pc$ values $|h_{m+1}| \approx 2A/C^{m-1}$. However, for somewhat smaller $pc$ values this estimate is too small, since

$$|h_{m+1}| \sim \frac{2A(pc+1)}{(m+1)(pc-1)} \left[ \frac{(pc+1)}{pc} \right]^{m-1} \frac{1}{C^{m-1}} < \frac{2^{m-1}A}{C^{m-1}}, \tag{26}$$

for all $(pc-1) \geq O(1)$ and $m \geq 3$. (See section 2.4 for the limiting case $(pc-1) \to 0$.) Substitution of (26) into (25) gives

$$\left| h_{m+1} \left\{ \frac{BM_t}{\sqrt{2}AC^2} \right\}^{m+1} \right| \leq \left| \frac{2^{m-1}A}{C^{m-1}} \left\{ \frac{BM_t}{\sqrt{2}AC^2} \right\}^{m+1} \right| = \left| \frac{2^{(m-3)/2}B^{m+1}M_t^{m+1}}{A^m C^{3m+1}} \right| \leq \varepsilon_t$$

$$\Rightarrow M_t(pc,m) \leq \left[ \frac{A^m C^{3m+1} \varepsilon_t}{2^{(m-3)/2} B^{m+1}} \right]^{1/(m+1)} = \frac{1}{\sqrt{\pi}} \left[ \frac{\varepsilon_t (pc^2-1)^m t^{(m-1)/2}}{2^{2(m-1)} pc^{(3m-1)/2}} \right]^{1/(m+1)}, \tag{27}$$

as an upper bound for the collection size. This bound gives a natural extension to the bound derived for the quadratic ($m = 2$) case [14, eq. 40], except for an extra $\sqrt{2}$ factor in the denominator. Substituting (27) into (21) determines a suitable definition for $X$

$$X = \left[ \frac{\varepsilon_t pc^2 (pc+1)^m}{t 2^{m-3}(pc-1)} \right]^{1/(m+1)}, \quad \Rightarrow \quad \frac{X}{C} = \left[ \frac{\varepsilon_t pc^2}{t 2^{m-3}(pc^2-1)} \right]^{1/(m+1)} < \frac{1}{2}, \tag{28}$$

provided $t \gg 16\varepsilon_t$ for all possible $pc$ values (since $pc(\alpha_E) \geq 1 + O(t^{-1/4})$ for all possible pivot values $\alpha_E$). Hence the values for $M_t(pc,m)$ and $X$ specified by (27) and (28) ensure that, a) integrand (18a) is negligible outside the range $u \in [-X, X]$ as per the supposition and b) inside this range the integrand can be approximated by the expansion (19), truncated at $l = m$, to an accuracy of at least $O(\varepsilon_t)$ in the relative error. Substituting (19) into (18a) and truncating the expression at $l = m$ gives

$$I^{\pm}(A, B^{\pm}, C) = e^{i\pi/4} f(0) e^{\pm iB^{\pm}j/C} \int_{-X}^{X} exp[w^{\pm}(u)]\, du + O(\varepsilon_t) \tag{29}$$

where

$$w^{\pm}(u) = qu - pu^2 + \sum_{l=3}^{m} s_l u^l,$$

$$q = \mp \frac{\mathrm{i}je^{\mathrm{i}\pi/4}B^{\pm}}{C^2}, \quad p = A \pm \frac{jB^{\pm}}{C^3} > 0, \quad s_l = h_l \pm \frac{\mathrm{i}B^{\pm}j\left(-e^{\mathrm{i}\pi/4}\right)^l}{C^{l+1}}. \qquad (30a, b)$$

Next the constituent integral in (29) must be estimated. To do this, first generalise $u \in \mathbb{C}$ and consider a disc $D$ centred on $u = 0$ of radius $X$ which encompasses the integration range of (29), as shown in Fig. 1. The polynomial $w(u)$ defined by (30) possesses a saddle point situated at $u = u_*$ such that $w'(u_*) = 0$, where

$$u_* = \frac{q}{2p} + \frac{3s_3q^2}{8p^3} + \frac{9s_3^2q^3}{16p^5} + \ldots = \frac{q}{2p}\left(1 + O\left(\frac{s_3q}{p^2}\right)\right) \le \frac{q}{2p}\left(1 + O\left(\left\{\frac{\varepsilon_t pc^2}{t(pc^2-1)}\right\}^{\frac{1}{4}}\right)\right). \qquad (31)$$

The estimate for $s_3q/p^2$ is obtained by setting $j = M_t$ in (29) and applying definitions (18b) and (27) alongside approximation (26). From (31) one can see that the first term completely dominates the remainder. Substituting (30b) into (31) and using (21) establishes that this saddle point is also an element of the disc, since $u_* \approx q/2p \approx \mp e^{\mathrm{i}3\pi/4}X/\sqrt{2} \in D$. This observation suggests that the best way to estimate integral (29) is deform the contour of integration from the real line $u \in [-X, X]$ by rotating it, clockwise through and an angle $\pi/4$ radians, so that it now passes through $u_*$, as shown in Fig. 1 This operation has the effect of creating two additional arcs along the edge of $D$ from $-X$ to $-Xe^{-\mathrm{i}\pi/4}$ and $Xe^{-\mathrm{i}\pi/4}$ to $X$, but they have no material effect since the integrand in (29) is exponentially small along the entirety of both these additional sections. Hence integral (29) will be dominated by its value along the contour joining $u = \pm Xe^{-\mathrm{i}\pi/4}$ traversing through $u_*$. Along this contour the polynomial (30a) can be expanded around the saddle to give

$$w^{\pm}(u) = w^{\pm}(u_*) - \left[p - \sum_{l=3}^{m}\binom{l}{2}s_l u_*^{l-2}\right](u-u_*)^2 + \sum_{k=3}^{m}\left[\sum_{l=k}^{m}\binom{l}{k}s_l u_*^{l-k}\right](u-u_*)^k. \quad (32)$$

Here the results already established ensure that the quadratic coefficient in (32) is completely dominated by the factor $p = O(t)$, since

$$\sum_{l=3}^{m}\binom{l}{2}\left(\frac{s_l}{p}\right)u_*^{l-2} < \sum_{l=3}^{m}\binom{l}{2}\left(\frac{h_l}{A}\right)\left(\frac{X}{\sqrt{2}}\right)^{l-2} < 2\sum_{l=3}^{m}\binom{l}{2}\left(\frac{\varepsilon_t}{t(pc+1)2^{m-2}}\right)^{(l-2)/(m+1)} \ll 1. \quad (33)$$

In addition, the proceeding analysis leading up to estimate (22) means that the higher order terms in (32) are very unlikely to contribute significantly to integral (29), since they are dominated at the endpoints, when the integrand is exponentially small in any case. If this assumption is correct (as justified below) and the higher order terms can be ignored, then the substitution $z = \sqrt{p}(u - u_*)$ transforms integral (29) to

$$
\begin{aligned}
I^{\pm}(A, B^{\pm}, C) &\approx \frac{f(0)e^{\pm \mathrm{i}B^{\pm}j/C + w^{\pm}(u_*)}}{e^{-\mathrm{i}\pi/4}\sqrt{p}} \int_{-\sqrt{p}\{Xe^{-\mathrm{i}\pi/4}+u_*\}}^{\sqrt{p}\{Xe^{-\mathrm{i}\pi/4}-u_*\}} e^{-z^2}\, dz + O(\varepsilon_t) \\
&= \frac{f(0)e^{\pm \mathrm{i}B^{\pm}j/C + w^{\pm}(u_*)}}{2e^{-\mathrm{i}\pi/4}\sqrt{p/\pi}} \left[\operatorname{erfc}\left\{-\sqrt{p}\{Xe^{-\mathrm{i}\pi/4} + u_*\}\right\} - \operatorname{erfc}\left\{\sqrt{p}\{Xe^{-\mathrm{i}\pi/4} - u_*\}\right\}\right] + O(\varepsilon_t) \\
&\approx \frac{f(0)e^{\pm \mathrm{i}B^{\pm}j/C + w^{\pm}(u_*)}}{2e^{-\mathrm{i}\pi/4}\sqrt{p/\pi}} \left[\operatorname{erfc}\left\{\left\{-\sqrt{p}Xe^{-\mathrm{i}\pi/4} - \frac{q}{2\sqrt{p}}\right\}\right\} - \operatorname{erfc}\left\{\left\{\sqrt{p}Xe^{-\mathrm{i}\pi/4} - \frac{q}{2\sqrt{p}}\right\}\right\}\right] + O(\varepsilon_t)
\end{aligned}
\tag{34}
$$

utilising (31) for $u_*$. Here $\operatorname{erfc}(z)$ represents the complementary error function [19, Chap. 7]. Analysis of the $z$ coefficients in the $\operatorname{erfc}(z)$ terms shows that

$$
\left|\frac{q}{2\sqrt{p}}\right| \approx \left|\frac{q}{2\sqrt{A}}\right| \lesssim \left|\frac{M_t B}{2\sqrt{A}C^2}\right| < \left|\frac{M_t B}{\sqrt{2A}C^2}\right| = \left|\sqrt{A}X\right| \approx \left|\sqrt{p}X\right| \approx \varepsilon_t^{1/(m+1)} \left\{\frac{t(pc^2 - 1)}{16pc^2}\right\}^{\frac{m-1}{2(m+1)}}. \tag{35}
$$

This means the values of $\operatorname{erfc}(z)$ are determined predominantly by the $\sqrt{p}Xe^{-\mathrm{i}\pi/4}$ factors and these factors are asymptotically large (q.v. eqs. 18a & 21). Hence

$$
\operatorname{erfc}\left\{\left\{-\sqrt{p}Xe^{-\mathrm{i}\pi/4} - \frac{q}{2\sqrt{p}}\right\}\right\} - \operatorname{erfc}\left\{\left\{\sqrt{p}Xe^{-\mathrm{i}\pi/4} - \frac{q}{2\sqrt{p}}\right\}\right\} = 2 + O\left(\frac{1}{\left|\sqrt{p}X\right|}\right) = 2 + O\left(t^{-\frac{(m-1)}{2(m+1)}}\right). \tag{36}
$$

Substituting this result into (34) gives the desired estimate

$$
I^{\pm}(A, B^{\pm}, C) = \frac{e^{\mathrm{i}\pi/4} f(0) e^{\pm \mathrm{i}B^{\pm}j/C + w^{\pm}(u_*)}\sqrt{\pi}}{\sqrt{p}} + O(\varepsilon_t) = \frac{e^{\mathrm{i}\pi/4} f(0) e^{\pm \mathrm{i}B^{\pm}j/C + w^{\pm}(u_*)}\sqrt{\pi}}{\sqrt{A \pm jB^{\pm}/C^3}} + O(\varepsilon_t)
$$

$$
\Rightarrow \quad I^{\pm}(A, B^{\pm}, C) = \frac{e^{\mathrm{i}\pi/4} f(0) e^{\pm \mathrm{i}B^{\pm}j/C + w^{\pm}(u_*)}\sqrt{\pi}}{\sqrt{A}} + O(\varepsilon_t). \tag{37}
$$

The final simplification is justified since $A = O(t/pc^2) \gg B^{\pm}j/C^3 = O\big(j\sqrt{t}/pc^{5/2}\big)$ as noted earlier. It remains to estimate the contribution of the higher order terms and justify their exclusion from (32) when formulating (34). The largest of these excluded terms are obviously those of third order. Following their inclusion, integral (34) would take the generic form

$$\int_{-\sqrt{p}\{Xe^{-i\pi/4}+u_*\}}^{\sqrt{p}\{Xe^{-i\pi/4}-u_*\}} e^{-z^2+s_3 z^3/p^{3/2}}dz \sim \sqrt{\frac{\pi}{p}}\sum_{k=0}^{\infty}\frac{(6k-1)!!}{(2k)!\,2^{3k}}\left(\frac{s_3^2}{p^3}\right)^k. \tag{38}$$

Since the factor $|s_3^2/p^3| \approx |h_3^2/A^3| = O(t^{-1})$ is so small, the inclusion of the third and indeed higher order terms in (32) would make no significant difference to estimate (37).

Together, these results establish the desired goal of this subsection, namely an estimate for the integral $B_C(\alpha_E, t, j) = (I^+ + I^-)/2$ itself. Substituting (37) into (14) and applying (15, 18 & 20), gives

$$B_C(\alpha_E,t,j) = \frac{e^{i\pi/4}}{(pc+1)}\sqrt{\frac{\pi pc^{3/2}}{t(pc-1)}}\,e^{\frac{i\pi(\alpha_E^2+j^2)}{4(pc+1)}+\frac{it}{2}log(pc)}\left[e^{ij\sqrt{\pi t/2pc}+w^+(u_*)} + e^{-ij\sqrt{\pi t/2pc}+w^-(u_*)}\right] + O(\varepsilon_t). \tag{39}$$

This approximation suffices to establish a formulation for $Z(t)$ in terms of generalized Gaussian sums. The final piece of analysis centres on the $w^{\pm}(u)$ polynomials (30), evaluated at their respective saddle points $u = u_*$, given by (31).

2.3 *The polynomial expressions* $w^{\pm}(u_*)$.

The calculation of these terms relies on a relatively simple recursive scheme, which in practice becomes algebraically very complicated. However, this turns out to be far less of a problem than one might imagine since, as will be seen, the crucial $\Phi_{1..m}$ coefficients governing the Gauss sums "crystalize" out into relatively simple and predictable pattern. The basic procedure relies on computing successively more precise values for the saddle point at $u = u_*$, up to the specified order $2,3,4,\ldots.,m$ of Gauss sums (6) one wishes to utilize. The relatively simple approximation (31) would be sufficient to express $Z(t)$ in terms of generalized $m = 3^{rd}$ order Gauss sums. To go beyond this specified higher order, the following algorithm can be employed.

*Set* $n = 4$, $u_{*,4} = u_{*,3} + \epsilon_4 = \frac{q}{2p} + \frac{3s_3q^2}{8p^3} + \frac{9s_3^2q^3}{16p^5} + \epsilon_4$ *and compute*

$$w'(u_{*,n}) = q - 2pu_{*,n} + \sum_{l=3}^{n} l s_l u_{*,n}^{l-1}. \tag{40}$$

*Establish a linear approximation to* $w'(u_{*,n}) = w_0 + \epsilon_n w_1 + O(\epsilon_n^2)$ *and set* $\epsilon_n = -w_0/w_1$.

*Then compute the first two terms of the Taylor series for* $\epsilon_n$ *about* $q = 0$, *typically giving*

$$\epsilon_n \approx \frac{n s_n q^{n-1}}{2^n p^n} + a\ term\ in\ q^n\ . \tag{41}$$

*Substitute (41) into* $u_{*,n}$ *above and compute*

$$w(u_{*,n}) = q u_{*,n} - p u_{*,n}^2 + \sum_{l=3}^{n} s_l u_{*,n}^l \tag{42}$$

*up to powers of* $q^n$*. This gives the next order approximation to* $w(u_*)$*. If* $n < m$ *then increment* $n\ by\ 1\ to\ n+1$*, set* $u_{*,n} = u_{*,n-1} + \epsilon_n$ *and repeat steps (40-42). If* $n = m$ *stop.*

Repeated mechanically, this procedure yields for example

$$u_{*,4} = u_{*,3} + \frac{s_4 q^3}{4p^4} + \frac{(135 s_3^3 + 120 s_3 s_4 p) q^4}{128 p^7}, \tag{43a}$$

.

$$u_{*,5} = u_{*,4} + \frac{5 s_5 q^4}{32 p^5} + \frac{3(189 s_3^4 + 252 s_3^2 s_4 p + 4(15 s_3 s_5 + 8 s_4^2) p^2) q^5}{256 p^9}. \tag{43b}$$

$$w^{\pm}(u_{*,3}) = \frac{q^2}{4p} + \frac{s_3 q^3}{8p^3}, \qquad w^{\pm}(u_{*,4}) = w^{\pm}(u_{*,3}) + \frac{(9 s_3^2 + 4 s_4 p) q^4}{64 p^5}, \tag{44a, b}$$

$$w^{\pm}(u_{*,5}) = w^{\pm}(u_{*,4}) + \frac{(27 s_3^3 + 24 s_4 s_3 p + 4 s_5 p^2) q^5}{128 p^7}. \tag{44c}$$

The correction terms decline since $\epsilon_n = O(s_n q/\{s_{n-1} p\}) \epsilon_{n-1} = O\left(t^{-1/(m+1)}\right) \epsilon_{n-1}$, meaning that each successive estimate of $u_{*,n}$ and $w^{\pm}(u_{*,n})$ becomes more and more precise. (The appropriate $\pm$ terms are expressed in the parameter definitions specified in eq. 30.)

So for example, suppose one wanted to express $Z(t)$ in terms of generalized $m = 3^{rd}$ order Gauss sums. Then one would substitute (44a) into (39), employ expressions (30) and simplify the resulting exponential phases up to powers of $j^3$. These routine operations yield

$$B_C(\alpha_E, t, j) = 2\sqrt{\frac{\pi p c^{3/2}}{t(pc-1)}} \frac{e^{\frac{i\pi}{4}\left\{\frac{\alpha_E^2}{(pc+1)}+1\right\}+\frac{it}{2}log(pc)}}{(pc+1)} e^{-\frac{i\pi j^2}{4(pc-1)}} cos\left[j\sqrt{\frac{\pi t}{2pc}} + \frac{j^3}{3\sqrt{t}(pc-1)^3}\left(\frac{\pi pc}{2}\right)^{\frac{3}{2}}\right] [1 + O(\varepsilon_t)]. \tag{45}$$

The quadratic equivalent of this result was derived previously [14, eq. 47]. In that case since the collection size (27) was specified to be $M_t(pc, 2) \sim \varepsilon_t^{1/3} t^{1/6}$, which meant that the $O(j^3)$ term in (45) was unnecessary. For the larger collection size $M_t(pc,3) \sim \varepsilon_t^{1/4} t^{1/4}$ associated with $m = 3$, it will make a significant contribution. The equivalent expression for $Z(t)$ then follows (almost) automatically (see next subsection).

However, there is no need to limit the expansion (45) to just cubic terms. Operations (40-42) can be repeated up to finite value of $m$ and the results utilised to obtain an expansion for $B_C(\alpha_E, t, j)$ up to powers of $O(j^m)$. Given the increasing algebraic complexity of these operations one might expect this to be reflected in the final results by a sequence of pseudo-random coefficients. But instead an intricate, but concise, pattern emerges. Replicating operations (40-42) up to $m = 9$, yields the following

$$B_C(\alpha_E,t,j) = 2\sqrt{\frac{\pi pc^{\frac{3}{2}}}{t(pc-1)}} \frac{e^{\frac{i\pi}{4}\left\{\frac{\alpha_E^2}{(pc+1)}+1\right\}+\frac{it}{2}log(pc)}}{(pc+1)} exp\left[-\frac{i\pi j^2}{4(pc-1)} - it\sum_{\substack{n=4\\ even}}^{m} \frac{L_{n-3}(pc)}{n(pc-1)^{2n-3}}\left(j\sqrt{\frac{\pi pc}{2t}}\right)^n\right] \times$$

$$cos\left[j\sqrt{\frac{\pi t}{2pc}} + t\sum_{\substack{n=3\\ odd}}^{m} \frac{L_{n-3}(pc)}{n(pc-1)^{2n-3}}\left(j\sqrt{\frac{\pi pc}{2t}}\right)^n\right][1+O(\varepsilon_t)]. \tag{46}$$

Here $L_n(pc)$ represents a series of polynomials in the portcullis variable of $O(pc^n)$. Explicit calculation of the first seven of these yields

$$L_0(pc) = 1,$$
$$L_1(pc) = pc + 1,$$
$$L_2(pc) = pc^2 + 3pc + 1,$$
$$L_3(pc) = pc^3 + 6pc^2 + 6pc + 1,$$
$$L_4(pc) = pc^4 + 10pc^3 + 20pc^2 + 10pc + 1,$$
$$L_5(pc) = pc^5 + 15pc^4 + 50pc^3 + 50pc^2 + 15pc + 1,$$
$$L_6(pc) = pc^6 + 21pc^5 + 105pc^4 + 175pc^3 + 105pc^2 + 21pc + 1. \tag{47}$$

The pattern of coefficients making up these polynomials echoes the structure of Pascal's triangle and one would expect them to satisfy some simple formulism. A little trial and error readily leads to the following result

$$L_n(pc) = \sum_{l=0}^{n} \frac{(n+1)}{(l+1)(n+1-l)} \binom{n}{l}^2 pc^{n-l}. \tag{48}$$

One would anticipate that repeated application of operations (40-42) extends the rather striking pattern exhibited in equations (46-48) for all positive $m$. (Numerical evidence suggests that the pattern does indeed continue far beyond $m = 9$, although it is unclear how one would actually *prove* this $\forall\ m$.)

The pattern of terms exhibited in (46) and in particular the appearance of reciprocal factors of $(pc-1)^{2m-3}$ in the coefficients, also highlights a significant issue, namely what happens in the regime when $(pc-1) < O(1)$. One would expect the size of leading term *excluded* (the $(m+1)^{th}$ term) from the phase of (46) to be of $O(\varepsilon_t)$, even at its maximum value when $j = M_t$. And indeed this is the case for $(pc-1) \geq O(1)$ since

$$\frac{L_{m-2}(pc)}{(m+1)(pc-1)^{2m-1}}\left(M_t\sqrt{\frac{\pi pc}{2t}}\right)^{m+1} \sim \frac{\varepsilon_t (pc+1)^m}{(m+1)pc(pc-1)^{m-1}2^{(5m-3)/2}} = O(\varepsilon_t), \qquad (49)$$

utilising the collection size $M_t$ specified by (27). But this collection size is inappropriate for use in the regime $(pc-1) < O(1)$, since (49) would be $O(\varepsilon_t L_{m-2}(1)/(pc-1)^{m-1}) \neq O(\varepsilon_t)$. This observation indicates that collection size (27) is simply too large for this (relatively small) regime and must be suitably reduced. An analogous problem arose in the *quadratic* Gauss sum formulation and a satisfactory modification to the collection size was derived in [14, section 4.3.2]. A similar modification to the collections size, suitable for the proposed generalized Gauss sum formulation, is derived in the next subsection.

2.4 *$Z(t)$ composed in terms of Generalized $\mathrm{m}^{th}$-order Gauss Sums*

The definition of $M_t(t,m)$ given by (27), implicitly assumes that $M_t(t,m) \geq 1$. Now if $\alpha_E \approx a$ this may not happen. If for instance $\alpha_E = a + \kappa \Longrightarrow (pc-1) \approx 2\sqrt{\kappa}(2\pi/t)^{1/4}$ and depending upon the value of $\varepsilon_t$, the constraint that $M_t(t,m) \geq 1$ is not guaranteed. (Although as mentioned above $M_t(t,m)$ is inappropriate for the regime $(pc-1) < O(1)$, this consideration will also apply to the modified collection size discussed below.) All this means is that there exists a constant, denoted by $b \geq O(1)$, such that if $\alpha \in (a, a+b)$ then the terms of (9) can only be added together one by one, rather than as amalgamative collections of length $M_t(t,m)$. In practice if $\alpha \in (a, a+b)$ then the value $j = 0$ must be substituted into equations (46) and (12) in what follows.

With this caveat, substitution of approximation (46) into (11) yields

$$\int_0^1 \frac{e^{\frac{\mathrm{i}\pi(\alpha_E^2+j^2)x}{4}+\mathrm{i}\left(\frac{t}{2}\right)log\left\{\frac{(1-x)}{x}\right\}}}{[x(1-x)]^{1/4}}\cos\left(\frac{\pi j\alpha_E x}{2}\right)dx = \sqrt{\frac{\pi pc^{\frac{3}{2}}}{t(pc-1)}}\frac{4e^{\mathrm{i}\pi/8}}{(pc+1)}\cos\left[j\sqrt{\frac{\pi t}{2pc}} + t\sum_{\substack{n=3\\ odd}}^{m}\frac{L_{n-3}(pc)}{n(pc-1)^{2n-3}}\left(j\sqrt{\frac{\pi pc}{2t}}\right)^n\right]$$

$$\times\cos\left[\frac{\pi\alpha_E^2}{4(pc+1)} + \frac{t}{2}log(pc) + \frac{\pi}{8} - \frac{\pi j^2}{4(pc-1)} - t\sum_{\substack{m=4\\ even}}^{m}\frac{L_{m-3}(pc)}{m(pc-1)^{2m-3}}\left(j\sqrt{\frac{\pi pc}{2t}}\right)^m\right][1+O(\varepsilon_t)].$$

(50)

(N.B. in equation 11 the term $e^{i\pi(\alpha_E+j)^2/4} \equiv e^{i\pi/4}$ as $\alpha_E$ is even and $j$ is odd.) Since every possible pivot $\alpha_E \geq O(t^{1/2})$, this means that $O(\alpha_E) > O(j)\ \forall\, j \in [1, M_t]$, and the amplitude of the $j\sin(\pi j \alpha_E x/2)$ term in (10) is much smaller $\ll O(\varepsilon_t)$ than the corresponding $\alpha_E\cos(\pi j\alpha_E x/2)$ term $\forall\, j \in [1, M_t]$ . Dropping this subdominant term from (10), applying approximation (50) and substituting the result into (9) gives

$$Z(t) = T_{a+\varepsilon} + \mathcal{H}(t)2\sqrt{2}\left\{\sum_{\substack{\alpha>a\\ odd}}^{a+b} \frac{\cos\left(\frac{t}{2}\left\{log(pc) + \frac{1}{pc} + 1\right\} + \frac{\pi}{8}\right)}{(\alpha^2 - a^2)^{\frac{1}{4}}}[1 + O(\Lambda_\alpha^{-1})] + R_{EM}\right\} +$$

$$\sum_{\substack{\alpha_E>a+b\\ step\ 2M_t+2}}^{N_\alpha} \frac{\mathcal{H}(t)4\sqrt{2}}{(\alpha_E^2 - a^2)^{\frac{1}{4}}} \sum_{\substack{j=1\\ odd}}^{M_t} \left\{\cos\left[\frac{t}{2}\left(log(pc) + \frac{1}{pc} + 1\right) + \frac{\pi}{8} - \frac{\pi j^2}{4(pc-1)} - t\sum_{\substack{n=4\\ even}}^{m} \frac{L_{n-3}(pc)}{n(pc-1)^{2n-3}}\left(j\sqrt{\frac{\pi pc}{2t}}\right)^n\right] \times \right.$$

$$\left. \cos\left[j\sqrt{\frac{\pi t}{2pc}} + t\sum_{\substack{n=3\\ odd}}^{m} \frac{L_{n-3}(pc)}{n(pc-1)^{2n-3}}\left(j\sqrt{\frac{\pi pc}{2t}}\right)^n\right]\right\}[1 + O(\varepsilon_t)]. \tag{51}$$

Here the phase of (51) has been simplified using (13) and similarly the amplitude simplified using the identity

$$\left(\frac{32\pi pc}{t(pc-1)^2}\right)^{1/4} = \frac{2\sqrt{2}}{(\alpha_E^2 - a^2)^{1/4}}. \tag{52}$$

The sum over the collections $j = 1,3, \ldots, M_t$ in (51) can be rewritten in the form

$$\frac{1}{2}\mathrm{Re}\left\{e^{-iX}\sum_{\substack{j=1\\ odd}}^{M_t} e^{i\Psi(j)}\left(e^{i\Omega(j)} + e^{-i\Omega(j)}\right)\right\}, \tag{53}$$

where

$$X = \frac{t}{2}\left(log(pc) + \frac{1}{pc} + 1\right) + \frac{\pi}{8}, \qquad \Psi(j) = \frac{\pi j^2}{4(pc-1)} + t\sum_{\substack{n=4\\ even}}^{m} \frac{L_{n-3}(pc)}{n(pc-1)^{2n-3}}\left(j\sqrt{\frac{\pi pc}{2t}}\right)^n,$$

$$\Omega(j) = j\sqrt{\frac{\pi t}{2pc}} + t\sum_{\substack{n=3\\ odd}}^{m} \frac{L_{n-3}(pc)}{n(pc-1)^{2n-3}}\left(j\sqrt{\frac{\pi pc}{2t}}\right)^n. \tag{54}$$

Substituting $(2k + 1)$ for $j$ into the above, with $k = 0, 1, 2, \ldots, (M_t - 1)/2 = M_t^-$, transforms (53) into

$$\frac{1}{2}\mathrm{Re}\{e^{\mathrm{i}\omega^+(pc)}S_{M_t^-}(\Phi_{1..m}^+) + e^{\mathrm{i}\omega^-(pc)}S_{M_t^-}(\Phi_{1..m}^-)\}, \tag{55}$$

where

$$\omega^\pm(pc) = -\frac{t}{2}\left(log(pc) + \frac{1}{pc} + 1\right) - \frac{\pi}{8} + \frac{\pi}{4}\left[\frac{1}{(pc-1)} \pm \frac{a}{\sqrt{pc}}\right] + \frac{\pi}{2}\sum_{n=3}^{m}\left(\pm\frac{2}{a}\right)^{n-2}\frac{pc^{n/2}L_{n-3}(pc)}{n(pc-1)^{2n-3}},$$

(56a)

$$\Phi_1^\pm(pc) = \frac{1}{2(pc-1)} \pm \frac{a}{4\sqrt{pc}} + \frac{1}{2}\sum_{n=3}^{m}\left(\pm\frac{2}{a}\right)^{n-2}\frac{pc^{n/2}L_{n-3}(pc)}{(pc-1)^{2n-3}}, \tag{56b}$$

$$\Phi_2^\pm(pc) = \frac{1}{2(pc-1)} + \frac{1}{2}\sum_{n=3}^{m}\binom{n-1}{1}\left(\pm\frac{2}{a}\right)^{n-2}\frac{pc^{n/2}L_{n-3}(pc)}{(pc-1)^{2n-3}}, \tag{56c}$$

$$\Phi_l^\pm(pc) = \frac{2^{l-1}}{2l}\sum_{n=l}^{m}\binom{n-1}{l-1}\left(\pm\frac{2}{a}\right)^{n-2}\frac{pc^{n/2}L_{n-3}(pc)}{(pc-1)^{2n-3}}, \qquad l = 3,4,\ldots,m. \tag{56d}$$

The corresponding parameter values for the quadratic case, that is with $m = 2$, were originally derived in [14, eq. 52]. Equations (55-56) show how these same values are modified in order for the Hardy function to be expanded in terms of generalized $m^{th}$-order Gauss sums.

From the pattern of coefficients established above, it is now possible to tackle the important issue raised in both previous subsections concerning the collection size in the regime $(pc-1) < O(1)$. The first contribution to the exponential phase neglected in (55-56) is the $(m+1)^{th}$ term, which must be of $O(\varepsilon_t)$. Hence one requires

$$|2\pi\Phi_{m+1}k^{m+1}| < \left|2\pi\Phi_{m+1}\left(\frac{M_t}{2}\right)^{m+1}\right| = \varepsilon_t. \tag{57}$$

Using equation (56d) one can see that, to leading order in $\Phi_{m+1}$, this constraint will be equivalent to

$$\left|2\pi\frac{2^{m-1}}{(m+1)}\left(\frac{2}{a}\right)^{m-1}\frac{pc^{(m+1)/2}L_{m-2}(pc)}{(pc-1)^{2m-1}}\left(\frac{M_t}{2}\right)^{m+1}\right| = \varepsilon_t,$$

$$\Rightarrow M_t = \frac{1}{\sqrt{pc}}\left[\left(\frac{\varepsilon_t}{\pi}\right) a^{(m-1)} \frac{(pc-1)^{2m-1}}{2^{m-2}}\right]^{1/(m+1)} \left[\frac{m+1}{L_{m-2}(pc)}\right]^{1/(m+1)}. \tag{58}$$

In order to establish an appropriate collection size, one needs to estimate the size of $L_{m-2}(pc)$ for $pc$ close to unity. Assuming the structure of the polynomials (47) continues to follow the formulation postulated by (48) for all values of $m \geq 3$, one can easily devise the following constraints

$$L_{m-2}(1) \leq \frac{2^{2m-2}}{\sqrt{\pi}(m-2)^{\frac{3}{2}}} < L_{m-2}(pc) \;\Rightarrow\; \frac{1}{L_{m-2}(pc)} < \frac{\sqrt{\pi}(m-2)^{\frac{3}{2}}}{2^{2m-2}}. \tag{59}$$

Substituting (59) into (58) and assuming $\left[(m-2)^{3/2}(m+1)\right]^{1/(m+1)} \sim 1$, gives an upper bound on the collection size in this regime of

$$M_t(pc, m) \leq \left[\left(\frac{\varepsilon_t}{\sqrt{\pi}}\right) a^{(m-1)} \frac{(pc-1)^{2m-1}}{2^{3m-4}}\right]^{1/(m+1)}. \tag{60}$$

However, because of technical issues arising from the proposed computational algorithm developed in the next section, this upper bound (60) on the collection size is a little too optimistic. Hence to be absolutely certain that $M_t(pc, m)$ is not too large for computational purposes, a slightly reduction in the estimate (60) is imposed, giving a new collection size

$$M_t(pc, m) = \left[\left(\frac{\varepsilon_t}{\pi}\right) \frac{a^{(m-1)}}{2^{5(m-2)/2}} \left(\frac{\left(\sqrt{pc}-1\right)}{pc^{1/4}}\right)^{2m-1}\right]^{1/(m+1)}, \tag{61}$$

in the regime $(pc-1) < O(1)$. The key algebraic difference between (61) and the original collection size (27) is the factor $(pc-1)^{(2m-1)}$ in the former, compared to $\left(\sqrt{pc}-1\right)^{(2m-1)}$ in the latter. Hence this new collection size, applicable in the regime $(pc-1) < O(1)$, is much smaller than (27), as anticipated at the end of subsection 2.3. Substitution of this smaller collection size into equation (21) in turn gives a new smaller parameter value for $X$. Utilising these two reduced values, one finds that all the subsequent analysis set out in subsections 2.2 and 2.3 is now applicable for $pc$ values close to unity.

This minor modification provides final link in the chain connecting the Hardy function (9) to a series of generalized $m^{th}$-order Gauss sums. Substituting (55) into (51) gives

$$Z(t) = T_{a+\varepsilon} + \mathcal{H}(t)2\sqrt{2}\left\{\sum_{\substack{\alpha>a \\ odd}}^{a+b} \frac{\cos\left(\frac{t}{2}\left\{log(pc) + \frac{1}{pc} + 1\right\} + \frac{\pi}{8}\right)}{(\alpha^2 - a^2)^{\frac{1}{4}}}[1 + O(\Lambda_\alpha^{-1})] + R_{EM}\right\} +$$

$$\sum_{\substack{\alpha_E>a+b \\ step\ 4M_t^-(m)+4}}^{N_\alpha} \frac{\mathcal{H}(t)2\sqrt{2}}{(\alpha_E^2 - a^2)^{\frac{1}{4}}} \mathrm{Re}\left\{e^{\mathrm{i}\omega^+(pc)} S_{M_t^-(m)}(\Phi_{1..m}^+) + e^{\mathrm{i}\omega^-(pc)} S_{M_t^-(m)}(\Phi_{1..m}^-)\right\}[1 + O(\varepsilon_t)]. \quad (62)$$

Here $m \geq 2$ and $M_t^-(m) = (M_t(pc, m) - 1)/2$, with the latter specified by equations (61) and (27) in the regimes $(pc - 1) < O(1)$ and $(pc - 1) \geq O(1)$ respectively. The relative error scale in the first term satisfies $\Lambda_\alpha = MIN\{t^{1/8}(\alpha - a)^{3/4}, \sqrt{t}\}$, whilst the associated Gaussian sum parameters are given by equations (56a-d) above. These coefficients follow an orderly pattern, with both $\Phi_{1,2} = O(1)$ modulo $2\pi$, whilst the remainder decrease regularly $\Phi_{l\geq 3} = O\left(t^{(2-l)/2}\right)$ as the order of the sum increases (assuming $(pc - 1) \geq O(1)$, with minor modifications in the $(pc - 1) < O(1)$ regime).

This asymptotic structure in the higher order coefficients is significant interest. Nominally it is similar in character to the higher order Gauss sum formulation for $Z(t)$ loosely formulated in [7, eqs. 2.4, 2.15-2.17], in that the coefficients $\Phi_{l\geq 3} = O\left(t^{(2-l)/2}\right)$. However, for these particular sums, across the main computational regime where $(pc - 1) \geq O(1)$, the relative size of the corresponding $\Phi_{l\geq 3}$ coefficients in (62) declines quite significantly as $pc$ increases. This characteristic makes them much more advantageous from a computational viewpoint. The work of both [7] and [14] demonstrated the potential of recursive algorithms for the computation of quadratic Gauss sums in something like $O\left(log^{\leq 3}(N/\varepsilon)\right)ops_{s;c}$. It is the proposition of this work that the particular characteristics exhibited by the generalized Gauss sums (namely the smallness and relative sizes of their higher order coefficients $\Phi_{l\geq 3}$) constituting formulation (62), makes them amenable to similar rapid computation using analogous recursive algorithms to those devised in [14] for the quadratic case. In which case it should be possible to reduce the operational count necessary for the computation of $Z(t)$ quite substantially. This proposition is the subject of the next section.

## 3. Efficient computation of the particular Generalized Gauss Sums constituting $Z(t)$

The subject of Gauss sums has been a focus for study for many authors (e.g. [1], [17], [21], and references therein). This section is devoted to the derivation and implementation of an algorithm for the rapid computation of a Generalized $m^{th}$-order Gauss sum of the form

$$S_N(\Phi_{1..m}) = S_N(\Phi_1, \Phi_2, \Phi_3, \ldots, \Phi_m) = \sum_{k=0}^{N} exp[2\pi\mathrm{i}(\Phi_1 k + \Phi_2 k^2 + \Phi_3 k^3 + \cdots + \Phi_m k^m)]$$

(63)

where $\Phi_1 \in (-1/2, 1/2]$, $\Phi_2 \in (-1/4, 1/4]$, and the remaining coefficients are asymptotically small $\Phi_j = O(t^{(2-j)/2})$ for $j = 3, 4, \dots, m$. Typically the length of the sum is restricted (27) so that $N = O(\varepsilon^{1/(m+1)} t^{(m-1)/2(m+1)})$ for some given tolerance parameter $0 < \varepsilon \ll 1$. This means the final term in the exponent $\Phi_m k^m \le \Phi_m N^m = O(\varepsilon^{m/(m+1)} t^{1/(m+1)})$, and any potential succeeding term $\Phi_{m+1} k^{m+1} \le \Phi_{m+1} N^{m+1} = O(\varepsilon)$ is negligible. These conditions characterise the $S_N(\Phi_{1..m})$ sums constituting the Hardy function.

The methodology is similar to that devised in [21] and implemented in [14] for the particular case of a quadratic Gauss sum. However, there are some significant differences. The quadratic case the method relies on the well-known result that any quadratic Gauss sum can be written (approximately) in terms of a shorter quadratic Gauss sum by means of a quadratic reciprocity formula (e.g. [1], [22]). Employed recursively this allows one to progressively reduce $N$ to a value $L_{n_{K+1}} < K$, where $K$ is some suitably prescribed cut-off integer small enough to make the exact computation of $S_{L_{n_{K+1}}}(\Phi_{1;n_{K+1};}, \Phi_{2;n_{K+1}})$ a trivial exercise. From this starting point one can quickly compute an estimate for the original $S_N(\Phi_1, \Phi_2)$. For this much more ambitious implementation, difficulties arise because the analogous reciprocity formula for a $m^{th}$-order Gauss sum, whilst reducing its length, *does not necessarily keep the order fixed* (in fact there is a tendency for the order to increase). In general, this increase in complexity of the reduced sum disrupts the utility of the method. However, here the situation is rescued by the fact that all the initial coefficients $\Phi_{j\ge 3}$ are all small, which has the effect of delaying the onset of the increase in order for many iterations of the reciprocity formula. By monitoring the size of a certain remainder term $E_m$ (see section 3.1.3), one can anticipate an impending increase in order and account for it. This allows for further iterations of the formula to be implemented, reducing $N$ down to a value below $K$, much as for the quadratic case.

3.1 *Derivation of the reciprocity formula for generalized m^th^-order Gauss sums*

This work closely follows the methodology discussed in [12, section 2.2], although computation of the various remainder terms must be carried out to greater precision for the purposes of providing a suitable approximation for $S_N(\Phi_{1..m})$. The starting point is the Euler – Maclaurin summation formula [11, section A5, eq. A23; 6, section 2.1; 4, section 6.2, eq. 3], which, when applied to (63), gives

$$S_N(\Phi_{1..m}) = \sum_{k=0}^{N} e^*(f(k)) = \int_0^N e^*(f(y))dy + \frac{[e^*(f(0)) + e^*(f(N))]}{2}$$
$$-2\pi \mathrm{i} \int_0^N (y - \lfloor y \rfloor - 1/2) f'(y) e^*(f(y)) dy. \qquad (64)$$

Here $e(f(k)) \equiv exp(2\pi \mathrm{i} f(k))$ and $e^*(f(k))$ its complex conjugate. The function

$$f(y) = -\sum_{q=1}^{m} \Phi_q y^q, \tag{65}$$

is defined so that $f(y) \in C^2[0, N]$ and $f''(y) < 0$ in $[0, N]$ (hence the negative sign in 65). In the formulation that follows the parameter $x \equiv 2\Phi_2$ must be positive, which in practice is easy to arrange. Substituting the Fourier expansion [12, eq. A27]

$$(y - \lfloor y \rfloor - 1/2) = -\frac{1}{\pi}\sum_{n=1}^{\infty} n^{-1} sin(2\pi n y), \qquad y \notin \mathbb{Z}, \tag{66}$$

into the last term of (64) gives

$$2i\sum_{n=1}^{\infty}\frac{1}{n}\int_0^N sin(2\pi n y) f'(y) e^*\big(f(y)\big)dy = \sum_{n=1}^{\infty}\frac{1}{n}\int_0^N [e(ny) - e(-ny)] f'(y) e^*\big(f(y)\big)dy$$

$$= \sum_{n=1}^{\infty}\frac{1}{n}\int_0^N \big[e\big(ny - f(y)\big) - e\big(-ny - f(y)\big)\big] f'(y) dy$$

$$= \sum_{n=1}^{\lfloor \xi \rfloor}\frac{1}{n}\int_0^N \{\Sigma_q q\Phi_q y^{q-1}\}\, e\big(\Sigma_q \Phi_q y^q - ny\big)dy \tag{67a}$$

$$+ \sum_{n=\lceil \xi \rceil}^{\infty}\frac{1}{n}\int_0^N \{\Sigma_q q\Phi_q y^{q-1}\}\, e\big(\Sigma_q \Phi_q y^q - ny\big)dy \tag{67b}$$

$$- \sum_{n=1}^{\infty}\frac{1}{n}\int_0^N \{\Sigma_q q\Phi_q y^{q-1}\}\, e\big(ny + \Sigma_q \Phi_q y^q\big)dy, \tag{67c}$$

where $\xi = \Sigma_q q\Phi_q N^{q-1}$. Hence the generalised Gaussian sum satisfies

$$S_N(\Phi_{1..m}) = \int_0^N e^*\big(f(y)\big)dy + \frac{\big[e^*\big(f(0)\big) + e^*\big(f(N)\big)\big]}{2} + (67a, b, c). \tag{68}$$

The next step is to estimate the sums of integrals (67a-c) in much more detail. Term (67a) is by far the most difficult, since the integrands possess saddle points lying within the integration range $[0, N]$. Hence this term will be dealt with first.

3.1.1 *Estimation of Term (67a)*

This term is a finite sum which can be rewritten in the form

$$\sum_{n=1}^{\lfloor\xi\rfloor}\frac{1}{n}\int_0^N\left\{\Sigma_q q\Phi_q y^{q-1}-n+n\right\}e\left(\Sigma_q\Phi_q y^q-ny\right)dy$$

$$=\sum_{n=1}^{\lfloor\xi\rfloor}\int_0^N e\left(\Sigma_q\Phi_q y^q-ny\right)dy+\sum_{n=1}^{\lfloor\xi\rfloor}\frac{1}{n}\int_0^N\left\{\Sigma_q q\Phi_q y^{q-1}-n\right\}e\left(\Sigma_q\Phi_q y^q-ny\right)dy$$

$$=\sum_{n=1}^{\lfloor\xi\rfloor}\int_0^N e\left(\Sigma_q\Phi_q y^q-ny\right)dy+\frac{1}{2\pi\mathrm{i}}\sum_{n=1}^{\lfloor\xi\rfloor}\frac{\left[e\left(\Sigma_q\Phi_q N^q-nN\right)-1\right]}{n}$$

$$=\sum_{n=1}^{\lfloor\xi\rfloor}\int_0^N e\left(\Sigma_q\Phi_q y^q-ny\right)dy+\frac{\left[e\left(\Sigma_q\Phi_q N^q\right)-1\right]}{2\pi i}\sum_{n=1}^{\lfloor\xi\rfloor}\frac{1}{n}$$

$$=\sum_{n=1}^{\lfloor\xi\rfloor}\int_0^N e^*\left(ny-\Sigma_q\Phi_q y^q\right)dy-\frac{\mathrm{i}\left[e^*\left(f(N)\right)-1\right]}{2\pi}\{\psi(\lfloor\xi\rfloor+1)+\gamma\}, \tag{69}$$

where $\psi(y)$ is the Psi function [19, chapter 5] and $\gamma$ = 0.57721 ... is Euler's constant. Up to now everything is exact, but to compute the remaining sum the constituent integrals in (69) must be estimated. Let $z\in\mathbb{C}$ and define

$$g(z)=nz-\sum_{q=1}^{m}\Phi_q z^q, \tag{70}$$

and let $c$ denote the position of the saddle point, i.e. $g'(c)=0$ for a given $n$ (since $n\in[1,\lfloor\xi\rfloor]\Longrightarrow 0<c<N$ and $c$ lies inside the range of integration). It is also useful to define the parameter $\beta=c/N\in(0,1)$. Then any one of the integrals in (69) can be computed by contour integration along a path shown in Fig. 2a such that

$$\int_0^N e\big(g(y)\big)dy=\int_0^{\mathrm{i}c}e\big(g(z)\big)dz+\int_{C\searrow}e\big(g(z)\big)dz+\int_{N-\mathrm{i}(N-c)}^{N}e\big(g(z)\big)dz, \tag{71}$$

where the contour $C\searrow$ represents the line segment $z=c+se^{-\mathrm{i}\pi/4}$ for $s\in\left[-c\sqrt{2},(N-c)\sqrt{2}\right]$ traversing the saddle point. In the vicinity of the saddle

$$g(z)=g(c)+\sum_{p=2}^{m}\frac{g^{(p)}(c)(z-c)^p}{p!}\approx g(c)-\Phi_2(z-c)^2=g(c)+\mathrm{i}\Phi_2 s^2, \tag{72}$$

since for $p \geq 3, g^{(p)}(c)/p! \approx \Phi_p = O(t^{1-p/2}) \Longrightarrow \Phi_2 s^2 \gg \Phi_p s^p$. This means

$$\int_{C\searrow} e(g(z))dz \approx e^{2\pi i g(c)-\frac{i\pi}{4}} \int_{-c\sqrt{2}}^{(N-c)\sqrt{2}} e^{-2\pi\Phi_2 s^2} ds = e^{2\pi i g(c)-\frac{i\pi}{4}} \int_{-\beta N\sqrt{2}}^{(1-\beta)N\sqrt{2}} e^{-2\pi\Phi_2 s^2} ds. \tag{73}$$

Now in most cases $c$ lies somewhere in the middle of the range $[0, N]$, in which case the limits of integration in (73) are of $O(N)$, which can be approximated $\pm\infty$ to an exponentially small error. However, complications arise if $n$ is small, so that $\beta \approx 0$, or if $n$ is close to $\lfloor \xi \rfloor$, in which case $\beta \approx 1$. In such an instance the integration range is only semi-infinite, and greater precision is required. Utilisation of standard integration formulae [6, eqs. 3.321.2 & 3.322.1] gives the following results

$$\int_{C\searrow} e(g(z))dz \approx \frac{e^{2\pi i g(c)-\frac{i\pi}{4}}}{\sqrt{2\Phi_2}} \left[1 - \frac{1}{2}\mathrm{erfc}\{2\sqrt{\pi\Phi_2}\beta N\}\right], \quad n = 1, 2, 3, \ldots$$

$$\int_{C\searrow} e(g(z))dz \approx \frac{e^{2\pi i g(c)-\frac{i\pi}{4}}}{\sqrt{2\Phi_2}} \left[1 - \frac{1}{2}\mathrm{erfc}\{2\sqrt{\pi\Phi_2}(1-\beta)N\}\right], \quad n = \lfloor \xi \rfloor, \lfloor \xi \rfloor - 1, \lfloor \xi \rfloor - 2, \ldots$$

(74a, b)

In practice the $\mathrm{erfc}(z)$ terms in (74) die off exponentially fast and are typically negligible for $n > k$ and $n < \lfloor \xi \rfloor - k$ for $k \geq 4$.

The first integral in (71) can be computed along a path $z = s\mathrm{i}$ for $s = [0, c]$, which gives

$$\int_0^{\mathrm{i}c} e(g(z))dz = \mathrm{i}\int_0^c e(g(\mathrm{i}s))ds \approx \mathrm{i}\int_0^c e^{-2\pi(n-\Phi_1)s} e^{2\pi \mathrm{i}\Phi_2 s^2} ds. \tag{75}$$

When $n$ is small, so is the saddle point $c$ and the range of integration is very short, whilst if $n$ is large the $e^{-2\pi(n-\Phi_1)s}$ term cuts off the rest integrand very quickly. These observations enable one to approximate (75) by

$$\mathrm{i}\int_0^c e^{-2\pi(n-\Phi_1)s} e^{2\pi \mathrm{i}\Phi_2 s^2} dq \approx \mathrm{i}\int_0^c e^{-2\pi(n-\Phi_1)s}[1 + 2\pi \mathrm{i}\Phi_2 s^2 + \cdots] ds$$

$$= \frac{\mathrm{i}[1 - e^{-2\pi(n-\Phi_1)c}]}{2\pi(n-\Phi_1)} - \frac{\Phi_2[1 - e^{-2\pi(n-\Phi_1)c}\{2\pi^2(n-\Phi_1)^2 c^2 + 2\pi(n-\Phi_1)c + 1\}]}{2\pi^2(n-\Phi_1)^3}. \tag{76}$$

In a similar vein, the third integral in (71) can be approximated by

$$\int_{N-\mathrm{i}(N-c)}^{N} e\big(g(z)\big)dz = -\mathrm{i}\int_{N-c}^{0} e\big(g(N-\mathrm{i}s)\big)ds$$

$$\approx -\mathrm{i}e^{2\pi \mathrm{i}g(N)}\int_{N-c}^{0} e^{-2\pi(\xi-n)s}e^{2\pi\mathrm{i}\Phi_2 s^2}ds. \qquad (77)$$

If $n$ is close to $\lfloor\xi\rfloor$ then $c \approx N$ and the integration range is very short, whilst for smaller $n$ the $e^{-2\pi s(\xi-n)}$ term cuts off the rest of the integrand very quickly. By substituting $e^{2\pi\mathrm{i}\Phi_2 s^2} \approx [1 + 2\pi\mathrm{i}\Phi_2 s^2 + \cdots]$, it is easy to show that (77) can be approximated by an analogous expression to (76), viz.

$$\frac{\mathrm{i}e^{2\pi\mathrm{i}g(N)}}{2\pi}\left\{\frac{\left[1-e^{-2\pi(\xi-n)(N-c)}\right]}{(\xi-n)} + \frac{\mathrm{i}\Phi_2\left[1-e^{-2\pi(\xi-n)(N-c)}\{2\pi^2(\xi-n)^2(N-c)^2+2\pi(\xi-n)(N-c)+1\}\right]}{\pi(\xi-n)^3}\right\}. \qquad (78)$$

There is one further issue which concerns the very first integral of equation (68). If $\Phi_1 < 0$, the integrand contains a saddle $c > 0$, which effectively introduces a $n = 0$ term into the summation (69). (The case $\Phi_1 > 0$ is discussed below.) It can be accounted for by utilising the step-function $H(\Phi_1)$ to fix the starting value of $n$. When the two most significant terms in (76) and (78) (those $O((n-\Phi_1)^{-1})$ and $O((\xi-n)^{-1})$ respectively), are substituted back into (69), they give rise to the following sums

$$\sum_{n=H(\Phi_1)}^{\lfloor\xi\rfloor}\frac{1}{(n-\Phi_1)} = \psi(\lceil\xi\rceil-\Phi_1)-\psi(H(\Phi_1)-\Phi_1), \qquad (79a)$$

$$\sum_{n=H(\Phi_1)}^{\lfloor\xi\rfloor}\frac{e^{2\pi\mathrm{i}g(N)}}{(\xi-n)} = \sum_{n=H(\Phi_1)}^{\lfloor\xi\rfloor}\frac{e^{2\pi\mathrm{i}(nN-\sum_q\Phi_q N^q)}}{(\xi-n)} = e\big(f(N)\big)[\psi(H(\Phi_1)-\xi)-\psi(\lceil\xi\rceil-\xi)]. \qquad (79b)$$

Here $\lceil x\rceil = \lfloor x\rfloor + 1$ is the ceiling function. The remaining terms in (76) & (78) die off rapidly as $n = H(\Phi_1), 2, 3, \ldots$ and $n = \lfloor\xi\rfloor, \lfloor\xi\rfloor-1, \lfloor\xi\rfloor-2, \ldots$ respectively in much the same manner as the erfc$(z)$ appearing in (74). Combining the results (74), (76), (78) & (79) and substituting into (69) one arrives at an approximation for the complex conjugate of term (67a) in (68), given by

$$Term\ (67a)^* \approx \frac{e^{-\mathrm{i}\pi/4}}{\sqrt{2\Phi_2}}\sum_{n=H(\Phi_1)}^{\lfloor\xi\rfloor} e^{2\pi\mathrm{i}g(c)} - \frac{\mathrm{i}\left[e^{-2\pi\mathrm{i}f(N)}-1\right]}{2\pi}\{\psi(\lceil\xi\rceil)+\gamma\} + W_1, \qquad (80)$$

where

$$W_1 = \frac{\mathrm{i}}{2\pi}\Big[\psi(\lceil\xi\rceil - \Phi_1) - \psi(H(\Phi_1) - \Phi_1) + e^{2\pi \mathrm{i} f(N)}[\psi(H(\Phi_1) - \xi) - \psi(\lceil\xi\rceil - \xi)]\Big]$$

$$-\left[\frac{e^{2\pi \mathrm{i} g(c) - \mathrm{i}\pi/4}}{2\sqrt{2\Phi_2}}\mathrm{erfc}\left\{2\sqrt{\pi\Phi_2}c\right\} + \frac{\Phi_2}{2\pi^2(n-\Phi_1)^3} + \frac{\mathrm{i}e^{-2\pi(n-\Phi_1)c}}{2\pi(n-\Phi_1)} - \frac{e^{-2\pi(n-\Phi_1)c}}{2\pi^2}\left(\frac{\Phi_2\{2\pi^2(n-\Phi_1)^2c^2 + 2\pi(n-\Phi_1)c + 1\}}{(n-\Phi_1)^3}\right)\right]_{n\geq H(\Phi_1)}$$

$$-\left[\frac{e^{2\pi \mathrm{i} g(c) - \mathrm{i}\pi/4}}{2\sqrt{2\Phi_2}}\mathrm{erfc}\left\{2\sqrt{\pi\Phi_2}(N-c)\right\} + \frac{\Phi_2 e^{2\pi \mathrm{i} f(N)}}{2\pi^2(\xi-n)^3} + \frac{\mathrm{i}e^{2\pi \mathrm{i} f(N)}e^{-2\pi(\xi-n)(N-c)}}{2\pi(\xi-n)} - \frac{e^{2\pi \mathrm{i} f(N)}e^{-2\pi(\xi-n)(N-c)}}{2\pi^2}\left(\frac{\Phi_2\{2\pi^2(\xi-n)^2(N-c)^2 + 2\pi(\xi-n)(N-c) + 1\}}{(\xi-n)^3}\right)\right]_{n\leq\lfloor\xi\rfloor}. \quad (81)$$

In principle one could sum the latter terms in (81) over all $n$. But in practice they die off very quickly and are typically negligible when $n > P$ and $n < \lfloor\xi\rfloor - P$ for $P = 3$.

3.1.2 *Estimation of Terms (67b) and (67c)*

These two terms are somewhat easier to deal with than (67a) since the saddles of the constituent integrands fall outside the range of integration. Consider the complex conjugate of the integrals making up (67b). They can be computed by contour integration along a path shown in Fig. 2b, so that

$$\int_0^N f'(y)e\big(g(y)\big)dy = \int_0^{\mathrm{i}R} f'(z)e\big(g(z)\big)dz + \int_{\mathrm{i}R}^{N+\mathrm{i}R} f'(z)e\big(g(z)\big)dz + \int_{N+\mathrm{i}R}^{N} f'(z)e\big(g(z)\big)dz. \quad (82)$$

Since the solutions to the equation $g'(c) = 0$ give values for $c > N$ for $n > \lfloor\xi\rfloor$, it means that the integrand of the second integral $\left|e\big(g(s+\mathrm{i}R)\big)\right| < exp\big(-2\pi R(n-\xi)\big)$. Hence the integral itself $< Nf'(N+\mathrm{i}R)exp\big(-2\pi R(n-\xi)\big) \xrightarrow[R\to\infty]{} 0$. So only the first and third integrals actually contribute anything significant to the computation of (82). One has

$$\int_0^N f'(y)e\big(g(y)\big)dy = \mathrm{i}\int_0^\infty f'(\mathrm{i}s)e\big(g(\mathrm{i}s)\big)ds + \mathrm{i}\int_\infty^0 f'(N+\mathrm{i}s)e\big(g(N+\mathrm{i}s)\big)ds. \quad (83)$$

Both of these integrals will be dominated by their behaviour close to $s = 0$. Near this point $f'(\mathrm{i}s) \approx -\Phi_1$, $e\big(g(\mathrm{i}s)\big) \approx e^{-2\pi(n-\Phi_1)s}$, so that

$$\int_0^N f'(y)e\big(g(y)\big)dy \approx -\mathrm{i}\Phi_1 \int_0^\infty e^{-2\pi(n-\Phi_1)s}ds + \mathrm{i}\int_\infty^0 f'(N+\mathrm{i}s)e\big(g(N+\mathrm{i}s)\big)ds,$$

$$= \frac{-\mathrm{i}\Phi_1}{2\pi(n-\Phi_1)} + \mathrm{i}\int_\infty^0 f'(N+\mathrm{i}s)e\big(g(N+\mathrm{i}s)\big)ds. \qquad (84)$$

The remaining integral requires more precision because the saddle may lie very close to $N$, in the particular case $n = \lceil\xi\rceil$. In general, near $s = 0$, $e\big(g(N+\mathrm{i}s)\big) \approx e^{2\pi\mathrm{i}g(N)-2\pi(n-\xi)s+2\pi\mathrm{i}\Phi_2 s^2}$, whilst $f'(N+\mathrm{i}s) \approx -\xi - 2\mathrm{i}\Phi_2 s$. If $n > \lceil\xi\rceil$, the linear exponent will dominate the quadratic part and the integral is trivial. But if $n = \lceil\xi\rceil$, then $(n-\xi)s \sim \Phi_2 s^2$ or indeed smaller, so a more precise estimate is necessary. Hence the second integral in (84) can be approximated by

$$\approx -\mathrm{i}\xi e^{2\pi\mathrm{i}g(N)}\left[\int_\infty^0 e^{-2\pi(n-\xi)s}ds\right]_{n\geq\lceil\xi\rceil+1} - \mathrm{i}e^{2\pi\mathrm{i}g(N)}\int_\infty^0 (\xi+2\mathrm{i}\Phi_2 s)e^{-2\pi(\lceil\xi\rceil-\xi)s+2\pi\mathrm{i}\Phi_2 s^2}ds$$

$$= \left[\frac{\mathrm{i}\xi e^{2\pi\mathrm{i}g(N)}}{2\pi(n-\xi)}\right]_{n\geq\lceil\xi\rceil+1} + \delta_{n,\lceil\xi\rceil}\mathrm{i}e^{2\pi\mathrm{i}g(N)}\left[\lceil\xi\rceil\int_0^\infty e^{-2\pi(\lceil\xi\rceil-\xi)s+2\pi\mathrm{i}\Phi_2 s^2}ds - \frac{1}{2\pi}\right]. \qquad (85)$$

Here $\delta_{n,\lceil\xi\rceil}$ is the Kronecker delta and the second term is obtained by integration by parts. The remaining integral in (85) can be written in terms of the complementary error function [6, eq. 3.322.2]

$$\int_0^\infty e^{-2\pi(\lceil\xi\rceil-\xi)s+2\pi\mathrm{i}\Phi_2 s^2}ds = \frac{e^{\mathrm{i}\pi/4}}{2\sqrt{2\Phi_2}}e^{\mathrm{i}\pi(\lceil\xi\rceil-\xi)^2/2\Phi_2}\operatorname{erfc}\left[e^{\mathrm{i}\pi/4}\sqrt{\frac{\pi}{2\Phi_2}}(\lceil\xi\rceil-\xi)\right]. \qquad (86)$$

Substituting results (86) into (85) means (84) can be conveniently written as

$$\int_0^N f'(y)e\big(g(y)\big)dy \approx \frac{-\mathrm{i}\Phi_1}{2\pi(n-\Phi_1)} + \left[\frac{\mathrm{i}\xi e^{2\pi\mathrm{i}g(N)}}{2\pi(n-\xi)}\right]_{n\geq\lceil\xi\rceil} - \frac{\mathrm{i}\xi e^{2\pi\mathrm{i}g(N)}}{2\pi(\lceil\xi\rceil-\xi)}\delta_{n,\lceil\xi\rceil}$$

$$+\delta_{n,\lceil\xi\rceil}\mathrm{i}e^{2\pi\mathrm{i}g(N)}\left[\lceil\xi\rceil\frac{e^{\mathrm{i}\pi/4}}{2\sqrt{2\Phi_2}}e^{\mathrm{i}\pi(\lceil\xi\rceil-\xi)^2/2\Phi_2}\operatorname{erfc}\left[e^{\mathrm{i}\pi/4}\sqrt{\frac{\pi}{2\Phi_2}}(\lceil\xi\rceil-\xi)\right] - \frac{1}{2\pi}\right]. \qquad (87)$$

Multiplying (87) by minus one and taking the complex conjugate (q.v. definitions (65) and (70)), gives the following result

$$Term(67b) \approx \frac{-\mathrm{i}\Phi_1}{2\pi}\sum_{n=\lceil\xi\rceil}^{\infty}\frac{1}{n(n-\Phi_1)} + \frac{\mathrm{i}e^{-2\pi\mathrm{i}f(N)}\xi}{2\pi}\sum_{n=\lceil\xi\rceil}^{\infty}\frac{1}{n(n-\xi)} + W_2, \tag{88}$$

where

$$W_2 = \mathrm{i}e^{-2\pi\mathrm{i}f(N)}\left[\frac{e^{-\mathrm{i}\pi/4}}{2\sqrt{2\Phi_2}}e^{-\mathrm{i}\pi(\lceil\xi\rceil-\xi)^2/2\Phi_2}\mathrm{erfc}\left[e^{-\mathrm{i}\pi/4}\sqrt{\frac{\pi}{2\Phi_2}}(\lceil\xi\rceil-\xi)\right] - \frac{1}{2\pi(\lceil\xi\rceil-\xi)}\right]. \tag{89}$$

Now consider the constituent integrals of the term (67c). In a similar manner to (82) they can be integrated along the contour path shown in Fig. 2b to give

$$\begin{aligned}\int_0^N f'(y)e\left(ny+\Sigma_q\Phi_q y^q\right)dy &= \int_0^{\mathrm{i}R} f'(z)e\left(nz+\Sigma_q\Phi_q z^q\right)dz + \int_{N+\mathrm{i}R}^{N} f'(z)e\left(nz+\Sigma_q\Phi_q z^q\right)dz,\\ &\approx \mathrm{i}\int_0^{\infty}(-\Phi_1-2\mathrm{i}\Phi_2 s)e^{-2\pi(n+\Phi_1)s-2\pi\mathrm{i}\Phi_2 s^2}ds + \frac{\mathrm{i}\xi e^{-2\pi\mathrm{i}f(N)}}{2\pi(n+\xi)}.\end{aligned} \tag{90}$$

In most cases the first integral in (90) can be approximated trivially by $-\mathrm{i}\Phi_1/2\pi(n+\Phi_1)$. However, as in (84) this may become very inaccurate for the particular case when $n=1$ and $\Phi_1<0$, when the corresponding saddle point only just falls outside the integration range. So in this particular case, greater precision is required. Integration by parts gives

$$\left[\mathrm{i}\int_0^{\infty}e^{-2\pi(1+\Phi_1)s-2\pi\mathrm{i}\Phi_2 s^2}ds - \frac{\mathrm{i}}{2\pi}\right] = \left[\frac{\mathrm{i}e^{-\mathrm{i}\pi/4}}{2\sqrt{2\Phi_2}}e^{-\mathrm{i}\pi(1+\Phi_1)^2/2\Phi_2}\mathrm{erfc}\left[e^{-\mathrm{i}\pi/4}\sqrt{\frac{\pi}{2\Phi_2}}(1+\Phi_1)\right] - \frac{\mathrm{i}}{2\pi}\right]. \tag{91}$$

Substituting (91) into (90), dividing by $n$ and summing gives the convenient form

$$Term(67c) \approx \frac{\mathrm{i}e^{-2\pi\mathrm{i}f(N)}\xi}{2\pi}\sum_{n=1}^{\infty}\frac{1}{n(n+\xi)} - \frac{\mathrm{i}\Phi_1}{2\pi}\sum_{n=1}^{\infty}\frac{1}{n(n+\Phi_1)} + W_3, \tag{92}$$

where

$$W_3 = H(-\Phi_1)\left[\frac{\mathrm{i}e^{-\mathrm{i}\pi/4}}{2\sqrt{2\Phi_2}}e^{-\mathrm{i}\pi(1+\Phi_1)^2/2\Phi_2}\mathrm{erfc}\left[e^{-\mathrm{i}\pi/4}\sqrt{\frac{\pi}{2\Phi_2}}(1+\Phi_1)\right] - \frac{\mathrm{i}}{2\pi(1+\Phi_1)}\right]. \tag{93}$$

Here the use of the step function $H(-\Phi_1)$ acts to remove the inaccurate approximation for $n = 1$ and $\Phi_1 < 0$ and replace it by its more precise formulation (91). One final point again arises from the first integral in eq. (68). The case $\Phi_1 < 0$ was dealt with in section 3.1.1 above. However, if $\Phi_1 > 0$ the saddle point now falls just outside the integration range, in which case it makes an additional $n = 0$ contribution to term (67c). In form it is identical to integral (90), but with $n = 0$ and $f'(y) \leftrightarrow 1$. Using the results (90) and (91) leads to the following result

$$W_4 = H(\Phi_1)\left[\frac{\mathrm{i}e^{-\mathrm{i}\pi/4}}{2\sqrt{2\Phi_2}} e^{-\mathrm{i}\pi\Phi_1^2/2\Phi_2}\mathrm{erfc}\left[e^{-\mathrm{i}\pi/4}\sqrt{\frac{\pi}{2\Phi_2}}\Phi_1\right] - \frac{\mathrm{i}e^{-2\pi\mathrm{i}f(N)}}{2\pi\xi}\right]. \tag{94}$$

Equations (88) and (92) provide two convenient expressions for terms (67b) and (67c) respectively. When the summation terms appearing in these respective equations are themselves added together, one obtains

$$\frac{\mathrm{i}e^{-2\pi\mathrm{i}f(N)}\xi}{2\pi}\left[\sum_{n=1}^{\infty}\frac{1}{n(n+\xi)} + \sum_{n=\lceil\xi\rceil}^{\infty}\frac{1}{n(n-\xi)}\right] = \frac{\mathrm{i}e^{-2\pi\mathrm{i}f(N)}\xi}{2\pi}\left[\sum_{n=1}^{\lfloor\xi\rfloor}\frac{1}{n(n+\xi)} + \sum_{n=\lceil\xi\rceil}^{\infty}\frac{2}{(n^2-\xi^2)}\right]$$

$$= \frac{\mathrm{i}e^{-2\pi\mathrm{i}f(N)}}{2\pi}\left[(\psi(\lceil\xi\rceil) - \psi(\lceil\xi\rceil+\xi) + \psi(\xi+1) + \gamma) + (\psi(\lceil\xi\rceil+\xi) - \psi(\lceil\xi\rceil-\xi))\right]$$

$$= \frac{\mathrm{i}e^{-2\pi\mathrm{i}f(N)}}{2\pi}[\psi(\lceil\xi\rceil) + \psi(\xi+1) + \gamma - \psi(\lceil\xi\rceil-\xi)], \tag{95$a$}$$

and

$$-\frac{\mathrm{i}\Phi_1}{2\pi}\left[\sum_{n=1}^{\infty}\frac{1}{n(n+\Phi_1)} + \sum_{n=\lceil\xi\rceil}^{\infty}\frac{1}{n(n-\Phi_1)}\right] = -\frac{\mathrm{i}[\psi(\lceil\xi\rceil) + \psi(\Phi_1+1) + \gamma - \psi(\lceil\xi\rceil-\Phi_1)]}{2\pi}$$

(95b)

Incorporating together these two expressions (95a,b) alongside equation (80) for term (67a), results in the factors involving $\{\psi(\lceil\xi\rceil) + \gamma\}$ cancelling out exactly. This means on combining the expressions (80), (88) and (92) into (68), one obtains the following approximation for $S_N(\Phi_{1..m})$.

$$S_N(\Phi_{1..m}) \approx \frac{e^{-\mathrm{i}\pi/4}}{\sqrt{2\Phi_2}}\sum_{n=H(\Phi_1)}^{\lfloor\xi\rfloor} e^{2\pi\mathrm{i}g(c)} + \frac{[1 + e^{-2\pi\mathrm{i}f(N)}]}{2} + W_1 + W_2 + W_3 + W_4 + W_5, \tag{96}$$

where

$$W_5 = \frac{\mathrm{i}}{2\pi}\left[e^{-2\pi\mathrm{i}f(N)}[\psi(\xi+1) - \psi(\lceil\xi\rceil-\xi)] - [\psi(\Phi_1+1) - \psi(\lceil\xi\rceil-\Phi_1)]\right], \tag{97}$$

and $W_1$, $W_2$, $W_3$ and $W_4$ are given by equations (81), (89), (93) and (94) respectively. (N.B. only one of $W_3$ or $W_4$ ever contributes to $S_N(\Phi_{1..m})$, depending upon the sign of $\Phi_1$.) This completes all the results necessary to estimate the $m^{th}$order Gauss sum $S_N(\Phi_{1..m})$ from (68). Since all of the $W_{1-5}$ terms are straightforward to calculate, the utility of the approximation (96) depends upon the characteristics of the secondary sum over $n \in [H(\Phi_1), \lfloor \xi \rfloor]$.

3.1.3 *Analysis of the secondary sum in (96)*

The secondary sum in (96) is characterised by an exponent term defined by

$$g(c) = nc - \sum_{q=1}^{m} \Phi_q c^q \quad with \;\; n - \Phi_1 = \sum_{q=2}^{m} q\Phi_q c^{q-1}. \tag{98}$$

In general, unlike for the quadratic case, equation (98) has no closed form solution. However, because the particular Gaussian sums of interest here (63) have higher order coefficients $\Phi_j = O(t^{1-j/2})$ for $j = 3, 4, \dots, m$, it is possible to make progress. Let $y = n - \Phi_1 > 0$, $x = 2\Phi_2 \in (0, 1/2)$ and define

$$c_{m-1} = \sum_{q=1}^{m-1} C_q(\Phi_{2,3,\dots q+1}) \frac{(q+1)y^q}{x^{2q-1}}, \tag{99}$$

to be the solution for $c$ satisfying

$$y = \sum_{q=2}^{m} q\Phi_q c^{q-1}, \tag{100}$$

truncated at the $m^{th}$- order in $y$. (The constant coefficients $C_q(\Phi_{2,3,\dots q+1})$ are complicated functions of the $\Phi_j$, but are relatively easy to calculate using any reasonably efficient computer algebra package. $C_1 = 1/2$, $C_2 = -\Phi_3$ and the specific values for $C_{3-9}$ are listed in the Appendix.) Substituting (99) into (98) one obtains

$$g_m(c_m) = yc_m - \sum_{q=2}^{m} \Phi_q c_m^q = \sum_{q=1}^{m-1} \frac{C_q y^{q+1}}{x^{2q-1}} + \sum_{q=m}^{\infty} E_q, \tag{101}$$

for the exponent $g(c)$, truncated at the $(m+1)^{th}$- order in $y$. For example, the two simplest cases for $m = 3$ & $4$ give

$$c_3 = \frac{y}{x} - \frac{3\Phi_3 y^2}{x^3} + \frac{4(9\Phi_3^2 - 2\Phi_4 x)y^3}{2x^5}; \quad c_4 = c_3 - \frac{5(27\Phi_3^3 - 12\Phi_3\Phi_4 x + \Phi_5 x^2)y^4}{x^7}, \tag{102}$$

$$g_3 = \frac{y^2}{2x} - \frac{\Phi_3 y^3}{x^3} + \sum_{q=3}^{\infty} E_q\,; \qquad g_4 = \frac{y^2}{2x} - \frac{\Phi_3 y^3}{x^3} + \frac{(9\Phi_3^2 - 2\Phi_4 x)y^4}{2x^5} + \sum_{q=4}^{\infty} E_q, \tag{103}$$

$$E_3 = \frac{(9\Phi_3^2)y^4}{2x^5};\; E_4 = -\frac{(27\Phi_3^3 - 12\Phi_3\Phi_4 x)y^5}{x^7};\; E_q = \frac{C_q(\Phi_{2,3,\ldots q}, \Phi_{q+1} = 0)y^{q+1}}{x^{2q-1}}. \tag{104}$$

Now the crucial point concerns the size of the remainder terms $E_q$. Provided these terms decrease sufficiently fast to ensure that $\sum_q E_q$ converges to a value below some specified *small* tolerance factor (typically $\varepsilon = O(\varepsilon_t)$), then (96) would give a suitable approximation for $S_N(\Phi_{1..m})$ in terms of another $m^{th}$-order Gauss sum of considerably shorter length $\lfloor \xi \rfloor$. The constraints imposed on $N = O\left(\varepsilon_t^{1/(m+1)} t^{(m-1)/2(m+1)}\right)$ and $\Phi_j = O\left(t^{1-j/2}\right)$ $j = 3, \ldots m$ initially, ensure that

$$\lfloor \xi \rfloor = \left\lfloor \sum_{q=1}^{m} q\Phi_q N^{q-1} \right\rfloor = O(2\Phi_2 N) = O(xN) \tag{105}$$

with $0 < x < 1/2$. So in principle, (96) could then be applied recursively, each iteration reducing the number of summands by at least a factor of two, until one reaches some appropriate termination scale $O(K) \ll N$ governing the length of the kernel sum in something like $n_K = \lfloor log(N/K)/log(2) \rfloor + 1$ steps. To see if this is practical, one needs to estimate the remainder terms.

The respective remainder terms $E_q$ for general $q$ depend upon nature of the $C_q(\Phi_{3,4,\ldots q})$ coefficients, which become increasingly complicated as $q$ increases (see Appendix). These coefficients give rise to a sequence of terms, commencing in the following specific fashion

$$E_q = \sum_{n=1}^{q-2} E_{q,n} = \binom{2q-2}{q-1}\frac{y^{q+1}(3\Phi_3)^{q-1}}{q(q+1)x^{2q-1}} - \binom{2q-2}{q-1}\frac{y^{q+1}(3\Phi_3)^{q-3}}{q(q+1)x^{2q-2}}\left[2\binom{q-2}{1}\Phi_4\right]$$
$$+\binom{2q-4}{q-2}\frac{y^{q+1}(3\Phi_3)^{q-5}}{(q-1)q(q+1)x^{2q-3}}\left[30\binom{q-2}{2}\Phi_3\Phi_5 + 48\binom{q-2}{3}\Phi_4^2\right] -$$
$$\binom{2q-4}{q-2}\frac{y^{q+1}(3\Phi_3)^{q-7}}{(q-1)q(q+1)x^{2q-4}}\left[54\binom{q-3}{2}\Phi_3^2\Phi_6 + 180\binom{q-3}{3}\Phi_3\Phi_4\Phi_5 + 128\binom{q-3}{4}\Phi_4^3\right] + \cdots$$

(106)

Initially, one would like to simplify the above expression in order to derive an estimate for $|E_q|$ (specifically the $max|2\pi E_q|$ when $y = \lfloor \xi \rfloor - \Phi_1$) in terms of $x, y$ and $\Phi_3$ alone, and subsequently demonstrate that this estimate is sufficiently small to ensure $\sum_{q \geq m} E_q$ converges. It would then be possible, for prescribed set of values of $N$ and $\Phi_{3,..m}$, to postulate an upper bound on $|2\pi E_q|$ to ensure that (101) can be safely truncated at the $(m+1)^{th}$- order in $y$.

To derive an upper bound on $|E_q|$, consider the case of a fourth order Gauss sum in which the coefficients $\Phi_{5,6,7,..} = 0$. In that case (106) simplifies to

$$E_q = \frac{\mathrm{E}(x, y, \Phi_3, q)}{2^{q-2}} \left[ \left(1 - \frac{2x\Phi_4}{(3\Phi_3)^2}\right)^{q-2} - \sum_{j=\lfloor (q+1)/2 \rfloor}^{q-2} \binom{q-2}{j} \left(-\frac{2x\Phi_4}{(3\Phi_3)^2}\right)^j \right] + O\left(\frac{\mathrm{E}}{q}\right),$$

where

$$E(x, y, \Phi_3, q) = \frac{y^2}{2x} \binom{2q-2}{q-1} \left(\frac{6\Phi_3 y}{x^2}\right)^{q-1} \frac{1}{q(q+1)}. \qquad (107a, b)$$

The modulus of $[\,]/2^{q-2}$ term in (107a) is always less than one with the proviso $|2x\Phi_4/(3\Phi_3)^2| < 1$. For cubic Gauss sums, to be utilised here, $\Phi_4 = 0$ and the result is trivial. But even in the case of quartic Gauss sums, substitution of expressions (56c & d) gives $|2x\Phi_4/(3\Phi_3)^2| = x(pc+1)(pc-1)/2pc$ for $pc < 2$ and $(pc+1)/2pc$ for $pc \geq 2$. Since $x < 1$ the proviso is always satisfied in this case too. Hence one concludes that $|E_q|$ is bounded by $|E(x, y, \Phi_3, q)|$ specified by (107b).

Application of Stirling's formula [19] provides an upper bound to the binomial coefficient in (107b), so that

$$E(x, y, \Phi_3, q) \underset{q \to \infty}{\longrightarrow} \frac{y^2}{2x\sqrt{\pi}} \left(\frac{24\Phi_3 y}{x^2}\right)^{q-1} \frac{1}{(q+1)(q-1/2)\sqrt{q}}, \qquad (108)$$

a result that means that $\sum_{q \geq m} E_q < \sum_{q \geq m} |E(x, y, \Phi_3, q)|$ converges, provided that $|24\Phi_3 y/x^2| < 1$. Now there are two distinct cases, depending upon the size of $pc$. First the case when $pc$ is large, in particular when $pc > 3$. In this regime $x = 2\Phi_2 \approx (pc-1)^{-1} \in (0, 1/2)$ and $y \leq \lfloor \xi \rfloor - \Phi_1 \approx x M_t^- \approx x M_t/2$. Substitution of these results, with $M_t$ specified by (27) and $\Phi_3$ by the dominant term of (56d), into (107) and (108) yields

$$\left(\frac{24\Phi_3 y}{x^2}\right) = \left[\sqrt{2^{m+9}} \frac{\varepsilon_t}{t} \frac{(pc+1)^m pc^2}{(pc-1)^{m+2}}\right]^{1/(m+1)} \ll 1, \qquad (109)$$

$$\left|\sum_{q=m}^{\infty} E_q\right| < \sum_{q=m}^{\infty} |E| < \frac{y^2}{2x\sqrt{\pi}} \sum_{q=m}^{\infty} \left[\sqrt{2^{m+9}} \frac{\varepsilon_t}{t} \frac{(pc+1)^m pc^2}{(pc-1)^{m+2}}\right]^{(q-1)/(m+1)} \Big/ q^{5/2}. \quad (110)$$

The initial term in this sum (when $q = m$) is clearly the largest, with the remainder dying off *at least* as fast as $m^{5/2}, (m+1)^{-5/2}, (m+2)^{-5/2}, \ldots$ Substituting the expressions for $x$ and $y$ into (110) gives

$$\left|\sum_{q=m}^{\infty} E_q\right| < \frac{\varepsilon_t}{2} \sqrt{\frac{2^{m-7}}{(\pi m)^3}} \left[\frac{(pc+1)^m}{(pc-1)^{m-1} pc}\right] \varsigma\left(\frac{5}{2}\right) = O(\varepsilon_t), \quad (111)$$

for $pc > 3$ up to at least a cut-off value of $m = 10$, more than sufficient for (101) to be applicable to cubic sums. (Reducing $M_t$ by a factor of $4\sqrt{2}$ extends the validity of this $O(\varepsilon_t)$ result to cover all cut-off values for $pc > 2$.)

Now consider the case when $1 < pc < 2$ and in particular in the limit $pc \longrightarrow 1$. The analysis is similar to the above, with the caveat that the approximation $x \approx (pc-1)^{-1}$ is no longer valid. Rather $x \in (0, 1/2)$ is now computed from (56c) as the absolute difference of $(pc-1)^{-1}$ from its nearest integer. Substituting $y \leq xM_t/2$ as before, but with the smaller collection size (61) appropriate for the case $O(pc-1) < 1$, one finds that

$$\left(\frac{24\Phi_3 y}{x^2}\right) = \frac{1}{x\left(\sqrt{pc}+1\right)^3} \left[\sqrt{2^{3(m+4)}} \frac{\varepsilon_t}{t} \frac{pc^{(4m+7)/4}}{\left(\sqrt{pc}-1\right)^{m+4}}\right]^{1/(m+1)}. \quad (112)$$

As with (109), this factor must be less than one to ensure that the sum $\sum_{q \geq m} E_q$ converges. Which presents a problem in the limit $pc \longrightarrow 1$, since the $\left(\sqrt{pc}-1\right)$ term in the denominator has the potential to make (112) arbitrarily large. In practice the problem is not quite so acute, since Gauss sum contributions in (62) only commence when $O(pc-1) = t^{-m/4(m+2)}$ for all $m \geq 3$; contributions associated with $pc$ values closer to unity being covered by the transition term $T_{a+\varepsilon}$ (see section 2.4 and eq. 128 below). However, this still means that one cannot commence any calculations with Gauss sums of order $m > 4$ - at least initially - since factor (112) would exceed unity. To overcome this issue, one must instead rewrite those terms in (4) with saddles satisfying $(pc-1) \in \left[t^{-m/4(m+2)}, t^{-1/(m+4)}\right]$ as *quadratic* Gauss sums, and then employ the second order scheme based upon [33] and implemented previously in [28] to estimate them rapidly. Indeed in this regime with $pc$ close to unity, condition (112) implies that one must employ the proposed iteration scheme judiciously, commencing with a set of second order sums, then third order, then fourth order…until reaching the full $m^{th}$-order sums desired. The schematic below highlights the organisation of the transitions from successive orders as $pc$ gradually increases from unity.

$$m = 3: \quad (pc-1) = O\left(t^{-3/20}\right) \xrightarrow[2nd-order\ scheme]{} O\left(t^{-1/7}\right) \xrightarrow[3rd-order\ onwards]{}$$

$$m = 4: \quad O\left(t^{-1/6}\right) \xrightarrow[2nd-order]{} O\left(t^{-1/7}\right) \xrightarrow[3rd-order]{} O\left(t^{-1/8}\right) \xrightarrow[4th-order\ onwards]{}$$

$$m = 5: \quad O\left(t^{-5/28}\right) \xrightarrow[2nd-order]{} O\left(t^{-1/7}\right) \xrightarrow[3rd-order]{} O\left(t^{-1/8}\right) \xrightarrow[4th-order]{} O\left(t^{-1/9}\right) \xrightarrow[5th-order\ onwards]{}$$

$$m = n: \quad O\left(t^{-n/4(n+2)}\right) \xrightarrow[2nd-order]{} O\left(t^{-1/7}\right) \xrightarrow[3rd-order]{} O\left(t^{-1/8}\right) \ldots O\left(t^{-1/(n+4)}\right) \xrightarrow[n^{th}-order\ onwards]{}$$

In practice for the case $m = 3$ the difference between $t^{-3/20}$ and $t^{-1/7}$ is insignificant for computationally feasible values of $t$ and in this instance one can estimate all the necessary sums using a 3rd-order scheme from the beginning. But for $m > 3$ it would be necessary to work up to a $m^{th}$-order scheme progressively, in the manner illustrated. When $O(pc-1) = t^{-1/(m+4)}$ is reached, the factor (112) will be a term of $O\left(t^{(2-m)/[(m+1)(m+4)]}\right)$, small enough to satisfy the convergence criterion for $\sum_{q\geq m} E_q$. Substituting (112) into (108) and following same analysis as for (110-111), one finds that for $1 < pc < 2$

$$\left|\sum_{q=m}^{\infty} E_q\right| < \varepsilon_t \frac{2^{(5m-4)/(m+1)}}{(\pi m)^{3/2}} \left[\left[\frac{2}{\sqrt{pc}+1}\right]\right]^{3(m-1)} \frac{pc^{4m}\varsigma(5/2)}{\left[x\left(\sqrt{pc}-1\right)\right]^{m-2}} \quad \forall\, m \geq 3. \quad (113)$$

In theory the $\left(\sqrt{pc}-1\right)$ term in the denominator of (113) has the potential to make this bound on the remainder terms unacceptably large. However, in practice one only begins employment of the approximation scheme when the collection size $M_t$ exceeds the prescribed cut-off integer $K$. By judicious choice of the latter, it is easy to ensure that the approximation scheme only commences when $E_m$ and hence $\sum_{q\geq m} E_q$ falls below the toleration factor $\varepsilon$ (typically after one or two pivot blocks – section 4.2 below). In summary this analysis demonstrates that formula (96), combined with approximation (101), can be utilised to estimate the specific Gauss sums that appear in (62).

The results presented above rest on the premise that the original $m^{th}$-order Gauss sum is to be estimated from (96) by the calculation of another (shorter) $m^{th}$-order Gauss sum of $O(xN)$ in length. Subsequently the estimation procedure is applied recursively, until one reaches some suitable kernel sum of $O(K) \ll N$, which is trivial to compute. But this premise unnecessarily restrictive. As the recursion process proceeds and the respective sum lengths are reduced, the ratio $\Phi_3 y/x^2$ in (108) gradually increases, which in turn will cause the size of the remainder term (101) to grow, to potentially unacceptable level. (For a concrete example of this see the discussion in section 3.4.2 below.) However, provided the condition $|24\Phi_3 y/x^2| < 1$ still holds, this growth can be regulated and the remainder term kept to a suitable size, simply by allowing the order of the sum to *increase*, when necessary, during the iterative process. This extra flexibility is perfectly possible, since (96) is valid for all orders. The proviso $|24\Phi_3 y/x^2| < 1$ is sufficient to ensure $E_m > E_{m+1} > \cdots > E_{m+n-1} > \varepsilon > E_{m+n}$, for some value of $n$, at which point the next iteration in the recursive procedure can be implemented for the shorter, but higher $(m+n)^{th}$-order sum. This does mean that for each iteration the size of $E_m$ has to be carefully monitored, which requires some extra calculation. But in practice this extra overhead is tiny, since the growth in $\Phi_3 y/x^2$ is relatively slow and the order only increases in ones or twos as the iteration process proceeds (see Tables 1-3 below).

3.2 *A recursive scheme for the efficient computation of a Generalized* $m^{th}$*-order Gaussian sum* $S_N(\Phi_{1..m})$*, satisfying coefficient conditions (63).*

This scheme in large part follows the formulation for the estimation of the quadratic Gauss sum discussed in [14]. As there let $S_N^{s=+1}(\Phi_{1..m}) \equiv$ the standard sum and $S_N^{s=-1}(\Phi_{1..m}) \equiv$ its complex conjugate and let $x = 2\Phi_2$. Now initially both $\Phi_1$ and $x$ are arbitrary real numbers, whilst $\Phi_{q=3..m} = O(t^{1-j/2})$ as prescribed after (63). To shift $\Phi_1$ and $x$ into the range $(-1/2, 1/2]$, set $x_1 = x - \lfloor x \rceil, \Phi_{2,1} = x_1/2$, $L_1 = N$, $m_1 = m$ and

$$\Phi_{1,1} = \text{sgn}(x_1)\{\Phi_1 - \lfloor \Phi_1 \rceil \pm (1 - (-1)^{\lfloor x \rceil})/4\}, \qquad \Phi_{2..m_1,1} = \text{sgn}(x_1)\Phi_{j=2..m}, \qquad (114)$$

where the choice of $\pm$ is made to ensure that $|\Phi_{1,1}|$is inside the range $(0, 1/2]$. Then $S_{L_0}^{s_1}(\Phi_{1..m_1,1})$ with $s_1 = \text{sgn}(-x_1)$ and $\Phi_2 > 0$, is in the form specified by equations (64-65).

Now on each successive iteration the initial default position is that the order of the next generation Gauss sum appearing on the right-hand side of (96) remains the same as previous generation on the left-hand side, that is $m_l = m_{l-1}$. Since $y = n - \Phi_1$ in equation (101), the coefficients of this next generation sum are given by

$$\Phi_{q,l} = \sum_{j=max(q,2)}^{m_l} \binom{j}{q} \left(-\Phi_{1,l-1}\right)^{j-q} \frac{C_{j-1}\left(\Phi_{2,3,\ldots j,l-1}\right)}{x_{l-1}^{2j-3}}, \quad \text{for } q = 0,1,2,..,m_l. \qquad (115)$$

So in particular if $m_l = 3$ one would have

$$\Phi_{0,l} = \frac{\Phi_{1,l-1}^2}{2x_{l-1}} + \frac{\Phi_{1,l-1}^3 \Phi_{3,l-1}}{(x_{l-1})^3}, \qquad \Phi_{1,l} = -\frac{\Phi_{1,l-1}}{x_{l-1}} - \frac{3\Phi_{1,l-1}^2 \Phi_{3,l-1}}{(x_{l-1})^3},$$

$$\Phi_{2,l} = \frac{x_l}{2} = \frac{1}{2x_{l-1}} + \frac{3\Phi_{1,l-1}\Phi_{3,l-1}}{(x_{l-1})^3}, \qquad \Phi_{3,l} = -\frac{\Phi_{3,l-1}}{(x_{l-1})^3}. \qquad (116)$$

The appearance of the $\Phi_{0,l>1}$ terms means that each successive next generation Gauss sum is multiplied by a constant phase factor of $e^{2\pi i \Phi_{0,l}}$ in (96). The reduced sum length $L_l$ associated with these next generation coefficients $\Phi_{q,l}$ is given by

$$L_l = \lfloor \xi_{l-1} \rfloor = \left\lfloor \sum_{q=1}^{m_{l-1}} q\Phi_{q,l-1} L_{l-1}^{q-1} \right\rfloor, \qquad \Phi_{2,l-1} > 0, \quad l \geq 2. \tag{117}$$

At this point the suitability of the order of the next proposed Gauss sum must be checked. This is done by first setting $y_{max} = L_l - \Phi_{1,l-1}$ and then computing

$$E_{m_l} = \frac{C_{m_l}(\Phi_{2,3,\ldots m_l,l-1}, \Phi_{m_l+1,l-1} = 0) y_{max}^{m_l+1}}{x_{l-1}^{2m_l-1}}. \tag{118}$$

Provided $|2\pi E_{m_l}| < \epsilon$ some pre-set tolerance level (typically $\epsilon \sim 10^{-3}$ is suitable), then the current order level is fine. However, if $|2\pi E_{m_l}| > \epsilon$ then this indicates that $m_l$ is currently too small for the next sum length $L_l$. At which point the order must be increased to $m_l \equiv m_l + 1$. A factor $\Phi_{m_l,l-1} = 0$ (which did not exist at the previous order) is then introduced, enabling operations (115) and (118) to be repeated for the new order value, until a point is reached when $|2\pi E_{m_l}| < \epsilon$. For coefficients satisfying the initial conditions prescribed here, $|2\pi E_{m_l}|$ increases incrementally, which means that $m_l$ increases, if at all, in single increments. However, as the sum length $L_l$ approaches $O(10 - 100)$ the value of $|2\pi E_{m_l}|$ starts to increase more rapidly, which means the order may need to be increased by two or more increments. Once a suitable value of $m_l$ has been established, the next generation factors $\Phi_{1,l}$ and $x_l$ are then shifted into the desired $(-1/2, 1/2]$ range by re-defining $x_l = x_l - \lfloor x_l \rceil$ and applying (114). One can see from equations (115-116) that the presence of powers of $x_l$ in the numerators means that each successive iteration results in a gradual increase in the corresponding values of the parameters $\Phi_{q=3\ldots m_l,l}$, although the ranking $\Phi_{3,l} \gg \Phi_{4,l} \gg \ldots \gg \Phi_{m_l,l}$ remains unchanged. Should any factor $\Phi_{3,l}$, $\Phi_{4,l}$ ... ever exceed unity at any point, then they too can be shifted into range $(-1/2, 1/2]$ by redefining $\Phi_{q,l} = \Phi_{q,l} - \lfloor \Phi_{q,l} \rceil$, without effecting the next value of the sum in the chain. Combining all these results gives rise to the following algorithm for the estimation of $S_N(\Phi_{1..m})$.

3.3. *An algorithm for the recursive computation of a Generalized $m^{th}$-order Gauss Sum $S_N(\Phi_{1..m})$, with $N \gg 1$, $m \geq 3$ and coefficients $\Phi_{1..m}$ satisfying the specific conditions set out in equation (1).* (*algorithm MGS*)

*Replace $S_N(\Phi_{1..m})$ by $S_{L_1}^{s_1}(\Phi_{1..m_1,1})$ using the definitions given in and around equation (114). Fix positive integer parameters $K > 30$, $P$ ($P = 3$ or $4$ are fine) and small real $\epsilon > 0$ (typically $\epsilon \leq 10^{-2}$ is suitable). Let $\omega = e^{i\pi/4}$.*

*Compute for $l = 2, 3, \ldots$*

- *a)* *Set $m_l = m_{l-1}$. Compute $L_l$ from $\Phi_{1\ldots m_l,l-1}$ and $L_{l-1}$ using (116). Set $y_{max} = L_l - \Phi_{1,l-1}$.*
- *b)* *Compute $\Phi_{0\ldots m_l,l}$ from $\Phi_{0\ldots m_l,l-1}$, using (115).*

c) $x_l = 2\Phi_{2,l}$; $x_l = x_l - \lfloor x_l \rfloor$; $\Phi_{2,l} = x_l/2$ ; $s_l = \mathrm{sgn}(-x_l)$; *Adjust each* $\Phi_{1,2\ldots m_l,l}$ *using (114) to ensure they all lie in the range* $(-1/2\,,1/2]$ *and multiply each one by* $\mathrm{sgn}(x_l)$, *so in particular* $\Phi_{2,l} > 0$.

d) *Compute the* $max\left|2\pi E_{m_l}\right|$ *using (118).*

e) *If* $\left|2\pi E_{m_l}\right| < \epsilon$ *carry on to f). Otherwise set* $m_l = m_l + 1$, $\Phi_{m_l,l-1} = 0$ *and repeat b)-e).*

f) *If* $L_{l+1} < L_l < K$ *and* $L_l > \lceil \Phi_{1,l-1} \rceil$, *set* $l_K = l$, *otherwise set* $l_K = l - 1$. *If* $L_{l+1} < 0 < K < L_l$, *set* $l_K = l$. *If* $K < L_l < L_{l+1}$, *set* $l_K = l$. *If none of these, repeat a)-f).*

g) *Compute the sum* $S_{L_{l_K}}^{s_{l_K}}\left(\Phi_{1..m_{l_K},l_K}\right)$ *exactly. N.B. if* $\Phi_{1,l_K-1} > 0$ *the summands commence at* $k = 1$ *rather than* $k = 0$, *to account for the step function in (96). This adjustment must also be after each iteration prescribed in h) below.*

h) *From this starting point compute, for* $l = l_K - 1, l_K - 2, \ldots, 1$, *the following iterations*

$$S_{L_l}^{s_l}\left(\Phi_{1..m_l,l}\right) = \left\{\frac{e^{2\pi i \Phi_{0,l+1}} S_{L_{l+1}}^{s_{l+1}}\left(\Phi_{1..m_{l+1},l+1}\right)}{\omega\sqrt{|x_l|}} + \frac{\left[1 + e^{-2\pi i f(L_l)}\right]}{2} + W_1 + W_2 + W_{3\ or\ 4} + W_5\right\}^{s_l},$$

*where* $W_{1-5}$ *are given by equations (81), (89), (93), (94) and (97) respectively.*

i) *The final iteration* $S_{L_1}^{s_1}\left(\Phi_{1..m_1,1}\right)$ *gives an estimate to* $S_N(\Phi_{1..m}) = S_{L_1}^{s_1}\left(\Phi_{1..m_1,1}\right) + \delta(K, N, \Phi_{1..m})$, *where* $\delta(K, N, \Phi_{1..m})$ *is the error term. The magnitude of this and the corresponding relative error, denoted by* $\varepsilon_{MGS}(K)$ *are discussed in Section 3.5.*

The termination conditions set out in point f) require some comment. The first set of conditions represent the ideal, (such an example is shown in Table 1). In this instance the initial exact sum (the kernel sum) computed in *g)* has a length $L_{l_K} < K$ the predesignated cut-off integer. However, this may not be attainable. Problems can occur when $L_{l-1}$ approaches values $K = O(10 - 100)$ in equation (117). Initially and certainly for most iterations $l \le l_K$ the terms in (117) satisfy $2\Phi_{2,l-1}L_{l-1} \gg 3\Phi_{3,l-1}L_{l-1}^2 \gg \ldots \gg q\Phi_{q,l-1}L_{l-1}^{q-1}$. Since $\Phi_2 > 0$ this means $L_l > 0$. But as the iteration procedure progresses the size of the respective $\Phi_{q=3\ldots m_l,l}$ parameters increases; both in magnitude and relative to the size of $\Phi_{2,l}$. Sometimes as $L_l$ approaches $K$, and in particular if the next $|x_l|$ is somewhat small (typically less than $10^{-1}$) producing a relatively steep drop from $L_l$ to $L_{l+1}$, this can mean $2\Phi_{2,l}L_l < 3\Phi_{3,l}L_l^2$ even though $\Phi_{2,l} > \Phi_{3,l}$. If in such cases $\Phi_{3,l} < 0$ this means $L_{l+1} < 0$ too (see Table 2). It is unclear whether such a Gauss sum over *negative* integers down to $L_{l+1}$ would be a correct interpretation for the kernel sum in such instances. Possibly, but since all the analysis presented in Section 3.1 assumes the summands are over non-negative integers, it is sensible to avoid such difficulties by moving the kernel sum one place further up the chain so that $L_{L_K} > 0$, even if this means $L_{L_K} > K$. Similar sorts of problems can lead to situations where $K < L_l < L_{l+1}$ if $K$ is set to a very small value. The stipulation that $K > 30$ is designed to screen them out, but if they do occur then the $(l+1)$th iteration is discarded and the kernel sum commences at the $l$th iteration.

3.4 *Sample computations using algorithm MGS*

3.4.1 *Cubic Gauss sums* $(m = 3)$

Application of algorithm *MGS* to the specific cubic $(m = 3)$ Gauss sums defined by (63) is obviously the easiest to implement. Tables 1 & 2 provide two detailed illustrative cases of how the estimation procedure works. The two particular Gauss sums shown in the tables were generated by first fixing an arbitrary value for $t = 10^{28}$ and then setting the associated pivot integers (see section 2.1) to $\alpha_E = NINT_E(2.9a/\sqrt{8})$ and $\alpha_E = NINT_E(5.5a/\sqrt{8})$ respectively. These choices illustrate the procedure for the cases when the portcullis variable $pc = 1.566726 \ldots$ is relatively small (Table 1) and $pc = 13.048362 \ldots$ relatively large (Table 2) respectively.

For a pre-set value of error tolerance and cut-off integer, in these examples $\varepsilon_t = 0.02$ and $K = 200$ respectively, the coefficients of the sum to be estimated (the progenitor) can then be calculated from equations (56b-d), along with the lengths $N = (M_t(pc, 3) - 1)/2$ derived from their respective collection sizes $M_t(pc, 3)$ given by (27). These values are shown in the first rows of the respective tables and define the sum to be estimated. The next set of rows show the associated coefficients of the respective Gaussian sums resulting from the systematic reduction in the sum length $L_l$ by means of the iterative procedure set out in equations (115-117). Together these form a hierarchical chain or tree of Gauss sums linking the progenitor to a kernel sum of vastly reduced length. Initially the chain consists of sums of the same order as the progenitor. However, as the chain proceeds the order starts to increase as the length is reduced, resulting in 5th order (Table 1) and 6th order (Table 2) kernel sums respectively. Increases in order are not necessarily linear but can jump by two (as illustrated in Table 2) or indeed more increments, particularly as the length of the sum falls into the $O(10 - 100)$ range. Once established, the kernel sum is computed exactly and then successive estimates of all the constituent sums of the hierarchical chain generated from (96). The final iteration yields the estimate for the desired progenitor.

The tables show a comparison of the various estimates with the exact values of the respective hierarchical sums, alongside the associated actual $\delta(K, N, \Phi_{1..m}, P)$ and relative $\varepsilon_{MGS}(K)$ error terms respectively. Typically, whilst the magnitude of the error $\delta(K, N, \Phi_{1..m}, P)$ builds up with each iteration, the relative error, after some initial fluctuations, stabilises (see section 3.5 below) to a value well below the tolerance. Notice that in Table 2 the length of the kernel sum is greater than $K$. When attempting to compute the next iteration, the relatively small size of $\Phi_{2,5}$ combined with fact $\Phi_{3,5} < 0$, means that the resulting sum length $\lfloor \xi_6 \rfloor = \Sigma_q q\Phi_{q,5}273^{q-1} < 0$. Consequently, the hierarchical chain has to commence with a kernel of length greater than $K$ (as discussed in section 3.3).

Another point arising from these tabulated results is that the actual number of iterations making up the hierarchical chains is much the same in both cases, despite the fact that the sum in Table 2 has around thirty times the summands as that in Table 1. In practice there is a tendency for the number of iterations to slowly decrease with portcullis variable. Although both the collection size and hence $N$ increase by a factor $pc^{1/2}$ (see 27), this is counterbalanced, beyond $pc > 3$, by the fact that the first iteration reduces $N$ by a factor $x_1 = 2\Phi_{2,1} \approx 1/(pc - 1)$ from equation (56c). This means that the first iteration produces an increasingly large initial reduction in $N$ of $O(pc^{-1/2})$, which has the effect of curtailing the length of the hierarchical chains as $pc \to \infty$. For values of $pc < 3$ (Table 1), $N$ is reduced in a more regular pattern $N/|a_1|, N/|a_1 a_2|, \dots, N/|a_1 a_2 \dots a_{l_K-1}|$, governed by the partial quotients $|a_0; a_1 a_2 \dots|$ of the nearest integer continued fraction representation [7] for $x_1$.

3.4.2 *Quartic Gauss Sums* ($m = 4$): *the reductive bound on cut-off* $K$.

Algorithm *MGS* clearly works well for the specific cubic Gauss sums prescribed in (63). The next question is to establish its applicability to higher order sums of similar specifications. Table 3 shows an illustrative example for a quartic Gauss sum centred upon the same pivot integer $\alpha_E = NINT_E(2.9a/\sqrt{8})$ utilised for the cubic example of Table 1. Application of algorithm *MGS* produces a reduction of the length of the progenitor sum by a factor $\sim 300$, over the first four iterations whilst increasing the order from four to eight. Operationally this is all analogous to the cubic examples shown in Tables 1 & 2 above. However, at this stage one encounters a serious issue. In order to reduce the length of the sum still further, the algorithm attempts to increase the order of the next progenitor to point where the remainder term falls below the specified tolerance, i.e. $|2\pi E_m| < \varepsilon$. Implicit in this step is the assumption that the magnitudes of $E_m$ decrease with $m$. But at this point in the iteration process, this is no longer the case. Substitution of the tabulated parameters defining the 8th - order Gauss sum of iteration four in Table 3, yields remainder values of $2\pi E_8 \approx 3.7 \times 10^9$ and $2\pi E_9 \approx 8.9 \times 10^9$, with further increases for larger values $m$.

In order for the iteration scheme to continue work down the hierarchical chain the condition $|24\Phi_3 y/x^2| < 1$ must be satisfied at each iteration step, to ensure that $E_m > E_{m+1}$ in equations (107-108). The combination of collection size (27), coupled with $\Phi_3$ defined by (56d), means that this condition is fine initially. However, subsequently whilst the sum length $L_l$ is decreased by a factor $x_{l-1}$ the value of $\Phi_{3,l}$ increases by a factor of $x_{l-1}^3$ (116), meaning that as one moves down the hierarchical chain the product $\Phi_{3,l} L_l$ increases by a net factor of around $x_{l-1}^{-2}$. Hence $\Phi_{3,l} L_l$ eventually approaches a value $O(x_l)$, as exemplified in Table 3, violating the condition necessary to proceed further. Ideally at this point $L_l < K$, the prescribed cut-off integer, in which case no further iterations would be necessary. Routine estimation shows that $\Phi_{3,l} L_l = O(x_l)$ occurs when

$$L_l = O\left(\left[\frac{\varepsilon_t (pc+1)^m}{(pc-1)pc^{m-1}}\right]^{3/2(m+1)} t^{(m-2)/2(m+1)}\right). \tag{119}$$

This scale effectively imposes a *reductive bound* on the value of $K$. If $K$ is fixed at a value below (119), the iteration process will cease before the length of the progenitor sum has been reduced sufficiently. Notice too that this is a problem that will ultimately impact on the reduction of cubic $m = 3$ sums, such as those illustrated in Tables 1 & 2. But for values of $t < 10^{40}$ currently computationally feasible, this is not going to be a significant issue.

The question now arises, can the iteration procedure be modified to allow further reductions in the current sum length $L_4$ shown in Table 3, to a point below (or a least close to) the specified cut-off integer $K = 200$? The answer is *qualified* yes, since one can always undertake one further iteration beyond the point when $\Phi_{3,l}L_l = O(x_l)$ is reached. Without going into all the mathematical details here, the procedure employs numerical methods to calculate the saddle points $c(n) \in [0, L_l]$ satisfying (100), rather than establishing closed form solutions in terms of $n$ directly, by means of (102) etc. Typically, $n \in [0, L_{l+1}]$ with $L_{l+1} \sim x_l L_l < L_l$ and so the number of summands is reduced, in much the same manner as previous iterations. (The situation is complicated somewhat by the relatively large size of $\Phi_3$, which sometimes means there are two saddle point solutions $c_1(n)$ and $c_2(n)$ arising from a single value of $n$. This in turn means that the contour of integration shown in Fig. 2a must be modified, to pass through both saddle points.) The advantage of the procedure is that the number and values of the Gauss sum coefficients remains the same as in the previous iteration $\Phi_1, \ldots, \Phi_l$. The downsize is the numerical computation of the various $c(n)$, alleviated to some extent by the implementation of efficient polynomial root finding algorithms. The lower section of Table 3 illustrates the implementation. The resulting extra iteration, number five here, creates a new kernel sum consisting of $n = 1108$ summands (compared to $n = 12020$ previously), where the function $g(c(n))$ in (100) is now only implicitly dependent upon the value of $n$. The value of this kernel sum is then utilised to estimate the original progenitor sum, to comparable accuracy. However, to further reduce the value of $n$ closer to the specified value of $K = 200$ is more problematic, since it necessitates the sub-division of the kernel sum in the first instance, a procedure that impacts on the accuracy of the final estimates. Indeed, for the $Z(t)$ calculations based upon cubic $m = 3$ progenitor sums, to be discussed in Section 4, even the single extra iteration procedure outlined here was not implemented. Rather it was simply more efficient to slowly raise the value of $K$ as $t$ increases.

3.5 *Error Analysis*

In [14] a detailed analysis of the errors that accumulate in an algorithm *QGS* designed to provide rapid estimates of generalised quadratic Gauss sums. Since algorithm *MGS* is essentially an extension of algorithm *QGS* (for general parameter values) and works in entirely analogous manner, most of the analysis in [14] carries over to this work. What follows is a brief synopsis of that work, assuming the hierarchal chain is made up of just *single* pre-cursor sums (so no sub-division as envisaged at the end of the previous sub-section). Hence in practice this section is primarily concerned with the estimation errors accruing in the estimates for the type of cubic Gaussian sums illustrated in section 3.4.1. It will also be assumed that the parameter $\Phi_2$ is a 'generic' irrational, characterised by the absence of any particularly large partial quotients *early on* in its nearest integer continued fraction representation. Specifically, early on means no partial quotient $|a_n| > L_l$ within the first $l \in [1, L_K]$ iterations needed to reduce $N$ to a number less than $K$. This is the most common scenario in which the algorithm reduces the progenitor sum length in fairly regular steps down to the kernel. The case when $\Phi_2$ 'non-generic' irrational in the quadratic case is discussed in [14]. In that instance algorithm *QGS* can terminate abruptly with relatively large reduction in the sum length from well above $K$ to well below $K$ in just a single iteration. Such a scenario is likely to be *far less* common for more complex sums because the new sum length $\lfloor \xi \rfloor$ also depends upon the parameter values $\Phi_3, \Phi_4$ etc.

Let $\varepsilon_l$ represent the respective error term associated with the $l^{th}$ iteration of the algorithm, defined by

$$\varepsilon_l(K) = S_{L_l}^{s_l}(\Phi_{1..m_l,l}) - \left\{\frac{e^{2\pi i\Phi_{0,l+1}}S_{L_{l+1}}^{s_{l+1}}(\Phi_{1..m_{l+1},l+1})}{\omega\sqrt{|x_l|}} + \frac{[1+e^{-2\pi i f(L_l)}]}{2} + CW\right\}^{s_l}, \tag{120}$$

where $CW = W_1 + W_2 + W_{3\ or\ 4} + W_5$ derived earlier. The $l = l_K$ sum is computed exactly so $\varepsilon_{l_K} = 0$. Moving back up the chain the next cumulative error is $|\varepsilon_{l_K-1}|$, followed by $\left||\varepsilon_{l_K-1}|/|x_{L_K-2}|^{1/2} + |\varepsilon_{l_K-2}|\right|$ etc. Continuing in this vein, one finds that the final error $\delta(K, N, \Phi_{1..m})$ in the estimate for the progenitor sum satisfies

$$|\delta(K,N,\Phi_{1-m})| \leq \frac{|\varepsilon_{L_K-1}|}{\sqrt{|x_1 \dots x_{L_K-2}|}} + \frac{|\varepsilon_{L_K-2}|}{\sqrt{|x_1 \dots x_{L_K-3}|}} + \dots + \frac{|\varepsilon_2|}{\sqrt{|x_1|}} + |\varepsilon_1|,$$

$$\leq \frac{max\{|\varepsilon_{1\dots L_K-1}|\}}{\sqrt{|x_1 \dots x_{L_K-2}|}}\left\{1 + \sqrt{|x_{L_K-2}|} + \dots + \sqrt{|x_1 \dots x_{L_K-2}|}\right\} < \frac{3.4max\{|\varepsilon_{1\dots L_K-1}|\}}{\sqrt{|x_1 \dots x_{L_K-2}|}}. \tag{121}$$

The last result comes about because all the $|x_i| < 1/2$, so the sum in { } brackets cannot exceed $\sqrt{2}/(\sqrt{2}-1) \approx 3.4$. Since $L_{n_K-1} > K$ the 'generic irrational restriction ensures that $\left|S_{L_{l_K-1}}^{s_{l_K-1}}(\Phi_{1..m_{l_K-1},l_K-1})\right| \geq \sqrt{K}$, in which case

$$|\delta(K,N,\Phi_{1-m})| < 3.4max\{|\varepsilon_{1\dots L_K-1}|\}\left|\frac{S_{L_1}^{s_1}(\Phi_{1..m_1,1})}{S_{L_{l_K-1}}^{s_{l_K-1}}(\Phi_{1..m_{l_K-1},l_K-1})}\right| \leq \frac{3.4max\{|\varepsilon_{1\dots L_K-1}|\}|S_N(\Phi_{1..m})|}{\sqrt{K}}. \tag{122}$$

The corresponding relative error

$$\varepsilon_{MGS} = \frac{|\delta(K,N,\Phi_{1-m})|}{|S_N(\Phi_{1..m})|} < \frac{3.4max\{|\varepsilon_{1\dots L_K-1}|\}}{\sqrt{K}}, \tag{123}$$

then follows automatically.

If one were to substitute the intermediate hierarchal sums in place of $|S_N(\Phi_{1..m})|$ in equations (122) and (123), then these results predict that whilst the actual errors in the intermediate sum estimate scales with sum length, the relative errors will remain fixed at something below $K^{-1/2}$. Those predictions are born out in the behaviour of the error terms for the examples illustrated in Tables 1-3. Indeed, the actual relative error is typically ten times smaller the general upper bound (123). To get a precise upper bound one would need to establish an estimate for the $max\{|\varepsilon_{1\dots L_K-1}|\}$ term. The greatest source of error usually lies in the estimate for the term $W_1$ given by equation (81), arising from the integrals in term (67a).

Since these get more complicated with the order of the sum, the $max\{|\varepsilon_{1\dots L_K-1}|\}$ is usually associated with the error in the first iteration $\varepsilon_{L_K-1}$ beyond the kernel sum. In the tables one can see the relative error jumping up in the first couple of iterations, before stabilising thereafter. In the quadratic case discussed in [14] it was possible to put a specific bound on this term. For the higher order sums discussed here it much more difficult to formulate such a bound. Rather than pursue this, one can take the pragmatic view that more precise estimates can be simply computed by including more terms in the approximations developed in section 3.1.1. As ever, greater precision imposes a higher operational count, and one must balance the demands for high computational speed against high accuracy.

In conclusion, the results presented in this section demonstrate that the specific type of Generalized Gaussian sums which constitute the Hardy function in (62) can be computed very efficiently by means of algorithm *MGS*. Large scale calculations of the Hardy function itself, utilising this methodology, are presented in the next section.

## 4. Sample Computations of $Z(t)$ for large $t$ utilising Cubic Gauss Sums

In this final section an algorithmic methodology for the computation of the Hardy function based on the representation (62) derived in section 2 in terms of generalized $m^{th}$-order Gauss sums. It follows very closely the methodology derived for the $m = 2$ case involving quadratic Gauss sums (*algorithm ZT13* of [14]). The basic idea is to utilise algorithm *MGS* to compute the respective generalized sums very rapidly. Typically for a sum of length $N$ and cut-off integer $K$, the number of operations required (excluding any subdivision as envisaged in section 3.4.2) is of $O(2K + \Omega^* log(N/K))ops_{s;c}$ to reach an estimate within the prescribed relative error $\varepsilon_t$. (An estimate for the constant $\Omega^* \approx 31$, was devised in [14].) Obviously, this is a vast saving over the $2N ops_{s;c}$ required to compute (1) term by term. The implication of such a saving is the potential to reduce the number of operations necessary to compute a single evaluation of $Z(t)$ from the $O(\sqrt{t})ops_{s;c}$ needed for the classical Riemann-Siegel formula, down to a value $O\left((t/\varepsilon_t)^{1/(m+1)}(log(t))^2\right) ops_{s;c}$ for a precision of $O(\varepsilon_t)$ in the relative error. However, this estimate is far too optimistic, due to the reductive bound issue, and in practice the actual performance is $O\left((t/\varepsilon_t)^{[0.25,\,0.3]}(log(t))^{2+o(1)}\right) ops_{s;c}$ when $m = 3$ (see section 4.3 below).

In what follows, the algorithmic methodology is devised for arbitrary $m \geq 3$. Sample computations of $Z(t)$, for computationally feasible values of $t \lesssim 10^{36}$, can be obtained for the case of $m = 3$ by applying algorithm *MGS* as it stands. Modifications to the algorithm *MGS* to alleviate the constrains imposed by the reductive bound on the value of $K$, which would be necessary for implementing the methodology for values beyond $m = 3$, will be the subject of future work.

### 4.1 *Formulation of a hybrid summation formula for the estimation of* $Z(t)$

Applying the estimation methodology developed in the previous section directly to the formula (62) linking the Hardy function to sequences of $m^{th}$- order Gauss sums would be computationally inefficient. Since the index $N_\alpha \geq INT_O(t/\pi) + 2$, which, loosely, means that the number of summands in (62) amounts to $N_\alpha/2M_t = O(t^{(m+3)/2(m+1)})$. This exceeds the $O(\sqrt{t})$ summands required to compute $Z(t)$ by classical means. However, this estimate is an oversimplification, which disguises the fact that formula (62) can be utilised to estimate the summation of the overwhelming majority of the terms making up the Riemann-Siegel formula (3) very efficiently. In particular, suppose one were to aim to estimate $Z(t)$ in something like an operational count of just the $O\left((t/\varepsilon_t)^{1/(m+1)}\left(log(t)\right)^2\right) ops_{s;c}$ proposed above. The crux of the method to try and achieve this is the utilisation of portions of *both* sums (3) and (62) together, in order to create a computationally efficient hybrid summation formula for $Z(t)$ in its entirety.

To begin, define $N_C(m) = \lfloor (t/\varepsilon_t)^{1/(m+1)} log(t)/\sqrt{\pi} \rfloor$ as suitable cut-off integer. This is used to truncate (3) into two sections. The first section consisting of the lower $N \in [1, N_C]$ summands goes to form the first part of the hybrid summation formula, a computation which, of course, requires $N_C ops_{s;c}$. The second part consists of a subsection of those terms in (62) whose sum corresponds to the sum over the upper $N \in [N_C + 1, N_t]$ summands of the Riemann-Siegel formula. This subsection of terms is obtained by truncating (62) at a cut-off integer $\alpha_C(m)$ corresponding to the value of $N_C$ defined above. During the course of the analysis [13] in the derivation of representation (4) for the Hardy function, the author was able to demonstrate how the $\alpha$ summands in (9) are linked to corresponding values of $N$ forming the summands of (3), via the portcullis variable. Specifically, the analysis of [13, Appendix B, eq. B12] identifies the connection $N = a/4\sqrt{pc(\alpha)}$ (for real $\alpha \geq a$). Combining this result with equation (13) and eliminating $pc$ produces a direct relation $\alpha(N)$, given by

$$\alpha = 2N\left(\frac{t}{2\pi N^2} + 1\right) \Rightarrow N = \frac{\alpha}{4}\left[1 - \sqrt{1 - \left(\frac{a}{\alpha}\right)^2}\right]. \qquad (124)$$

Substituting $N_C$ into equation (124) gives the corresponding cut-off integer $\alpha_C(m)$ as

$$\alpha_C(m) = \left\lfloor 2N_C\left(\frac{t}{2\pi N_C^2} + 1\right) \right\rfloor \Rightarrow \alpha_C(m) \approx \frac{(\varepsilon_t t^m)^{1/(m+1)}}{\sqrt{\pi} log(t)}. \qquad (125)$$

Hence the hybrid summation formula for the Hardy function is defined by

$$Z(t) = 2\sum_{N=1}^{N_C(m)} \frac{cos\left(\theta(t) - tlog(N)\right)}{\sqrt{N}} - (-1)^{N_t}\left(\frac{2\pi}{t}\right)^{\frac{1}{4}} \sum_{r=0}^{M}\left(-\sqrt{\frac{2\pi}{t}}\right)^r \Psi_r(p) + ZP(t,m), \qquad (126)$$

where the second partial sum, $ZP(t,m)$, is given by a truncation of (62) in the form

$$ZP(t,m) = T_{a+\varepsilon} + \mathcal{H}(t)2\sqrt{2}\left\{\sum_{\substack{\alpha>a \\ odd}}^{a+b} \frac{\cos\left(\frac{t}{2}\left\{log(pc) + \frac{1}{pc} + 1\right\} + \frac{\pi}{8}\right)}{(\alpha^2 - a^2)^{\frac{1}{4}}}[1 + O(\Lambda_\alpha^{-1})]\right\} +$$

$$\sum_{\substack{\alpha_E > a+b \\ step\ 4M_t^-(m)+4}}^{\alpha_C(m)} \frac{\mathcal{H}(t)2\sqrt{2}}{(\alpha_E^2 - a^2)^{\frac{1}{4}}} \mathrm{Re}\{e^{\mathrm{i}\omega^+(pc)} S_{M_t^-(m)}(\Phi_{1..m}^+) + e^{\mathrm{i}\omega^-(pc)} S_{M_t^-(m)}(\Phi_{1..m}^-)\}[1 + O(\varepsilon_t)], \quad (127)$$

for $m \geq 2$. Now the number of summands in the second of the sums making up $ZP(t,m)$ is approximately $\alpha_C/4M_t^- \approx \alpha_C/2M_t$. Provided the pivot integers $\alpha_E$ are sufficiently greater than $a$ to ensure $(pc - 1) \geq O(1)$, then the sum lengths $M_t$ are specified by (27), which means that

$$\frac{\alpha_C}{2M_t(m)} \sim \left[\frac{2^{(m-3)} pc^{(3m-1)/2}}{(pc^2 - 1)^m}\right]^{1/(m+1)} \frac{\sqrt{t}}{log(t)} \sim \left(\frac{\sqrt{2^{(3m-5)}t}}{\varepsilon_t}\right)^{1/(m+1)}. \quad (128)$$

This final estimate is obtained by substituting the portcullis value at the cut-off integer, i.e. $pc(\alpha_C) \approx 4\alpha_C^2/a^2$. Since the number of summands in (127) is $O\left((t/\varepsilon_t)^{1/(m+1)}\right)$ and estimates for the respective Gaussian sums $S_{M_t^*(m)}\left(\Phi_{1..m}^{\pm}\right)$ can be computed using the methodology of section 3 in $O\left(2K + 2\Omega^* log(M_t^*/K)\right) ops_{s;c}$, computational estimates for $ZP(t,m)$ can be achieved at the expense of a relatively small operational count.

Before continuing, it should be emphasized the proposed hybrid representation (126-127) of the Hardy function is not exact. The formulation of $Z(t)$ in terms of Gaussian sums in section 2 is itself only accurate to $O(\varepsilon_t)$ in the relative error, whilst the computational methodology for the Gaussian sums themselves will introduce further errors of a similar size. So compared to the high precision calculations presented in the papers of [2, 7, 9], which are typically accurate to 5-6 decimal places, the methodology presented here is somewhat crude. However, the algorithms of [6, 19] require at least an $O\left(t^{1/3+o_\lambda(1)}\right)$ operational count to achieve this level of accuracy. The simplest formulation of the hybrid (126-127) with $m = 3$, upon which the current sample computations are based, already provides superior performance in terms of a reduced operational count (section 4.3). At the same time it offers the potential for faster algorithms based upon larger values of $m$.

4.2 *Computational scheme for the hybrid summation formula*

The partial sum $ZP(t,m)$, which forms the core of the hybrid summation formula consists of three main terms, the first two of which are easily dealt with. The transition term $T_{a+\varepsilon}$ is negligible (exponentially small) except on those rare occasions when the parameter $a$ happens to lie close to an odd integer; specifically when $|NINT_O(a) - a| < \varepsilon \leq t^{-1/6}$. If $\varepsilon > 0$ then the summands in the second term of $ZP(t,m)$ commence at $\alpha = INT_O(a) + 4$, rather than $INT_O(a) + 2$, the usual scenario. The 'missing' $INT_O(a) + 2$ summand is then replaced by $T_{a+\varepsilon}$. For $\varepsilon = g t^{-1/6} > 0$ and small $g \in [0, 1/4]$ one can derive a specific algebraic form for the transition term, given by

$$T_{a+\varepsilon} = \mathcal{H}(t)\frac{2^{3/4}\Gamma(1/3)exp\left(-(32\pi^3)^{1/4}\,g^{3/2}/3\right)}{3^{2/3}\pi^{1/4}t^{1/12}}\left[\cos\left(t+\sqrt{\frac{\pi}{2}}t^{1/3}g+\frac{\pi}{24}\right)+O(g)\right]. \quad (129)$$

However, if $g \in [1/4\,,1]$ or $\varepsilon < 0$ then $T_{a+\varepsilon}$ must be calculated numerically by (simple) numerical integration. Full details of the form and calculation of $T_{a+\varepsilon}$ can be found in [13, Appendix A].

The next step is to fix the value of the upper summand of the second term in $ZP(t,m)$. There is some degree of flexibility regarding the scale of parameter $b$ (it must be an even integer). For consistency with computations presented in the earlier work [14] it will be defined as $b = 2\lfloor t^{1/(m+2)}\rfloor$. Given this value the sum of the second term can easily be computed. This leaves the evaluation of the crucial third term of $ZP(t,m)$, based around the Gaussian sums. Since $b = 2\lfloor t^{1/(m+2)}\rfloor$, the lowest value pivot integers $\alpha_E$ commence in a regime where $(pc-1)\sim O\left(t^{-m/4(m+2)}\right)$. As discussed at the end of section 3.1, the order of the Gaussian sums must be stepped up in a regulated fashion, from two to three to four etc. up to the designated value of $m$, at which point $(pc-1)\sim O\left(t^{-1/(m+4)}\right)$. However, for the $m = 3$ cubic case, which is the main focus here, the difference between $t^{-3/20}$ and $t^{-1/7}$ is insignificant (at least for the practical large values of $t \in 10^{23-36}$). Hence in this instance, $3^{rd}$ order Gauss sums are employed immediately.

The first stage of the procedure covers those pivot integers lying in the range $\alpha_E \in [a+b, Ya]$, where $Y$ is a constant parameter set close to unity, but sufficiently large to ensure that for any $\alpha_E > Ya \Rightarrow (pc-1) \geq O(1)$. In practice any value of $Y > 1.06$ is usually sufficient, and in these particular calculations a value of $Y = 1.11$ was employed. Within this pivot range $(pc-1) < O(1)$, which means that the smaller of the two proposed collection sizes (61) must be employed. Now define a parameter $X = exp[\{(m+2)^{-1} - h_t\}/(m+1)]$ where $h_t = log\left(\varepsilon_t/\sqrt{2^m}\right)/log(t)$. This parameter enables one to subdivide the interval $[a+b, Ya]$ into a series of regularly increasing pivot 'Blocks', numbered $p = 1,2,\ldots L$. Each block consist of $\lfloor X^p t^{1/(m+2)}\rfloor$ distinct pivots $\alpha_E$ points and associated with each Block is a fixed collection size $M_{t,p} = NINT_O(M_t)$. This means within each Block, pivot points are separated from their nearest neighbours by a factor of $2M_{t,p}+2$. The boundaries of Block $p$ are defined by the interval $\left[\alpha_{E,p-1}\alpha_{E,p}\right]$, where

$$\alpha_{E,0} = INT_O(a) + 2 + b + M_{t,1} = INT_O(a) + 2 + 2\lfloor t^{1/(m+2)}\rfloor + M_{t,1},$$

$$\alpha_{E,p} = \alpha_{E,p-1} + 2\lfloor X^p t^{1/(m+1)}\rfloor\left(M_{t,p}+1\right), \quad p = 1,2,\ldots L$$

$$M_{t,1} = NINT_O[M_t(pc,m)], \quad M_{t,p} = NINT_O[M_t(pc,m)], \qquad p > 1.$$

(130a, b, c)

Within the Block interval $\left[\alpha_{E,p-1}\alpha_{E,p}\right]$ itself, the pivot integers are sited at the values $\alpha_{E,p-1} + \left(M_{t,p}+1\right)$, $\alpha_{E,p-1} + 3\left(M_{t,p}+1\right)$, $\alpha_{E,p-1} + 5\left(M_{t,p}+1\right)$, etc. It is helpful to rewrite the collection size (61) in terms of the pivot integers, using (13), to give

$$M_t(pc(\alpha_E),m) = \left[\left(\frac{\varepsilon_t}{\pi}\right)\frac{2^{3(3-m)/2}}{\sqrt{a}}\right]^{\frac{1}{(m+1)}} \times (\alpha_E - a)^{\frac{(2m-1)}{2(m+1)}} = \mathrm{C} \times (\alpha_E - a)^{\mu}. \qquad (131)$$

Then successively substituting (131) into (130b) starting from $p = 1$ up to $p = L$, one finds that

$$\alpha_{E,L} \gtrsim \left\{2\mathrm{C}t^{1/(m+2)}\right\}^{\sum_{p=1}^{L}\mu^{p-1}} \times X^{\sum_{p=1}^{L}p\mu^{L-p}}, \qquad (132)$$

where $\mu = (2m-1)/2(m+1)$. Since $\mu < 1$, the sums forming the indices of (132) converge for large $L$ and can be approximated by

$$\sum_{p=1}^{L}\mu^{p-1} \approx \frac{1}{1-\mu} = \frac{2(m+1)}{3}, \quad \sum_{p=1}^{L}p\mu^{L-p} \approx \frac{L}{1-\mu} - \frac{\mu}{(1-\mu)^2} \approx \frac{2(m+1)L}{3}, \qquad (133)$$

provided $L \gg m$. By the time the calculation has reached the $L \approx log(t)$ block, the pivot integers will have reached the value

$$\alpha_{E,log(t)} \gtrsim \left\{2\mathrm{C}t^{1/(m+2)}\right\}^{\frac{2(m+1)}{3}} X^{2(m+1)log(t)/3} \approx$$

$$2^{2(m+1)/3}\left[\left(\frac{\varepsilon_t}{\pi}\right)\frac{2^{3(3-m)/2}}{\sqrt{a}}\right]^{\frac{2}{3}} t^{\frac{2(m+1)}{3(m+2)}} \times exp\left[\frac{2log(t)}{3(m+2)} - \frac{2log\left(\varepsilon_t/\sqrt{2^m}\right)}{3}\right] = 2^{5/3}a. \qquad (134)$$

Hence the condition $\alpha_{E,L} > Ya$ is guaranteed to be met after $L \approx log(t)$ blocks. This value represents an upper bound, and in practice the value $Ya$ will be exceeded after a smaller number of blocks. So, for example, the $m = 3$ computations summarised in the next sub-section typically only require $L \approx log(t)/2$ blocks to reach this point. To take this discrepancy into account, it is useful to define a parameter $\lambda \in (0,1]$ such that $L = \lambda log(t)$. The number of operations required to reach (134) will be given by

$$2\sum_{p=1}^{L}\left\lfloor X^p t^{\frac{1}{(m+2)}}\right\rfloor \times GSum\ ops_{s;c} < log(t) \times X^{\lambda log(t)} t^{\frac{1}{(m+2)}} \times [2K + 2\Omega^* log(max(M_t^*)/K)]$$

$$\leq 2\left(\frac{\sqrt{2^m}}{\varepsilon_t}\right)^{\frac{\lambda}{(m+1)}} t^{\frac{(m+1+\lambda)}{[(m+1)(m+2)]}} log(t)\left[K + \frac{(m-1)\Omega^*}{2(m+1)}log(t) - \Omega^* log(K) + \varphi_1\right] ops_{s;c}. \qquad (135)$$

Employing collection size (27) bounds the constant $|\varphi_1| < \Omega^* \left|log\left\{Y\varepsilon_t^{1/(m+1)}/5\right\}\right|$.

The second stage of the procedure covers those pivot integers lying in the range $\alpha_E \in [Ya, \alpha_C(m)]$, which forms the main bulk of the calculation. Since in this region $(pc-1) \geq O(1)$ one can safely substitute (27) for the collection size. Once again, the regime is sub-divided into pivot Blocks numbered $p = L+1, L+2, \dots, L_C$, but here the Block size remains fixed at a value of $\lfloor X^L t^{1/(m+2)} \rfloor$. As in (131) above, it useful to rewrite the collection size (27) directly in terms of pivot integers, which gives

$$M_t(pc(\alpha_E), m) = \frac{\varepsilon_t^{1/(m+1)} t^{(m-1)/2(m+1)}}{2^{2(m-1)/(m+1)}\sqrt{\pi}}\left(\frac{2\alpha_E}{a} - \frac{a}{2\alpha_E} + O\left(\frac{a^3}{\alpha_E^3}\right)\right). \qquad (136)$$

Now by the end of Block $L_C$ the pivot integers will have reached values of

$$\alpha_{E,L+L_C} = \alpha_{E,L} + 2\lfloor X^L t^{1/(m+2)} \rfloor \sum_{p=1}^{L_C} (M_{t,L+p} + 1). \qquad (137)$$

Starting with $p = 1$ one has

$$M_{t,1}(pc(\alpha_E), m) \approx \frac{\varepsilon_t^{1/(m+1)} t^{(m-1)/2(m+1)}}{2^{2(m-1)/(m+1)}\sqrt{\pi}} 2Y\left(1 - \frac{1}{4Y^2}\right),$$

$$\Rightarrow \alpha_{E,L+1} \approx Ya + 4\left(\frac{\sqrt{2^m}}{\varepsilon_t}\right)^{\frac{\lambda}{(m+1)}} t^{\frac{(m+1+\lambda)}{[(m+1)(m+2)]}} \times \frac{\varepsilon_t^{1/(m+1)} t^{(m-1)/2(m+1)}}{2^{2(m-1)/(m+1)}\sqrt{\pi}} Y\left(1 - \frac{1}{4Y^2}\right)$$

$$= Ya\left[1 + \left\{1 - \frac{1}{4Y^2}\right\}\left\{\frac{\varepsilon_t^{(1-\lambda)}}{2^{[(3-\lambda)m-5]/2} t^{(1-\lambda)/(m+2)}}\right\}^{\frac{1}{(m+1)}}\right] = YaF_1. \qquad (138)$$

Successively repeating this process for $p = 2,3, \dots L_C$, gives an estimate for $\alpha_{E,L+L_C} \approx YaF_1F_2 \dots F_{L_C}$, where

$$F_p = \left[1 + \left\{1 - \frac{1}{4Y^2 F_{1\dots}F_{p-1}}\right\}\left\{\frac{\varepsilon_t^{(1-\lambda)}}{2^{[(3-\lambda)m-5]/2} t^{(1-\lambda)/(m+2)}}\right\}^{\frac{1}{(m+1)}}\right]. \qquad (139)$$

Now by the end of Block $L_C$, the pivot integers need to have reached the cut-off integer, i.e. $\alpha_{E,L+L_C} \sim \alpha_C(m)$ given by (125). Since $F_1 < F_2 < \cdots < F_{L_C}$, this means $\alpha_{E,L+L_C} > YaF_1^{L_C} = Yat^{\rho log(F_1)}$, after setting $L_C = \rho log(t)$ where $\rho > 0$ is parameter to be determined. To be certain of reaching the cut-off integer by the end of Block $L_C$ one must have

$$YaF_1^{L_C} = Yat^{\rho log(F_1)} \sim \alpha_C(m) \approx \frac{(\varepsilon_t t^m)^{1/(m+1)}}{\sqrt{\pi} log(t)},$$

$$\Rightarrow \qquad t^{\rho log(F_1)} \approx \left[\frac{\varepsilon_t^{1/(m+1)}}{Y\sqrt{8}}\right] t^{(m-1)/2(m+1) - log(log(t))/log(t)},$$

$$\Rightarrow \rho \approx \frac{(m-1)}{2(m+1)log(F_1)} + \frac{1}{log(F_1)log(t)}\left\{log\left[\frac{\varepsilon_t^{1/(m+1)}}{Y\sqrt{8}}\right] - log\big(log(t)\big)\right\}. \tag{140}$$

Hence the number of Blocks needed to reach the cut-off point will be

$$L_C = \rho log(t) \approx \left\{\frac{(m-1)}{2(m+1)log(F_1)}\right\} log(t). \tag{141}$$

With $\lambda = 1$, its theoretical maximum, (141) gives a limiting value for $L_C \leq 3log(t)/2$ for large $m$. To be certain that the parameter $\rho$ remains of $O(1)$ for practical computations in which $\lambda < 1$, one must impose a condition on the factor $F_1$, to ensure it does not get too close to unity. In practice this means the term $t^{(1-\lambda)/[(m+1)(m+2)]}$ appearing in (138) cannot be allowed to get too large. Theoretically this can be achieved by specifying that $m > O\big(\sqrt{log(t)}\big)$ in the limit $t \to \infty$. However, for the $m = 3$ scheme, which is of primary interest here, this will not be an issue until $t \gg 10^{40}$, well beyond the range of anything currently computationally feasible. For the sample computations involving practical values of $t$ discussed in the next sub-section, all values of $\rho$ fell within a range $1.75 - 2.30$. It follows that the number of operations required to reach Block $L_C$ is given by

$$2L_C\left[X^L t^{\frac{1}{(m+2)}}\right] \times GSum\, ops_{s;c} < \rho log(t) \times X^{\lambda log(t)} t^{\frac{1}{(m+2)}} \times [2K + 2\Omega^* log(max(M_t^*)/K)]$$

$$\leq 2\rho\left(\frac{\sqrt{2^m}}{\varepsilon_t}\right)^{\frac{\lambda}{(m+1)}} t^{\frac{(m+1+\lambda)}{[(m+1)(m+2)]}} log(t)\left[K + \frac{(m-1)\Omega^*}{(m+1)} log(t) - \Omega^* log\big(K log(t)\big) + \varphi_2\right] ops_{s;c}. \tag{142}$$

With collection size (27), the constant $|\varphi_2| < \Omega^*\left|log\left\{\varepsilon_t^{2/(m+1)}/40\right\}\right|$. Note too that in both estimates (135) and (142) the cut-off value $K$ is treated as a fixed constant, when in the limit as $t \to \infty$ the reductive bound (119) will manifest itself.

The sum of counts (135) and (142) gives a bound on the number of operations necessary to compute $ZP(t,m)$. In theory, the maximum number of operations required occurs when $\lambda = 1$, but in practice $\lambda$ is usually close to a half. As one can see from these estimates, completion of the second stage of the procedure typically requires something close to $2\rho$ more operations compared to the first stage. The final stage of procedure for the estimation of $Z(t)$ utilising the hybrid formula (126) involves the calculation of the sum of the first $N_C(m)$ terms of the main sum of the Riemann-Siegel formula. Obviously, this amounts to just

$$N_C(m) ops_{s;c} \approx \frac{(t/\varepsilon_t)^{1/(m+1)} log(t)}{\sqrt{\pi}} ops_{s;c}, \tag{143}$$

somewhat below the operational count required for $ZP(t,m)$.

4.3 *Sample Computations*

Thanks to the invaluable web-based database [9] of sample calculations of the Hardy function, it is possible to compare the performance of the $m = 3$ hybrid summation scheme discussed above, with the results obtained by utilising the $O\big(t^{1/3+o(1)}\big)$ operational algorithm of [2, 7]. Table 4 displays a summary of the accuracy of new hybrid scheme across one thousand Hardy function evaluations for values of $t \in [5.01 - 15.00] + 10^{24}$. For these relatively modest $t$ values the calculations can be performed in a reasonable time on a PC workstation, using a sequential version of the code. Pre-set values of $K = 400$ and $\varepsilon_t = 0.005$ were employed. The table shows the average absolute and relative errors arising from these one thousand function evaluations, when compared to the equivalent results using the algorithms [2, 7]. (The relative errors statistics are restricted to those values of $t$ for which $|Z(t)| > 1/2$, to avoid excessive values obtained when dividing numbers near zero.) Overall, the agreement is excellent. Typically, the hybrid scheme yields an estimate within $\pm 0.01$ of the corresponding values listed in [9]. In the worst-case scenario, the error reaches a maximum of $-0.0198$. Within this computational range, $Z(t)$ exhibits 84 zeros, each bracketed by relatively small $|Z(t)| < 1/2$ values, of different sign. The largest discrepancy at a zero point occurs for the zero situated just above $t = 11.84 + 10^{24}$. The hybrid scheme calculations gives an estimate of its position at $t = 11.84845 + 10^{24}$, compared to a value $t = 11.84450 + 10^{24}$ predicted by [9]. This illustration is important because, unlike for the quadratic $m = 2$ case formulated in [14], there is no guarantee that the associated errors in the estimates of each of the constituent cubic sums in (62) will not adversely combine to impinge significantly on the estimate for $Z(t)$ itself, particularly for calculations near zeros. However, is clear from these results that the highly regimented correlation between such errors, which would be necessary to bring about such an impingement, does not occur in practice. Between the zeta zeros, i.e. $|Z(t)| > 1/2$, both the average and the worst-case scenario of the relative error lie well below the upper bound $\varepsilon_t$ prescribed for these calculations.

The total operational count for the hybrid scheme is given by the combination of the estimates (135, 142-143). Since $\Omega^* \approx 31$, this means that $\Omega^* log(t)/2$ factor appearing in (142) will comfortably exceed a value of $K = 400$ by at least a factor of two for $t \geq 10^{23}$, and that (142) is the dominant, first order, contribution to the total operational count. A series of computations were carried out to test both the veracity of this estimate and provide an illustration of the superior computational speed of the hybrid scheme vis-à-vis the $O\left(t^{1/3+o(1)}\right)$ algorithm of [7]. Hybrid computations were made for a series of $g = log_{10}(t)$ values with $g > 23$, for fixed parameter values $K = 400$ and $\varepsilon_t = 0.01$, utilising a sequential version of the code. The corresponding run times were then compared to an identical computation when $t = 10^{23}$. Since run time acts as a proxy for operational count, the respective run times should increase (to first order) by an operation factor $OF1(g)$ given by

$$OF1(g) = \left(\frac{g}{23}\right)^2 10^{(g-23)/4}, \tag{144}$$

if the operational count estimate (142) is accurate. Utilising the $O\left(t^{1/3+o(1)}\right)$ algorithm for computations over the same range of $t$ values, one would expect the run times to scale upwards by a factor

$$OF2(g) = \left(\frac{g}{23}\right)^{\kappa_{10}} 10^{(g-23)/3}, \tag{145}$$

devised from the operational count estimate of Theorem 2.3 of [7]. No estimate for the absolute constant $\kappa_{10}$ is put forward in [7], but it cannot be any less than two and is probably a little larger. Figure 3 shows a plot of these two operational factors (144-5), the latter using $Min(\kappa_{10}) = 2$, and the corresponding values devised from the series of hybrid computation run times. For values of $g \in [23, 25]$ the observed increase in computational run times closely matches that predicted by $OF1(g)$, with errors small enough to be accounted for by corrections arising from the second order terms in (135, 142). However, beyond $g = 25$ the run times are clearly starting to drift upwards, away from $OF1(g)$, when they ought to converge as the second order terms become less significant. This drift is a manifestation of the reductive bound issue discussed in Section 3.4.2. As $t$ gets larger, the respective collection sizes $M_t(pc, 3)$ increase, which results in longer progenitor sums. As the length of the progenitor sums increases, it becomes more difficult to reduce these down to kernel sums of arbitrary fixed length $K$ before the reductive bound (119) is reached. In practical terms, this means that the kernel sums computed by algorithm *MGS* become lengthier and take longer to calculate, slowing down the speed of the overall computations for $Z(t)$. Hence the drift away from the predicted behaviour (144) for $g > 25$ seen in Figure 3. Nevertheless, as one can see, all the run times remain *considerably faster* than the (minimal) prediction (145) derived from the $O\left(t^{1/3+o(1)}\right)$ algorithm. An optimal, least squares fit of the run time data gives an exponent estimate of $t^{0.274}$ across this range of $t$. (And although this exponent will increase with $t$, estimates indicate it should not exceed a value 0.3 until $t > 10^{35}$.) Indeed, since the difference between the operational factors (144) and (145) increases exponentially with $g$, the hybrid scheme should extend its superiority over the $O\left(t^{1/3+o(1)}\right)$ algorithm in terms of computational efficiency. So, for example, when $g = 25$ the actual hybrid run time is 28% faster than $OF2(g)$ prediction, whilst for $g = 26.247$ (the longest run shown on the figure) the corresponding saving increases to 34%.

Table 5 shows a similar comparison of the accuracy of the new hybrid scheme for a selection of $t$ values where $Z(t)$ exhibits an unusually large local maximum/minimum, compared to other such points in its vicinity. (Since $Z(t)$ is unbounded, much greater values of $|Z(t)|$ certainly exist, although for larger, computationally unfeasible $t$ values.) Details of the methodology utilised to identify narrow intervals along the critical line where such points are likely (but not guaranteed) to occur, are discussed in [9]. The latter reference also provides details of many examples where $|Z(t)| > 10^3$, whilst [23] highlights some calculations of *especially large* values of $|Z(t)| > 10^4$, utilising the $O\left(t^{1/3+o(1)}\right)$ algorithm. The large $|Z(t)|$ computations shown in Table 5 are for $t > 10^{28}$, covering the highest possible range of computationally feasible $t$ values.

For these particular computations, a parallelised version of the hybrid summation code was developed for implementation on the ARCHER2 UK National Supercomputing Service (http://www.archer2.ac.uk) hosted by the University of Edinburgh. Open-source variants of this hybrid code are now freely available from the repository [15], along with documentation describing its implementation and some sample output files. ARCHER2 is a HPE Cray EX supercomputing system with a total of 5860 nodes, each consisting of two 2.25GHz, 64 core AMD EPYC7742 processors with access to 256GB of memory (https://www.archer2.ac.uk/user-guide/hardware/). Most of the computations listed in Table 5 were run using 4-8 nodes and required between 3-12 hours of CPU time. However, the highest calculations, when $t \geq 10^{32}$, were carried out during one of the periodic ARCHER2 Capability Days sessions, which give users the opportunity to run large scale jobs, at no budgetary cost, using many nodes. The $t \geq 10^{32}$ calculations highlighted in the Table 5 employed 512-1024 nodes and typically required between 1-3 hours of CPU time.

Looking at the results shown in Table 5, one can see the agreement between the hybrid estimates and the corresponding values published in [9, 23] is excellent, with all the relative errors lying far below the prescribed value of $\varepsilon_t = 0.005$. The absolute errors are also very low, which is noteworthy since such errors tend to reach a maximum in and around points where $|Z(t)|$ is large. Of course, if more accurate estimates are desired, then one can simply reduce the user prescribed value of $\varepsilon_t$, at the expense of a marginal increase in run times. One of the chief difficulties encountered for computations where $t > 10^{28}$ is the potential for catastrophic round-off error, resulting from the fact that Fortran double precision variables can only represent numbers to an accuracy of some thirty-three significant figures. In order to circumvent this issue, the main Fortran codes in the repository [15] utilise certain C subroutines from the PariGP computer library [20], accurate to many thousands of digits, for those particular calculations where round-off error is problematic. Indeed, it is anticipated that the repository code should be capable of carrying out computations for $t \sim 10^{40}$ or even a little beyond, before the need for quadruple precision reals becomes imperative. However, although feasible, calculations on this scale would exhaust the supercomputer budgetary resource allocation very generously awarded in support of this work.

## 5. Conclusions

The main focus of this work has been on the development of a new, novel methodology for the computation of the Hardy function $Z(t)$ for high values of $t$ along the critical line. Currently the fastest state of the art algorithms, developed by [2, 7], for such computations have an operational complexity $O(t^{1/3+o(1)})$ in time, or a slightly superior $O(t^{4/13+o(1)})$ algorithm in time and space [7]. These algorithms depend principally upon the observation that the main sum terms of the Riemann-Siegel formula (3) can be amalgamated into blocks of quadratic Gauss sums $S_N(\Phi_1, \Phi_2)$ of collection size $N = O(t^{1/6})$, which can themselves be rapidly computed in something like $O(log(N))$ operations using the recursive methods of [7, 14, 18]. Theoretically it should be possible to extrapolate this idea and formulate (3) in terms of longer, higher order, generalised Gauss sums. However, in general the computational utility of such a procedure would be futile, because such sums do not satisfy the quadratic reciprocity property underlying the recursive methodology utilised for the estimation of $S_N(\Phi_1, \Phi_2)$. This observation leads to the proposition posted on [9] that there exists a "natural barrier" to the complexity exponent for the computation of $Z(t)$ at, or just below, the value of $1/3$ and obtaining, for instance, a $O(t^{1/4+o(1)})$ algorithm will be "challenging as it requires new ideas".

The work in this paper describes the formulation of a new algorithm, with an operational complexity $O\left((t/\varepsilon_t)^{[0.25,\,0.3]}(log(t))^{2+o(1)}\right) ops_{s;c}$ for $t \in [10^{23-35}]$, corresponding to the $m = 3$ summation scheme (q.v. eqs. 135 & 142, assuming the worst-case scenario when $\lambda = 1$). The algorithm is then implemented for the practical computations described in Section 4. The basis of the development of this scheme is the exact formulation of the Hardy function in terms of a sum of Kummer's functions (4), the starting point of the "new ideas" requirement specified above. The extrapolation, carried out in Section 2, of the approximation methods initially developed in [13, 14] to estimate these specific Kummer's functions, yields a representation (62) for $Z(t)$ in terms of arbitrary, $m^{th}$-order Gauss sums, of ever-increasing collection size $M_t$. The asymptotic structure of the coefficients of these particular Gauss sums makes them amenable for rapid computation in $O(2K + 2\Omega^* log(\, M_t/2K)) ops_{s;c}$, utilising *algorithm MGS* devised in Section 3. Comparisons of the algorithm's performance *vis-à-vis* those developed in [2, 7], clearly show it is superior in operational speed, allowing for computations at $t$ values up to $10^4$ times larger than previously attempted. Open-source code that achieves this is now available from the repository [15]. True the reductive bound (119) on the cut-off integer $K$ will, ultimately, impact and increase the number of operations required to evaluate $Z(t)$ for arbitrarily large values of $t$. But in practice this will not set in until at least $t > 10^{40}$, well beyond anything currently computationally feasible.

Moving forward, the next step is to move towards the formulation of a genuine, practical $O(t^{1/4+o(1)})$ algorithm based around quartic $m = 4$ Gauss sums. Such a development would require an efficient means of truncating through the point in the iteration process when the reductive bound issue starts to manifest itself. The example computation shown in Table 3, illustrates how the increase in the (relative) size of the cubic sum coefficient $\Phi_3$ eventually prevents further truncation down to a kernel sum of an arbitrary length $K = O(10 - 100)$, beyond the bound (119). The extra iteration methodology briefly discussed at the end of section 3.4.2 shows promise. However, at this stage implementation of this idea for quartic sums (see [15] for a sample quartic code) yields only a faster, if slightly less accurate, version of the $O\left(t^{[0.25,\,0.3]}\right)$ cubic sum algorithm already developed. But to go much beyond this point and significantly further reduce kernel sum is currently problematic, without the introduction of a sub-division scheme which is likely to be both operationally expensive and less accurate. Be that as it may, this work still opens up a number of promising avenues of future research directed towards the development of yet more efficient algorithms for computation of the Hardy function.

**Acknowledgements:** The authors would like to acknowledge the help and advice of the ARCHER2 UK National Supercomputing Service User Support Team in facilitating the computational results presented in this paper. In particular, the help of Dr Mario Antonioletti (EPCC, The University of Edinburgh) is acknowledged. Establishment of the GitHub computational repository [15] was funded under the EPSRC/ARCHER2-eCSE grant scheme (grant no. 11-7).

## Figure Captions

Figure 1. Contour integration path utilised for the estimation of integral $I(A, B^{\pm}, C)$ given by equation (29). The original range of integration range $[-X, X]$ along the real axis (dot-dash arrow), must be rotated clockwise through an angle $\varphi = \pi/4$ as shown (dashed arrow), so that the contour passes through one of the two saddle points situated at $u_* \approx \mp e^{i3\pi/4} X/\sqrt{2}$ corresponding to $B^{\pm}$ respectively.

Figures 2a and 2b. Contour integration paths utilised to estimate integrals given by equations (71: Fig. 2a) and (82: Fig. 2b). In each case the path of integration along the real line (solid arrow) is substituted by the complex contour path traced out by the dashed arrows. Estimates for the constituent line segments forming the respective contour paths are derived in the text.

Figure 3. Illustration of the comparative performance of the hybrid scheme devised here and the $O\left(t^{1/3+o(1)}\right)$ algorithm of [7]. The two operational factor scales [144-5] illustrate how run times are predicted to increase with $g = log_{10}(t)$ beyond an (arbitrary) normalisation scale of $g = 23$. The dots represent corresponding operational factors derived from sample computational runs of the hybrid scheme.

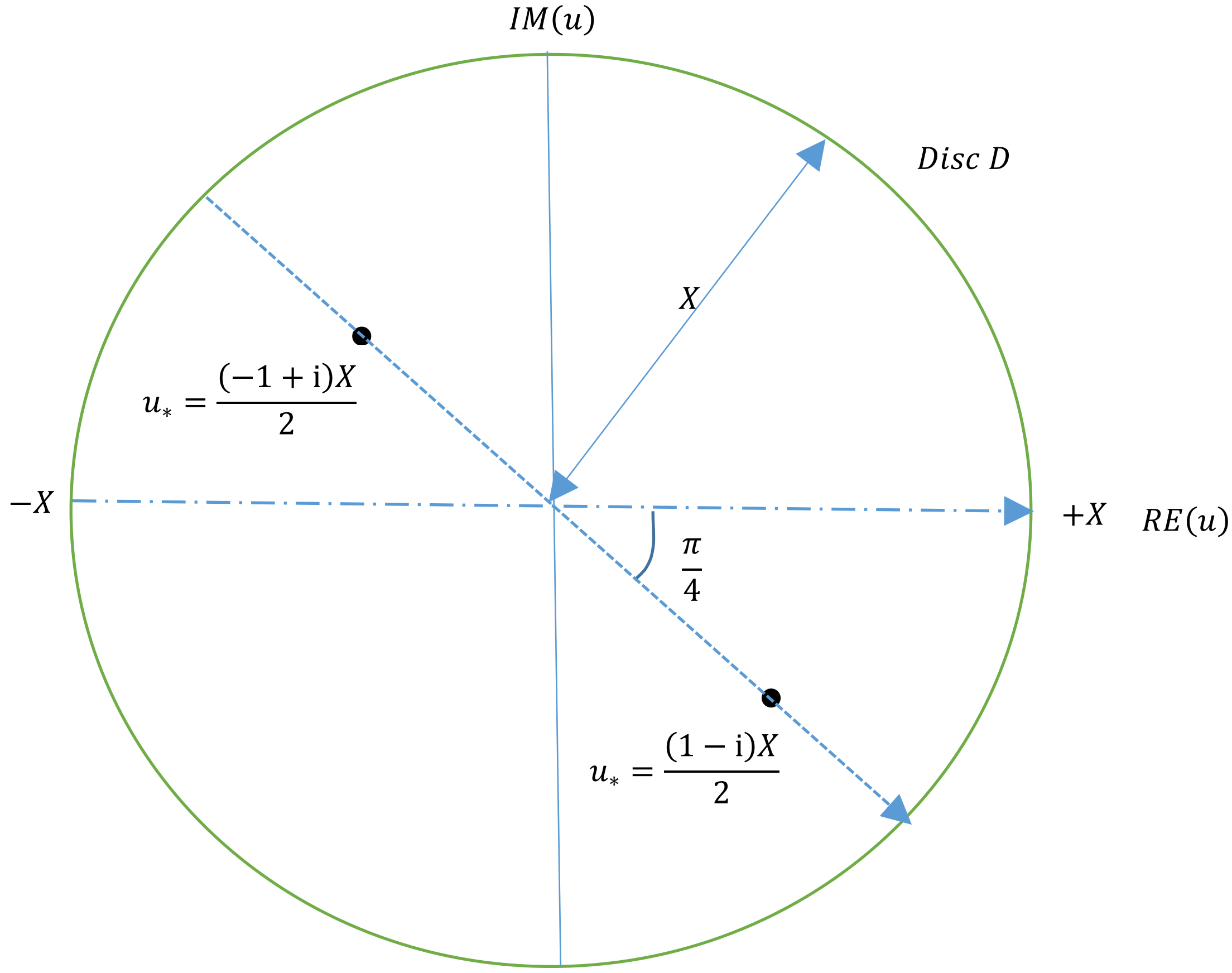
$IM(u)$
$Disc\ D$
$X$
$u_* = \frac{(-1+\mathrm{i})X}{2}$
$-X$
$+X$
$RE(u)$
$\frac{\pi}{4}$
$u_* = \frac{(1-\mathrm{i})X}{2}$

Figure 1

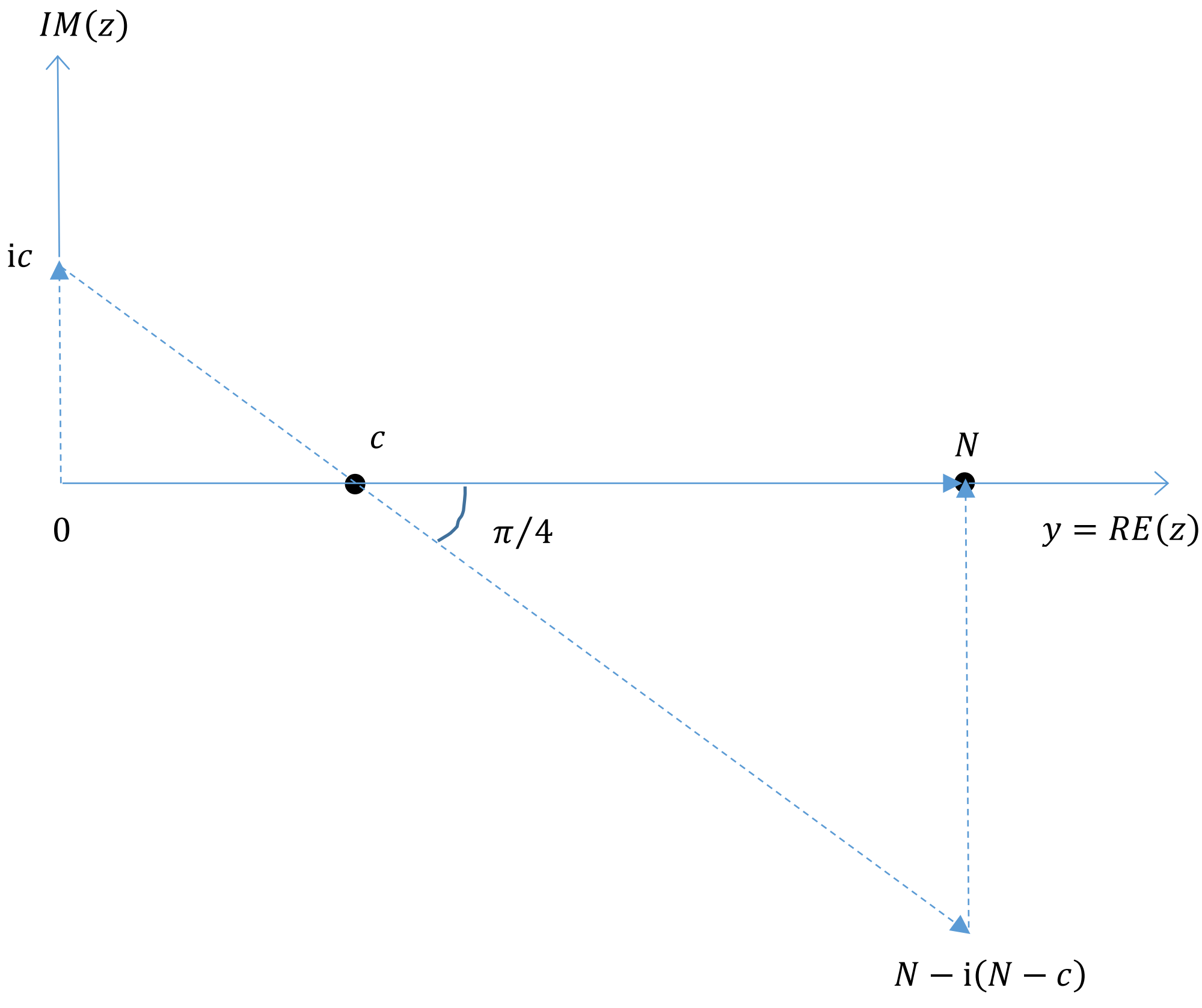

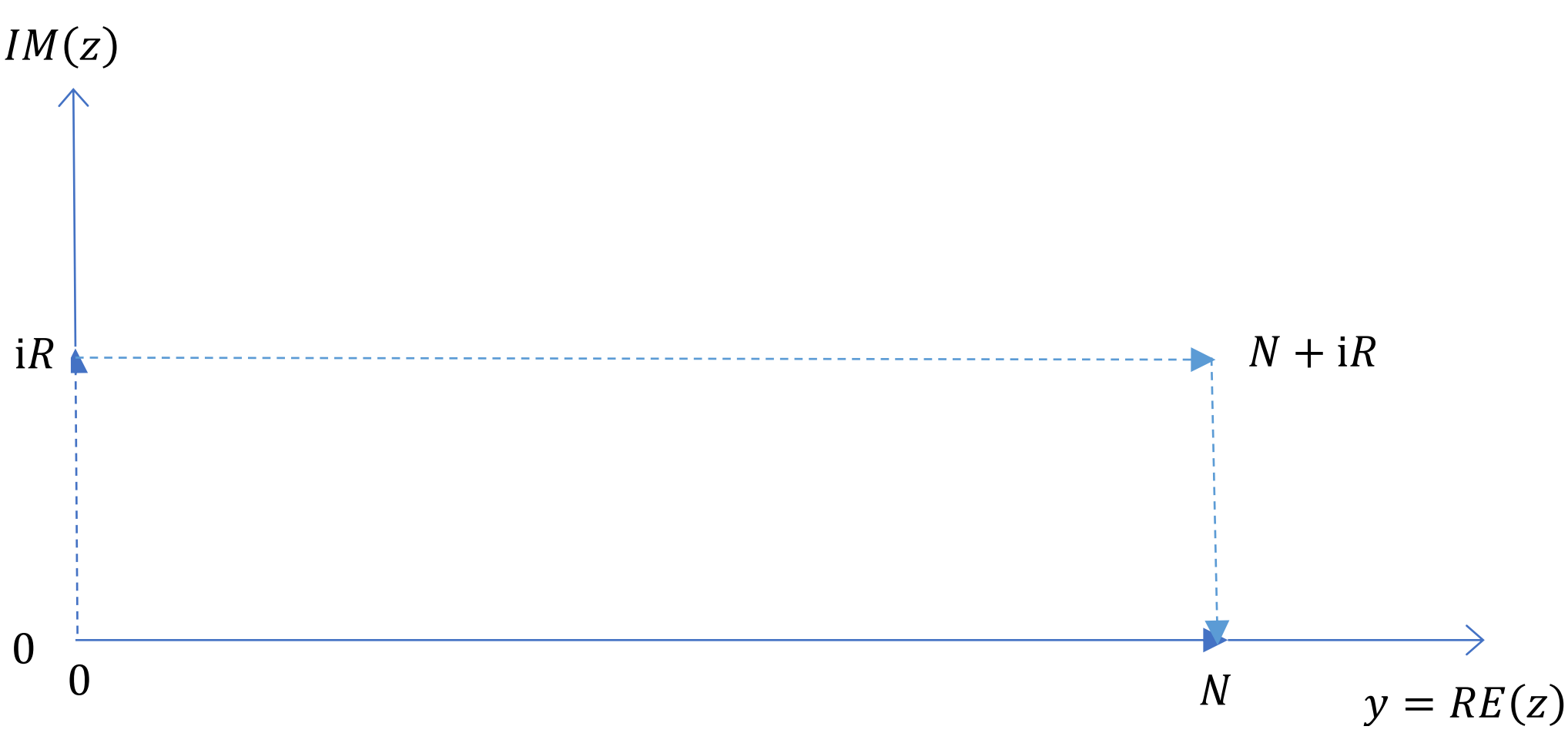

Figure 2a & 2b

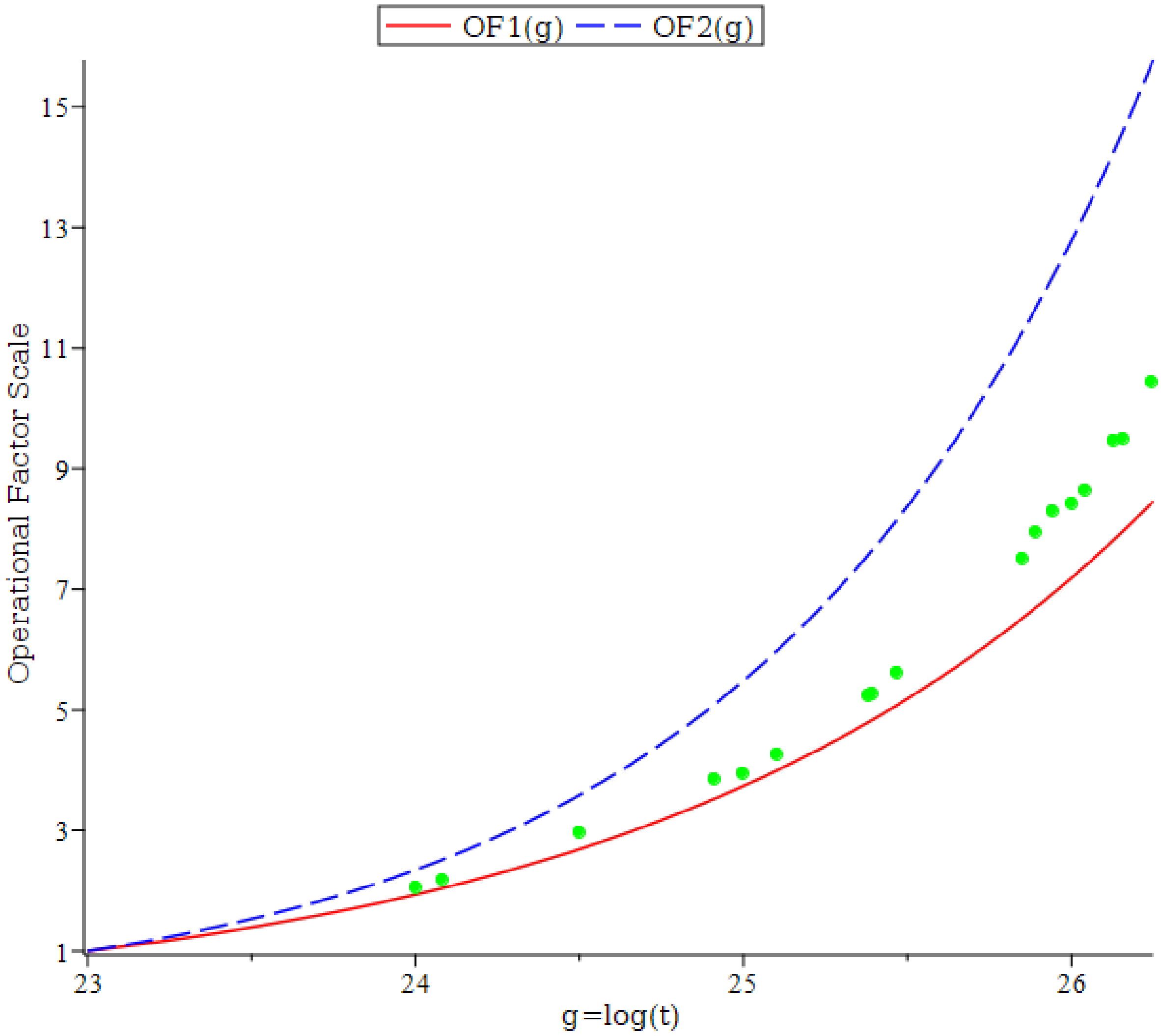
OF1(g)
OF2(g)
Operational Factor Scale
15
13
11
9
7
5
3
1
23
24
25
26
g=log(t)

Figure 3

## Mathematical Appendix

A list of the first nine coefficients $C_q(\Phi_{2,3,\ldots q+1})$ which form the $m^{th}$-order saddle point solutions $c_{m-1}$ satisfying (99-100), generated using the Maple algebraic software. They fix the size and scale of the remainder terms (104) that govern the iteration scheme implemented in *algorithm MGS* of section 3.3. They quickly become increasingly complex and there is no simple recursion scheme governing their behaviour.

In the expressions below $x = 2\Phi_2 \in (0, 1/2)$.

$$C_1(\Phi_2) = 1/2\,.$$

$$C_2(\Phi_2, \Phi_3) = -\Phi_3.$$

$$C_3(\Phi_2, \Phi_3, \Phi_4) = \frac{(9\Phi_3^2 - 2x\Phi_4)}{2}.$$

$$C_4(\Phi_2, \Phi_3, \Phi_4, \Phi_5) = -(27\Phi_3^3 - 12x\Phi_3\Phi_4 + x^2\Phi_5).$$

$$C_5(\Phi_2, \Phi_3, \Phi_4, \Phi_5, \Phi_6) = [189\Phi_3^4 - 126x\Phi_3^2\Phi_4 + x^2(15\Phi_3\Phi_5 + 8\Phi_4^2) - x^3\Phi_6].$$

$$C_6(\Phi_2, \Phi_3, \Phi_4, \Phi_5, \Phi_6, \Phi_7) = -\big[1458\Phi_3^5 - 1296x\Phi_3^3\Phi_4 + 12x^2\Phi_3(15\Phi_3\Phi_5 + 16\Phi_4^2) \\ -2x^3(9\Phi_3\Phi_6 + 10\Phi_4\Phi_5) + x^4\Phi_7].$$

$$C_7(\Phi_2, \Phi_3, \Phi_4, \Phi_5, \Phi_6, \Phi_7, \Phi_8) = [24057\Phi_3^6 - 26730x\Phi_3^4\Phi_4 + 810x^2(8\Phi_3^2\Phi_4^2 + 5\Phi_3^3\Phi_5) \\ -6x^3(81\Phi_3^3\Phi_6 + 180\Phi_3\Phi_4\Phi_5 + 32\Phi_4^3) + x^4\big(42\Phi_3\Phi_7 + 48\Phi_4\Phi_6 + 25\Phi_5^2\big) - 2x^5\Phi_8\big]/2.$$

$$C_8(\Phi_2, \Phi_3, \Phi_4, \Phi_5, \Phi_6, \Phi_7, \Phi_8, \Phi_9) = -[104247\Phi_3^7 - 138996x\Phi_3^5\Phi_4 \\ +1485x^2(32\Phi_3^3\Phi_4^2 + 15\Phi_3^4\Phi_5) - 110x^3(27\Phi_3^3\Phi_6 + 90\Phi_3^2\Phi_4\Phi_5 + 32\Phi_3\Phi_4^3) \\ +5x^4\big(63\Phi_3^2\Phi_7 + 144\Phi_3\Phi_4\Phi_6 + 75\Phi_3\Phi_5^2 + 80\Phi_4^2\Phi_5\big) - x^5(24\Phi_3\Phi_8 + 28\Phi_4\Phi_7 + 30\Phi_5\Phi_6) + x^6\Phi_9\big].$$

$$C_9(\Phi_2, \Phi_3, \Phi_4, \Phi_5, \Phi_6, \Phi_7, \Phi_8, \Phi_9, \Phi_{10}) = [938223\Phi_3^8 - 1459458x\Phi_3^5\Phi_4 \\ +81081x^2\Phi_3^4(3\Phi_3\Phi_5 + 8\Phi_4^2) - 1287x^3(27\Phi_3^4\Phi_6 + 120\Phi_3^3\Phi_4\Phi_5 + 64\Phi_3^2\Phi_4^3) \\ +11x^4\big(378\Phi_3^3\Phi_7 + 1296\Phi_3^2\Phi_4\Phi_6 + 675\Phi_3^2\Phi_5^2 + 1440\Phi_3\Phi_4^2\Phi_5 + 128\Phi_4^4\big) \\ - 22x^5\big(18\Phi_3^2\Phi_8 + 42\Phi_3\Phi_4\Phi_7 + 45\Phi_3\Phi_5\Phi_6 + 24\Phi_4^2\Phi_6 + 25\Phi_4\Phi_5^2\big) + x^6(27\Phi_3\Phi_9 \\ + 32\Phi_4\Phi_8 + 35\Phi_5\Phi_7 + 18\Phi_6^2) - x^7\Phi_{10}\big].$$

## Tables

| $l$ | $L_l$ | $m_l$ | $Generalized\ Gauss\ sum\ factors\ \Phi_{0-m_l}$ | $s_l$ |
|---|---|---|---|---|
| 1 | 71487 | 3 | $\Phi_0 = 0.0$<br>$\Phi_1 = 6.824248433 \times 10^{-2}$<br>$\Phi_2 = 1.177409528 \times 10^{-1}$<br>$\Phi_3 = -9.00195855 \times 10^{-14}$ | $+1$ |
| 2 | 16833 | 3 | $\Phi_0 = 9.88831106 \times 10^{-3}$<br>$\Phi_1 = -2.89799269 \times 10^{-1}$<br>$\Phi_2 = 1.233053918 \times 10^{-1}$<br>$\Phi_3 = 6.893890544 \times 10^{-12}$ | $-1$ |
| 3 | 4150 | 3 | $\Phi_0 = 1.70275637 \times 10^{-1}$<br>$\Phi_1 = 1.75128133 \times 10^{-1}$<br>$\Phi_2 = 2.748635921 \times 10^{-2}$<br>$\Phi_3 = -4.59650983 \times 10^{-10}$ | $-1$ |
| 4 | 228 | 4 | $\Phi_0 = 2.789552919 \times 10^{-1}$<br>$\Phi_1 = -1.85727835 \times 10^{-1}$<br>$\Phi_2 = 9.541922607 \times 10^{-2}$<br>$\Phi_3 = 2766854649 \times 10^{-6}$<br>$\Phi_4 = 1.893794426 \times 10^{-12}$ | $-1$ |
| 5 | 43 | 5 | $\Phi_0 = 9.037448961 \times 10^{-2}$<br>$\Phi_1 = 4.731790001 \times 10^{-1}$<br>$\Phi_2 = 1.197952734 \times 10^{-1}$<br>$\Phi_3 = -3.97996774 \times 10^{-4}$<br>$\Phi_4 = 1.346150552 \times 10^{-7}$<br>$\Phi_5 = -6.073686206 \times 10^{-11}$ | $-1$ |
| $l$ | Estimate of $S_{L_l}^{s_l}(\Phi_{1..m_l,l})$ | Exact Value $S_{L_l}^{s_l}(\Phi_{1..m_l,l})$ | $\lvert error \rvert$ | Relative error |
| 5 | As exact value | $-3.23983 + 4.21204i$ | 0.0 | 0.0 |
| 4 | $-5.59315 - 12.22387i$ | $-5.55197 - 12.16960i$ | 0.0681 | 0.0051 |
| 3 | $29.66146 + 49.50084i$ | $29.58776 + 49.22094i$ | 0.2894 | 0.0050 |
| 2 | $29.50811 - 113.55207i$ | $29.52304 - 112.94740i$ | 0.6048 | 0.0052 |
| 1 | $-108.45643 - 217.59156i$ | $-107.59957 - 216.68530i$ | 1.2472 | 0.0052 |

Table 1. Example illustration of the reduction and computation of a $3^{rd}$- order cubic Gaussian sum, satisfying the specific constraints (63), achieved using *algorithm (MGS)*. Parameter values: $t = 10^{28}$, pivot integer $ae = 163614984609380.0 \equiv\ pc = 1.566726 \ldots$, $K = 200$, and $P = 4$. Since $pc$ is close to one, the smaller collection size (61) is employed with $\varepsilon_t = 0.02$. The initial sum is reduced to a $5^{th}$- order kernel sum of length $N = 43$, and the iterations reversed to produce the estimate.

| $l$ | $L_l$ | $m_l$ | $Generalized\ Gauss\ sum\ factors\ \Phi_{0-m_l}$ | $s_l$ |
|---|---|---|---|---|
| 1 | 2100836 | 3 | $\Phi_0 = 0.0$<br>$\Phi_1 = 4.577951577 \times 10^{-1}$<br>$\Phi_2 = 4.149941700 \times 10^{-2}$<br>$\Phi_3 = 2.251742112 \times 10^{-16}$ | $-1$ |
| 2 | 174367 | 3 | $\Phi_0 = 1.262526209 \times 10^{0}$<br>$\Phi_1 = 4.843182053 \times 10^{-1}$<br>$\Phi_2 = 2.418101414 \times 10^{-2}$<br>$\Phi_3 = -3.93824450 \times 10^{-13}$ | $-1$ |
| 3 | 8433 | 4 | $\Phi_0 = 2.42508567 \times 10^{0}$<br>$\Phi_1 = -4.85568797 \times 10^{-1}$<br>$\Phi_2 = 1.613103826 \times 10^{-1}$<br>$\Phi_3 = -3.481682424 \times 10^{-9}$<br>$\Phi_4 = -2.63812477 \times 10^{-18}$ | $+1$ |
| 4 | 2719 | 4 | $\Phi_0 = 3.654090031 \times 10^{-1}$<br>$\Phi_1 = 5.076154596 \times 10^{-3}$<br>$\Phi_2 = 4.980739785 \times 10^{-2}$<br>$\Phi_3 = 1.036841408 \times 10^{-7}$<br>$\Phi_4 = 1.585086768 \times 10^{-14}$ | $-1$ |
| 5 | 273 | 6 | $\Phi_0 = 1.293349449 \times 10^{-4}$<br>$\Phi_1 = -5.09578458 \times 10^{-2}$<br>$\Phi_2 = 1.933629021 \times 10^{-2}$<br>$\Phi_3 = -1.048917155 \times 10^{-4}$<br>$\Phi_4 = 4.770974356 \times 10^{-9}$<br>$\Phi_5 = -2.89012763 \times 10^{-13}$<br>$\Phi_6 = 2.042123848 \times 10^{-17}$ | $-1$ |
| $l$ | Estimate of $S_{L_l}^{s_l}(\Phi_{1..m_l,l})$ | Exact Value $S_{L_l}^{s_l}(\Phi_{1..m_l,l})$ | $\|error\|$ | Relative error |
| 5 | As exact value | 1.57504 + 14.86040i | 0.0 | 0.0 |
| 4 | 38.69952 − 29.34762i | 38.73012 − 29.23127i | 0.1203 | 0.0025 |
| 3 | 55.47567 + 65.28449i | 55.34872 + 65.34816i | 0.1420 | 0.0017 |
| 2 | −361.23492 − 148.67401i | −361.27895 − 148.02249i | 0.6530 | 0.0017 |
| 1 | −421.06045 + 1288.25289i | −422.86094 + 1286.81109i | 2.3066 | 0.0017 |

Table 2. Further illustration of the reduction and computation of a $3^{rd}$- order cubic Gaussian sum achieved using *algorithm (MGS)*. Parameter values: $t = 10^{28}$, pivot integer $ae = 310304270951266.0 \equiv pc = 13.048362 \ldots$, and others as for Table 1. Since $pc$ is relatively large, collection size (27) is used. The initial sum is reduced to a $6^{th}$- order kernel sum of length $N = 273$, and the iterations reversed to produce the estimate.

| $l$ | $L_l$ | $m_l$ | $Gauss\ sum\ factors\ \Phi_{0-m_l}$ | $s_l$ |
|---|---|---|---|---|
| 1 | 3745088 | 4 | $\Phi_1 = 1.672394213 \times 10^{-1}$<br>$\Phi_2 = 1.177409528 \times 10^{-1}$<br>$\Phi_3 = 9.00195855 \times 10^{-14}$<br>$\Phi_4 = -1.69284335 \times 10^{-26}$ | $+1$ |
| 2 | 881904 | 4 | $\Phi_0 = 5.938677957 \times 10^{-2}$<br>$\Phi_1 = 2.89799269 \times 10^{-1}$<br>$\Phi_2 = 1.233053918 \times 10^{-1}$<br>$\Phi_3 = -6.89389054 \times 10^{-12}$<br>$\Phi_4 = 5.586689046 \times 10^{-23}$ | $-1$ |
| 3 | 217471 | 5 | $\Phi_0 = 1.70275637 \times 10^{-1}$<br>$\Phi_1 = -1.7512813 \times 10^{-1}$<br>$\Phi_2 = 2.748635921 \times 10^{-2}$<br>$\Phi_3 = 4.596509835 \times 10^{-10}$<br>$\Phi_4 = 2.193620923 \times 10^{-17}$<br>$\Phi_5 = 1.389215833 \times 10^{-28}$ | $-1$ |
| 4 | 12020 | 8 | $\Phi_0 = 2.789552919 \times 10^{-1}$<br>$\Phi_1 = 1.857278353 \times 10^{-1}$<br>$\Phi_2 = 9.541922607 \times 10^{-2}$<br>$\Phi_3 = -2.76685466 \times 10^{-6}$<br>$\Phi_4 = 1.869772964 \times 10^{-12}$<br>$\Phi_5 = -1.68472802 \times 10^{-18}$<br>$\Phi_6 = 1.770999255 \times 10^{-24}$<br>$\Phi_7 = -2.05165786 \times 10^{-30}$<br>$\Phi_8 = 2.541854919 \times 10^{-36}$ | $-1$ |
| 5 | 1108 | 8 | As above | $-1$ |
| $l$ | Estimate of $S_{L_l}^{s_l}(\Phi_{1..m_l,l})$ | Exact Value $S_{L_l}^{s_l}(\Phi_{1..m_l,l})$ | $\lvert error \rvert$ | Relative error |
| 4 | As exact value | $-23.66524 + 69.02917i$ | 0.0 | 0.0 |
| 3 | $-298.94781 - 84.02281i$ | $-298.40153 - 83.58318i$ | 0.7151 | 0.0023 |
| 2 | $-531.15710 + 331.88097i$ | $-530.33628 + 330.72472i$ | 1.4171 | 0.0023 |
| 1 | $-728.22660 + 1065.14613i$ | $-726.56506 + 1061.95502i$ | 3.5977 | 0.0028 |
| 5 | As exact value | $13.67259 - 28.45408i$ | 0.0 | 0.0 |
| 4 | $-23.62772 + 69.10024i$ | $-23.66524 + 69.02917i$ | 0.0823 | 0.0011 |
| 3 | $-298.69827 - 84.22845i$ | $-298.40153 - 83.58318i$ | 0.7102 | 0.0023 |
| 2 | $-530.56836 + 332.14394i$ | $-530.33628 + 330.72472i$ | 1.4376 | 0.0023 |
| 1 | $-726.89743 + 1065.15463i$ | $-726.56506 + 1061.95502i$ | 3.2167 | 0.0025 |

Table 3. Illustration of the reduction and computation of a $4^{th}$- order quartic Gaussian sum achieved using *algorithm (MGS)*. Parameter values: $t = 10^{28}$, pivot integer $ae = 163614984609380.0 \equiv pc = 1.566726 \ldots$, and other as for Table 1. Since $pc$ is near one, collection size (61) is used. Only a partial reduction of the initial sum is possible, as discussed in Section 3.4.2.

| Computational regime: $Z(t = 10^{24} + 5.00 + n \times 0.01)$ with $n \in [1,1000]$ | |
|---|---|
| Average absolute error | $7.53\times 10^{-3}$ |
| Worst Case: $Z(10^{24} + 14.21)$<br>[9] result : hybrid estimate | $1.98\times 10^{-2}$<br>10.2064716 : 10.1866861 |
| | |
| Average relative error (for $\|Z(t)\| > 1/2$) | $4.19\times 10^{-3}$ |
| Worst Case: $Z(10^{24} + 14.33)$<br>[9] result : hybrid estimate | $3.52\times 10^{-2}$<br>0.5151024 : 0.4969901 |
| Largest errors bracketing a zero:<br>$Z(10^{24} + 11.84)$, [9] result : hybrid estimate<br>$Z(10^{24} + 11.85)$, [9] result : hybrid estimate<br>$t_{zero}$ , [9] result : hybrid estimate | <br>$-0.0265649$ : $-0.0074452$<br>0.0217327 : 0.0408843<br>$10^{24} + 11.84450$ : $10^{24} + 11.84845$ |

Table 4. Illustrative error analysis in the calculations of the Hardy function when utilising the hybrid formulation methodology (with $K = 400$ and $\varepsilon_t = 0.005$), in comparison with the corresponding values listed in [9]. The table shows the average and worst case scenarios of the respective absolute and relative errors resulting from one thousand calculations of $Z(t)$ for $t \in [5.01, 15.00] \times 10^{24}$.

| $t$ value | $Z(t)$ Hybrid algorithm of sections 4.1-2. | $Z(t)$ $t^{1/3}$ algorithm of [2, 7]; $\pm 5 \times 10^{-7}$ | Relative error $\lvert\varepsilon\rvert$ |
|---|---|---|---|
| 2.059936512320112518074691006892E28 | 3803.86727 | 3803.86539 | $4.30 \times 10^{-7}$ |
| 3.177469531676391818363765436511E28 | −9549.91673 | −9549.88868 | $2.93 \times 10^{-6}$ |
| 4.670914185466097236850548903274E28 | 9587.87895 | 9587.85455 | $2.54 \times 10^{-6}$ |
| 1.0835660710102772561572011153186E29 | 7297.78564 | 7297.75658 | $3.98 \times 10^{-6}$ |
| 2.8928607671932530771838054905026E29 | 10907.57778 | 10907.55189 | $2.37 \times 10^{-6}$ |
| 5.5216641000993128888680863234651E29 | −13541.70579 | −13541.67744 | $2.09 \times 10^{-6}$ |
| 6.9815628897151991613594294046033E29 | 11187.59076 | 11187.56821 | $2.02 \times 10^{-6}$ |
| 8.0362572859234436312381421877820E29 | 10282.41019 | 10282.40085 | $9.13 \times 10^{-7}$ |
| 1.42060805696869950116950900347991E30 | 5045.59913 | 5045.59024 | $1.76 \times 10^{-6}$ |
| 1.90791528718078622313186060719755E30 | 10242.78033 | 10242.78831 | $7.78 \times 10^{-7}$ |
| 2.40997280874481941027683455628022E30 | −9268.18171 | −9268.20237 | $2.23 \times 10^{-6}$ |
| 3.96223170348329766133173296323307E30 | −7483.25069 | −7483.25629 | $7.49 \times 10^{-7}$ |
| 6.08302869527654580706324834685591E30 | 9729.73977 | 9729.76572 | $2.66 \times 10^{-6}$ |
| 6.22949918584108157280509258048941E30 | −14335.913 | −14335.948 | $2.39 \times 10^{-6}$ |
| 6.63237818782358897400245791070660E30 | 12010.61542 | 12010.63042 | $1.24 \times 10^{-6}$ |
| 6.82518053590400464945958005027695E30 | 12561.854 | 12561.880 | $2.02 \times 10^{-6}$ |
| 9.66225949824542125462912689965804E30 | −12070.174 | −12070.186 | $1.07 \times 10^{-6}$ |
| 9.83228440804649950062286954013174E30 | −10123.63767 | −10123.62451 | $1.30 \times 10^{-6}$ |
| 1.472969364244678383540127734056387E31 | −7707.03675 | −7706.96360 | $9.49 \times 10^{-6}$ |
| 1.970307186839852852378315796522770E31 | −13870.960 | −13870.969 | $6.84 \times 10^{-7}$ |
| 2.845670144939668837484722465519674E31 | −15484.372 | −15484.420 | $3.12 \times 10^{-6}$ |
| 3.557586000421470624922724880597724E31 | 13337.02398 | 13337.12691 | $7.71 \times 10^{-6}$ |
| 3.924676458989430915525116928410405E31 | 16244.39356 | 16244.47429 | $4.96 \times 10^{-6}$ |
| 4.589001484792927188496155886462819E31 | 9967.75442 | 9967.80244 | $4.80 \times 10^{-6}$ |
| 5.005475723107396211588045467161740E31 | −10621.48019 | −10621.51326 | $3.12 \times 10^{-6}$ |
| 6.080550187116167583688613601408264E31 | 14110.817 | 14110.835 | $1.28 \times 10^{-6}$ |
| 7.039106631049132430879196955445325E31 | −14055.05211 | −14055.11312 | $4.33 \times 10^{-6}$ |
| 9.066617388021913893282064021886284E31 | 15601.555 | 15601.620 | $4.16 \times 10^{-6}$ |
| 1.2783429043337574414365607070650632E32 | −15566.535 | −15566.607 | $4.64 \times 10^{-6}$ |
| 2.8609411594551963691214046991910957E32 | 16093.291 | 16093.351 | $3.73 \times 10^{-6}$ |
| 2.9615715957422120018598884000190662E32 | −16477.936 | −16478.028 | $5.61 \times 10^{-6}$ |
| 3.1067883362908396566754057659368205E32 | 16874.148 | 16874.202 | $3.20 \times 10^{-6}$ |
| 4.2114691908666607058133050702649423E32 | −15207.804 | −15207.882 | $5.14 \times 10^{-6}$ |
| 6.9619002836022067265947731886112914E32 | 11358.636 | | |
| 7.1071732114916294610536133167221954E32 | −8632.336 | | |
| 9.1724411211204221645707071038294101E32 | −11557.987 | | |
| 1.04460781211065043835431081519805492E33 | 11303.206 | | |
| 1.06269161656636941548485858301763355E33 | −7331.794 | | |
| 1.42368234667160116783366845609851195E33 | 11563.860 | | |
| 2.13961597512187919855165627131694812E33 | −6474.999 | | |
| 5.08596628180943122832041195042272642E33 | 4904.773 | | |

Table 5. Illustrative calculations of the Hardy function for $t > 10^{28}$ utilising the hybrid formulation methodology, compared to the corresponding values listed in [9, 23]. At these particular $t$ values the modulus of $Z(t)$ is unusually large. For all the hybrid computations, the precision setting was $\varepsilon_t = 0.005$ with $K \in [598, 992]$.